\documentclass[12pt]{article}

\usepackage{amsthm,amsmath,amsfonts,amssymb}
\usepackage{natbib} 
\usepackage[english]{babel}
\usepackage{csquotes}
\usepackage{mathtools}
\usepackage{bm}    
\usepackage{bbm}        %needed for bold ones
\usepackage{graphicx}   %pictures
\usepackage{tcolorbox}
\usepackage{hyphenat}
\usepackage{tikz}
\usepackage{authblk}
\usetikzlibrary{shapes, decorations, graphs, quotes, arrows.meta}

\usepackage[colorlinks,citecolor=blue,urlcolor=blue]{hyperref}
\hypersetup{urlcolor=blue, citecolor=blue, colorlinks=true, linkcolor=blue} 

\usepackage{lmodern}

\theoremstyle{plain}
    \newtheorem{theorem}{Theorem}[section]
    \newtheorem{lemma}[theorem]{Lemma}
    \newtheorem{proposition}{Proposition}
    \newtheorem{corollary}{Corollary}

\theoremstyle{remark}
    \newtheorem{definition}[theorem]{Definition}
    \newtheorem*{example}{Example}

\newcommand{\norm}[1]{\left\lVert #1 \right\rVert}
\newcommand{\round}[1]{\left( #1 \right)}
\renewcommand{\square}[1]{\left[ #1\right]}
\newcommand{\curly}[1]{\left\{ #1 \right\}}
\newcommand{\abs}[1]{\left| #1 \right|}

\newcommand{\Exp}[1]{\mathbb{E}\round{#1}}

\newcommand{\Cov}[2]{\mathbb{C}{\rm ov}\round{#1,#2}}
\newcommand{\Corr}[2]{\mathbb{C}{\rm or}\round{#1,#2}}

\newcommand{\N}[0]{\mathbb{N}}

\newcommand{\R}[0]{\mathbb{R}}
\newcommand{\indOne}[0]{\mathbbm{1}}
\newcommand{\ul}[1]{\underline{#1}}
\newcommand{\ol}[1]{\overline{#1}}	
\newcommand{\Mv}[3]{\mathcal{M}_{#1}\begin{pmatrix}
	#2\\
	#3
\end{pmatrix}}
\newcommand{\Mh}[3]{\mathcal{M}_{#1}\round{#2,#3}}
\newcommand{\onen}[0]{{\mathbf{1}_{n}}}
\newcommand{\onenp}[0]{{\mathbf{1}_{np}}}
\newcommand{\onenn}[0]{{\mathbf{1}_{n\times n}}}
\newcommand{\onep}[0]{{\mathbf{1}_{p}}}
\newcommand{\onepp}[0]{{\mathbf{1}_{p\times p}}}
\newcommand{\ie}{\textit{i.e. }}
\newcommand{\eg}{\textit{e.g. }}

\newcommand{\dVec}[1]{{\bm{#1}}}            %used for deterministic vectors
\newcommand{\latMatr}[1]{\mathbf{#1}}          %used for latin matrices
\newcommand{\grMatr}[1]{{\boldsymbol{#1}}}   %used for greek matrices
\newcommand{\matr}[1]{\mathbf{#1}}          %temp
\newcommand{\rVec}[1]{{\bm{#1}}}            %used for random vectors

\newcommand{\tobi}[1]{\textcolor{teal}{#1}}
\DeclareMathOperator{\diag}{diag}
\DeclareMathOperator{\Var}{\mathbb{V}ar}

\DeclareMathOperator{\rank}{rank}

\begin{document}
\title{Vector-Valued Gaussian Processes and their Kernels on a Class of Metric Graphs}
%\runtitle{Kernels on Metric Graphs}
%\thankstext{T1}{A sample of additional note to the title.}

\author[A]{Tobia Filosi}%\ead[label=e1]{tobia.filosi@unitn.it}
\author[B]{Emilio Porcu*}%\ead[label=e2]{emilio.porcu@ku.ac.ae}
\author[C]{Xavier Emery}%\ead[label=e3]{xemery@ing.uchile.cl}
\author[A]{Claudio Agostinelli}%\ead[label=e4]{claudio.agostinelli@unitn.it}
\author[D]{Alfredo Alegría}%\ead[label=e5]{alfredo.alegria@usm.cl}
%%%%%%%%%%%%%%%%%%%%%%%%%%%%%%%%%%%%%%%%%%%%%%
%% Addresses                                %%
%%%%%%%%%%%%%%%%%%%%%%%%%%%%%%%%%%%%%%%%%%%%%%
\affil[A]{Department of Mathematics, University of Trento}%\printead{e1,e4}}
\affil[B]{Department of Mathematics, Khalifa University $\&$ \textsc{ADIA LAB}}%\printead[]{e2}}
\affil[C]{Department of Mining Engineering, Universidad de Chile }%\printead[]{e3}}
\affil[D]{Department of Mathematics, Universidad T\'ecnica Federico Santa Mar\'ia}% \printead[presep={,\ }]{e5}}

\maketitle
% \begin{keyword}[class=MSC]
%     \kwd[Primary ]{60G60, 05C12, 30L05; }
%     \kwd[secondary ]{20H20, 05C22}
% \end{keyword}

% \begin{keyword}
%     \kwd{Embeddings}
%     \kwd{Matrix-valued kernels}
%     \kwd{Metric Graphs}
%     \kwd{Positive semidefiniteness}
%     \kwd{Quasi-Laplacian}
% \end{keyword}
\begin{abstract}
    Despite the increasing importance of stochastic processes on linear networks and graphs, 
    current literature on multivariate (vector-valued) Gaussian random fields on metric graphs is elusive. This paper challenges several aspects related to the construction of proper matrix-valued kernels structures. We start by considering matrix-valued metrics that can be composed with scalar- or matrix-valued functions to implement valid kernels associated with vector-valued Gaussian fields. We then provide conditions for certain classes of matrix-valued functions to be composed with the univariate resistance metric and ensure positive semidefiniteness. Special attention is then devoted to Euclidean trees, where a substantial effort is required given the absence of literature related to multivariate kernels depending on the $\ell_1$ metric. Hence, we provide a foundational contribution to certain classes of matrix-valued positive semidefinite functions depending on the $\ell_1$ metric. This fact is then used to characterise kernels on Euclidean trees with a finite number of leaves. Amongst those, we provide classes of matrix-valued covariance functions that are compactly supported.
\end{abstract}

\section{Introduction}
\subsection{Context: Random Fields on Networks}
    Networks analysis has become ubiquitous in several branches of theoretical and applied sciences, including mathematics, statistics, machine learning, artificial intelligence, with applications connected to the data planet. Networks can be used, either, to represent connections between random variables, or to provide a topological structure where a given process (called a random field) is observed. Our paper deals with the second situation, under the condition that the observed process is a vector-valued Gaussian random field that is second-order isotropic in the sense that the matrix-valued covariance kernel depends on distance between the points in the network. \par
    Random fields can be continuously \citep{anderes_isotropic_2020} or discretely indexed over a network; for the second case, the common nomenclature is that of point processes on networks \citep{moradi, BADDELEY2021100435, BADDELEY}. For both cases, defining a random field over a network is a challenging task.\par
    The electrical engineering and machine learning communities have been very active on the subject of networks, and the reader is referred to {\cite{Ortega2018} and} \cite{chami2022machine} for a comprehensive review. \par
    For this manuscript, the network defines a topological structure representing the domain of definition of a given random field. This fact is exploited by \cite{anderes_isotropic_2020}: endowing the network with a metric provides a metric space that is then (quasi) embeddable over Hilbert spaces. This is the crux of the argument that allows to build isotropic covariance functions over networks, where isotropy is understood with respect to two alternative metrics (the geodesic and the resistance metrics). A relevant remark is that  networks    endowed with metrics can be {suitably} represented through metric graphs. This approach is adopted by \cite{bolin}, \cite{anderes_isotropic_2020} and \cite{bolin2022gaussian}. We shall make use of this fact subsequently. Scalar {random fields} that evolve temporally over graphs have been challenged by \cite{porcu2023stationary}, \cite{filosi2023temporally} and \cite{tang}. \par    
    The way a metric graph can be specified is certainly not unique. \cite{anderes_isotropic_2020} work over graphs with Euclidean edges --- called generalised networks in \cite{porcu2023stationary} --- which extend linear networks to non-linear edges. Further, the {random field} defined over such structures can have realisations over any point over the edges, and not only in the nodes. Roughly, these are graphs where each edge is associated with an abstract set in bijective correspondence with a segment of the real line. This provides each edge with a Cartesian coordinate system to measure distances between any two points on that edge. \par
    The fact that linear networks are overly limited to provide a suitable topological structure to {random fields} is well understood in the machine learning literature  \citep{alsheikh2014machine,  hamilton2017representation,  borovitskiy2022isotropic}, as well as in spatial statistics \citep{cressie2006spatial,  ver2006spatial, peterson2013modelling, peterson2007geostatistical, montembeault2012impact, xiao2017modeling, perry2013point, deng2014ginibre, baddeley2017stationary}.

\subsection{Challenge}
    Our paper is unique according to the state of the art. The literature on vector-valued random fields on (linear or generalised) networks is elusive. This is not surprising: kernels on metric graphs are a very recent subject, and so far the main efforts have been focused to the case of scalar-valued random fields. \par
    A subset of the graphs used in this paper is represented by the so-called Euclidean trees with a given number of leaves, which are substantially linear networks. For such a case, since a linear network is embedded into the two-dimensional Euclidean space, any matrix-valued covariance functions depending on the $\ell_1$ norm (the so-called Manhattan distance) can be used as a covariance function for a vector-valued field on a linear network. The reader is referred to \cite{zastavnyi_positive_2000} for relevant results in terms of isometric embeddings. Unfortunately, there are two problems related to this fact: 
    \begin{enumerate}
        \item To the knowledge of the authors, there are no models for matrix-valued covariances depending on the $\ell_1$ distance;
        \item Such constructions are no longer valid for the case of generalised networks.
    \end{enumerate}
    Another relevant aspect is related to what we term {\em cross metrics} in the paper. To illustrate the principle, consider two {random fields} --- call them $Z_1$ and $Z_2$. Since the network does not necessarily represent a geographic space, there is no reason to believe that the distance between $Z_1(u)$ and $Z_2(v)$, for two different points $u,v$ {on the network}, should be the same as the distance between $Z_1(u)$ and $Z_1(v)$ for instance. We take this aspect into account and devote attention to matrix-valued metrics.

\subsection{Our Constructions: Road Map}
    Our constructions are based on the following steps. We provide here a concise description. Careful notation and terminology are established in Section \ref{sec:mathBackground}.\par
    \noindent {\bf The Road Map}
    \begin{enumerate}
        \item The purpose of the work is to construct a $p$-variate Gaussian {random field}, denoted $\rVec{Z}:= \{ \rVec{Z}(u)=(Z_1(u),\ldots,Z_p(u))^{\top}, \; u \in {\cal G} \}$ throughout and defined over a metric graph, $\mathcal G$, with a given matrix-valued covariance function that is isotropic, that is it depends on the {\em distance} between the points of the graph. Hence, the topological space $\mathcal G$, endowed with a proper matrix-valued (semi) metric $\latMatr D: \mathcal G \times \mathcal G \to [0,+\infty)^{p \times p}$, becomes a (semi) metric space $(\mathcal G, \latMatr D)$. 
        \item Such a (semi) metric space can be (quasi) embedded on a suitable Hilbert space, which in turns allows to resort classical machinery for covariance functions on metric spaces \citep{schoenberg}.
        \item Hence, a substantial part of the job is to build the multivariate metric $\latMatr D$.
        \item We then take advantage of isometric embedding arguments as in \cite{anderes_isotropic_2020} to claim that, for a certain class of continuous mappings $\psi: [0,+\infty) \to (0,+\infty)$, the {element-wise} composition
        \begin{equation} 
            \label{eq:covariance_general}
            \latMatr K(u,v):= \psi \left ( \latMatr D(u,v) \right ), \qquad u,v \in \mathcal G, 
        \end{equation}
        provides a valid matrix-valued covariance function.
        \item We then generalise this construction by considering a matrix-valued class of continuous mappings $\grMatr{\Psi}:{[0,+\infty)}^{p \times p  } \to {[0,+\infty)}^{p \times p}$ having elements $\Psi_{ij}$, such that the {element-wise} composition 
        \begin{equation} 
            \label{eq:covariance_general_2}
            \latMatr{K}(u,v):= \grMatr{\Psi} \left ( \latMatr{D}(u,v) \right ) = \left  [ \Psi_{ij}(D_{ij}(u,v)) \right ]_{i,j=1}^p, \qquad u,v \in \mathcal G, 
        \end{equation}
        is a valid matrix-valued covariance function.
    \end{enumerate} 
    In turn, the construction of the multivariate metric $\latMatr D$ requires a collection of steps, which entails several technical challenges: 
    \begin{description}
        \item[(a)] Following the construction in \cite{anderes_isotropic_2020}, the vector-valued random field $\bm{Z}$ is defined as the sum of two random fields, $\bm{Z}_V$ and $\bm{Z}_E$, where the subscripts $V$ and $E$ stand, respectively, for the sets containing vertices and the edges associated with the metric graph, $\mathcal G$. Both {random fields} are continuously indexed over the graph;
        \item[(b)] The {random field} $\bm{Z}_V$ at any point of an edge, $e \in E$, is constructed through interpolation of the {random field} defined at the extremes of each edge (which are obviously vertices in $V$);
        \item[(c)] The multivariate metric is then checked for desiderata. Specifically, that the marginal metrics (those in the diagonal of $\latMatr{D}$) respect the ingenious univariate condition provided by \cite{anderes_isotropic_2020}. Further, the cross-elements in $\latMatr{D}$ are required to be homogeneous, in the sense that $D_{ij} = D_{i'j'}$ for $i \ne j $ and $i' \ne j'$, for $D_{ij}$ a generic element of $\latMatr{D}$.
    \end{description}
    For a subclass termed {\em Euclidean trees}, we provide multivariate kernels with compact support depending on the so-called Manhattan distance. \par 
    The structure of the paper is the following. Section \ref{sec:mathBackground} provides the necessary background and notation. Results are provided in Sections \ref{sec:results} and \ref{sec:linearNetworks}. Section \ref{sec:conclusions} concludes the paper. As proofs are rather technical and lengthy, they are deferred to Appendix \ref{A_proofs}. Further, Appendix \ref{A_B} contains a worked example that show some relevant property of the proposed construction. \par
\begin{figure}[ht]
        \centering
        \begin{tikzpicture}[thick, every node/.style = {draw, ellipse}, scale=0.8] 
            \node (L_eigQ) at (3,6) {L\ref{lem:eigenvaluesQ}};
            \node (L_useId) at (6,6) {L\ref{lem:usefulIdentities}};
            \node (P_exprK) at (9,6) {P\ref{prop:explicit_K_ZV_and_K_ZE}};
            \node (T_Zast) at (12,6) {T\ref{theo:MultivariateSchoenberg}};
            \node (T_covFB) at (12,4) {T\ref{theo:covFunctionBuilding}};
            \node (P_qL_prop) at (0,4) {P\ref{prop:quasiLaplacianProperties}};
            \node (L_XplusY) at (1, 6) {L\ref{lem:sumOfMatricesPD}};
            \node (L_ThetaInv) at (4,4) {L\ref{lem:KinverseTheta}};
            \node (P_ourDprop2) at (8,4) {P\ref{prop:ourDistProperties}};
            \node (T_gen1) at (12,2) {T\ref{theo:la_raja}};
            \node (P_PDequi) at (0, 2) {P\ref{prop:equivalentFormPositDef}};
            \node (P_ThetaQL) at (2,4) {P\ref{prop:ThetaQuasiLaplacian}};
            \node (P_ourDprop1) at (4,2) {P\ref{prop:actuallyQuasiMetricSpace}};
            \node (P_classPhi) at (0,0) {P\ref{prop:p-Schoen}};
            \node (T_gen2) at (9,0) {T\ref{theo:el_patron}};
            \node (T_gen3) at (11,0) {T\ref{theo:dios}};
            \node (T_gen4) at (13,0) {T\ref{theo:dios2}};
            \node (T_gen5) at (15,0) {T\ref{theo:dios3}};
            \node (P_indPD) at (2,2) {P\ref{prop:covZ_V_isACov}};
            \node (P_classPhiIT) at (2,0) {P\ref{prop:PleaseEmilioChooseWiseLabels}};
            \node (P_genKisPD) at (4,0) {P\ref{prop:infarto_miocardio}};
            \node (P_asympt) at (6, 4) {P\ref{prop:asymptoticResults}};
            \node (C_exprKZ) at (10, 4) {C\ref{cor:explicitExprK_Z}};
            
            \draw [->] 
                (L_eigQ) edge (L_useId)
                (L_eigQ) edge (L_ThetaInv)
                (L_eigQ) edge (P_ThetaQL)
                (L_useId) edge (L_ThetaInv)
                (L_useId) edge (P_asympt)
                (P_qL_prop) edge (P_ThetaQL)
                (L_XplusY) edge (P_ThetaQL)
                (L_ThetaInv) edge (P_ThetaQL)
                (L_ThetaInv) edge (P_ourDprop1)
                (L_ThetaInv) edge (P_indPD)
                (P_ThetaQL) edge (P_ourDprop1)
                (P_ThetaQL) edge (P_indPD)
                (P_PDequi) edge (P_indPD)
                (P_exprK) edge (P_ourDprop2)
                (T_Zast) edge (T_covFB)
                (T_covFB) edge (T_gen1)
                (T_gen1) edge (T_gen2)
                (T_gen1) edge (T_gen3)
                (T_gen1) edge (T_gen4)
                (T_gen1) edge (T_gen5)
                (P_classPhi) edge (P_classPhiIT)
                (P_classPhiIT) edge (P_genKisPD)
                (P_exprK) edge (C_exprKZ);
        \end{tikzpicture} 
        \caption{Road-map of the main results in this manuscript. L, P, T and C stand for Lemma, Proposition, Theorem and Corollary, respectively. The suffix  indicates that the result is stated in Appendix \ref{A_proofs}. The missing results have either no direct connection with the above results (L\ref{lem:tree}, P\ref{prop:explicitWritingMultivMetric}) or are {essential for many others} (L\ref{lem:MpMp_expr}).} 
        \label{fig:resultsRoadMap}
    \end{figure}
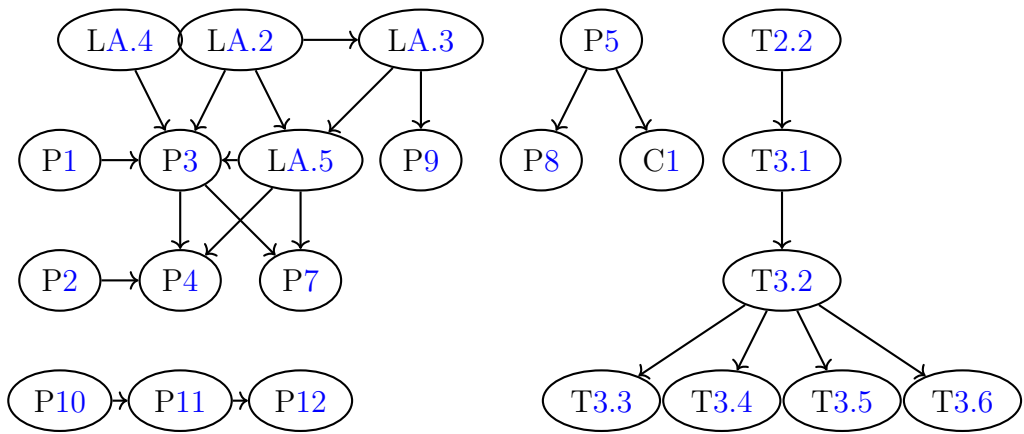

    Finally, to guide the reader into the structure of the results of this manuscript, we provide a road-map in Figure \ref{fig:resultsRoadMap}.
\section{Mathematical Background} 
\label{sec:mathBackground}
%\emi{We need to add linear networks and linear trees. We also need to add a couple of pictures to represent the graphs. } \tobi{I can take care of the pictures. Is it ok to recycle the one in the paper of Alfredo?} \emi{Yes it is okay. Please add a definition of Euclidean tree - just half a line. }
\subsection{Notation handout}
    \label{ssec:notation}
    %We start by explicitly stating the relevant notation. 
    Throughout, $p$ is an integer greater than $1$ (the case $p=1$ has been treated by \cite{anderes_isotropic_2020}).
    Bold lowercase letters denote (deterministic) vectors (\eg $\dVec{\delta}$). Bold uppercase letters denote (deterministic) matrices (\eg $\latMatr{A}$). Bold uppercase italic letters denote random vectors (\eg $\rVec{Z}$). $\latMatr{A}^\top$, $\dVec{\delta}^\top$ and $\rVec{Z}^\top$ are the transposes of $\latMatr{A}$, $\dVec{\delta}$ and $\rVec{Z}$, while $\latMatr{A}^-$ stands for the Moore-Penrose inverse of $\latMatr{A}$. Non-bold letters represent scalar (real-valued) quantities or non-numerical objects (sets, graphs, vertices, edges, or points). The $(i,j)$-th entry of $\latMatr{A}$ is denoted ${A}_{ij}$, while the $i$-th components of $\dVec{\delta}$ and $\rVec{Z}$ are ${\delta}_i$ and ${Z}_i$, respectively. The all-ones vector of size $n$, all-ones matrix of size $n \times n$, {and identity matrix of size $n\times n$} are denoted as $\onen$, $\onenn$ and $\latMatr{I}_n$, respectively. Finally, we define $\latMatr{J}_n := \diag{[   0, 0, \dots, 0, 1]^\top}$ ($n \times n$ diagonal matrix with only one non-zero entry).
    
    Continuity, differentiation, integration and scalar- or matrix-valued functions (except power functions) involving matrices are understood as being performed element-wise. Concerning the power functions, $\latMatr{A}^2$ is $\latMatr{A}$ multiplied by itself, $\sqrt{\latMatr{A}}=\latMatr{A}^{1/2}$ is the principal square root of $\latMatr{A}$, and $\latMatr{A}^{-1}$ is the inverse of $\latMatr{A}$ in the matrix sense. Products (\eg $\latMatr{A} \latMatr{B}$ or $\latMatr{A} \dVec{\delta}$) and Kronecker products (\eg $\latMatr{A} \otimes \latMatr{B}$) are also understood in the matrix sense, not as element-wise operations.

\subsection{Gaussian {random fields} on graphs with Euclidean edges}
    \label{ssec:graphsWithEE}

    The paper deals with a vector-valued Gaussian random field ${\rVec{Z}}$ %$= \round{Z_{1}, \ldots, Z_p}^{\top}$ 
     defined on a metric graph as exposed below. We assume a {random field} %weakly stationary \xave{IS IT NECESSARY TO STATE STATIONARITY, OR IS IT SUFFICIENT TO STATE ISOTROPY AS DONE NEXT?} \emi{We do not cover fields with infinite variance.} \xave{This is excluded next because we state (one line below) 
    %that there is an isotropic covariance, so the variance exists}, 
    with zero-mean vector and with a matrix-valued covariance matrix $\latMatr K: \mathcal G \times \mathcal G \to \R^{p \times p}$ that is isotropic in the sense there exists a pair $(\grMatr \Psi, \latMatr D)$, with $\grMatr \Psi: [0,+\infty)^{p \times p} \to \R^{p \times p}$ and $\latMatr D : \mathcal G \times \mathcal G \to {[0,+\infty)}^{p \times p}$, both continuous, such that Equation (\ref{eq:covariance_general_2}) is satisfied. 
    We work with a specific class of metric graphs, termed graphs 
     {\em with Euclidean edges} as being proposed by \cite{anderes_isotropic_2020} to generalise linear networks. The main idea underlying this construction is to associate every edge of the graph to a segment of the real line. As a consequence, it is possible to consider both the vertices and the points over the edges as actual points of the graph, where distances can be computed and {random fields} can be defined. As this is natural on linear networks as well, we stress that graphs with Euclidean edges do not have any restriction on the topology, \ie the edges' lengths are free to vary and so are the connections between the vertices. Here we rephrase their original definition.
   % \emi{Tobia, ricicla la definizione di grafo dal paper con Alegria-Emery. Grazie. } \tobi{It should be ok now.}

    \begin{definition}[Graph with Euclidean edges]
        \label{def:graphWithEE}
        The type of topology considered by \cite{anderes_isotropic_2020} (to which the reader is referred for additional details and motivation) is called graph {\em with Euclidean edges}, denoted with a triple ${\cal G}=\left( V,E,\{\varphi_e\}_{e\in E} \right)$ throughout, where the elements are blended in the following way:
        \begin{enumerate}
            \item[(a)] {{$(V, E)$}} has a graph structure, \ie $V$ is the set of vertices and {$E\subset V\times V$} accounts for the edges. We assume that this graph is simple and connected, \ie $V$ is finite, the graph has not repeated edges or edges that join a vertex to itself and every pair of vertices is connected by a path.
            \item[(b)] Each edge $e\in E$ is provided with a \emph{length} $\ell(e)>0$ and a \emph{weight} $w(e):=1/\ell(e)$; furthermore, it is associated with a unique abstract set, also denoted $e$, such that $V$ and all the edge sets are mutually disjoint. 
            \item[(c)] Let  $v_1$ and $v_2$ be vertices connected by $e\in E$. Then, $\varphi_e$ is a continuous and bijective mapping defined on $e \cup \{v_1,v_2\}$, such that $\varphi_e$ maps $e$ onto an open interval $(0,\ell(e)) \subset \R$ and $\{v_1,v_2\}$ onto $\{0,\ell(e)\}$.
            %\item[(d)] Let $d_\mathcal{G}:\mathcal{V}\times \mathcal{V}\rightarrow [0,+\infty)$ be the shortest path weighted graph metric on the vertices of $(\mathcal{V},\mathcal{E})$, with edge weights given by $\overline{e}-\underline{e}$ for every $e\in\mathcal{E}$.  For each $e\in\mathcal{E}$ connecting $u,v\in\mathcal{V}$, we have $d_\mathcal{G}(u,v) = \overline{e}-\underline{e}$. \color{red}{Tobia: Here I think that we could drop this condition. Anderes et al. use it just as a bridge from the geodesic distance and the resistance one, yet I think we do not need geodesic distance. However, if we drop this, we should adjust some notation thoughtout the paper.} \color{brown}{EP: I agree with Tobia. If you only use the resistance metric, you do not need to adopt this technical condition.}
        \end{enumerate}
        % \todo{this is a copypaste}
        % Consider a simple, connected and weighted graph $G=(V, E, w)$, where $w:E\to \R^+$ represents the weight mapping. Then, $G$ is called a \emph{graph with Euclidean edges} provided that the following conditions hold.
        % \begin{enumerate}
        %     \item Edge sets: Each edge $e\in E$ is associated to the compact segment (also denoted by $e$) $[0,\ell(e)]$, where $\ell(e):=w(e)^{-1}$ may be interpreted as the \emph{length} of the edge $e$. 
        %     \item Linear edge coordinates: Each point $u\in e=(\ul u,\ol u)$ is uniquely determined by the endpoints $\ul u$ and $\ol u$ of $e$ and its relative distance $\delta(u):=\frac{u}{\ell(e)}=u\,w(e)$ from $\ul u$, that is $u=\round{\ul u, \ol u, \delta(u)}$, so that $\ul u=(\ul u,\ol u,0)=(\ol u,\ul u,1)$ and $\ol u=(\ul u,\ol u,1)=(\ol u,\ul u,0)$.
        % \end{enumerate}
    \end{definition}
    Henceforth, we shall assume the existence of a total order relation on the set of vertices $V$ and that every edge is represented through the ordered pair $(v_1,v_2)$, where $v_1<v_2$. {By abuse of notation, we write} {$u\in\mathcal G$ to denote a point on the graph: it can be either a vertex or a point on an edge. Each point $u\in \mathcal G$ will be identified with the triple $(\ul u, \ol u, \delta)$, where: $\ul u < \ol u$ are the endpoints of the edge $e$ containing $u$ and $\delta \in [0,1]$ is the relative distance of $u$ from $\ul u$; formally: $\delta=\delta(u):=w(e)\varphi_e(u)$. Notice that, whenever $u$ is a vertex, it is always possible to write $u=(\ul u, \ol u, \delta)$ for $\delta \in \curly{0,1}$, though the choice of $\ul u$ and $\ol u$ may not be unique.} Finally, we use the notation $e(u)$, where $e:\mathcal G \to E$, to denote the edge that contains $u$. If $u\in V$, then $e(u)$ is any of the edges that have an endpoint in $u$.  It is in order to note that our definition slightly deviates from the original one in \cite{anderes_isotropic_2020}: we do not require any distance consistency property, as we {do not} use {the} geodesic metric.\par
\begin{figure}
        \centering
        \includegraphics[width = \textwidth]{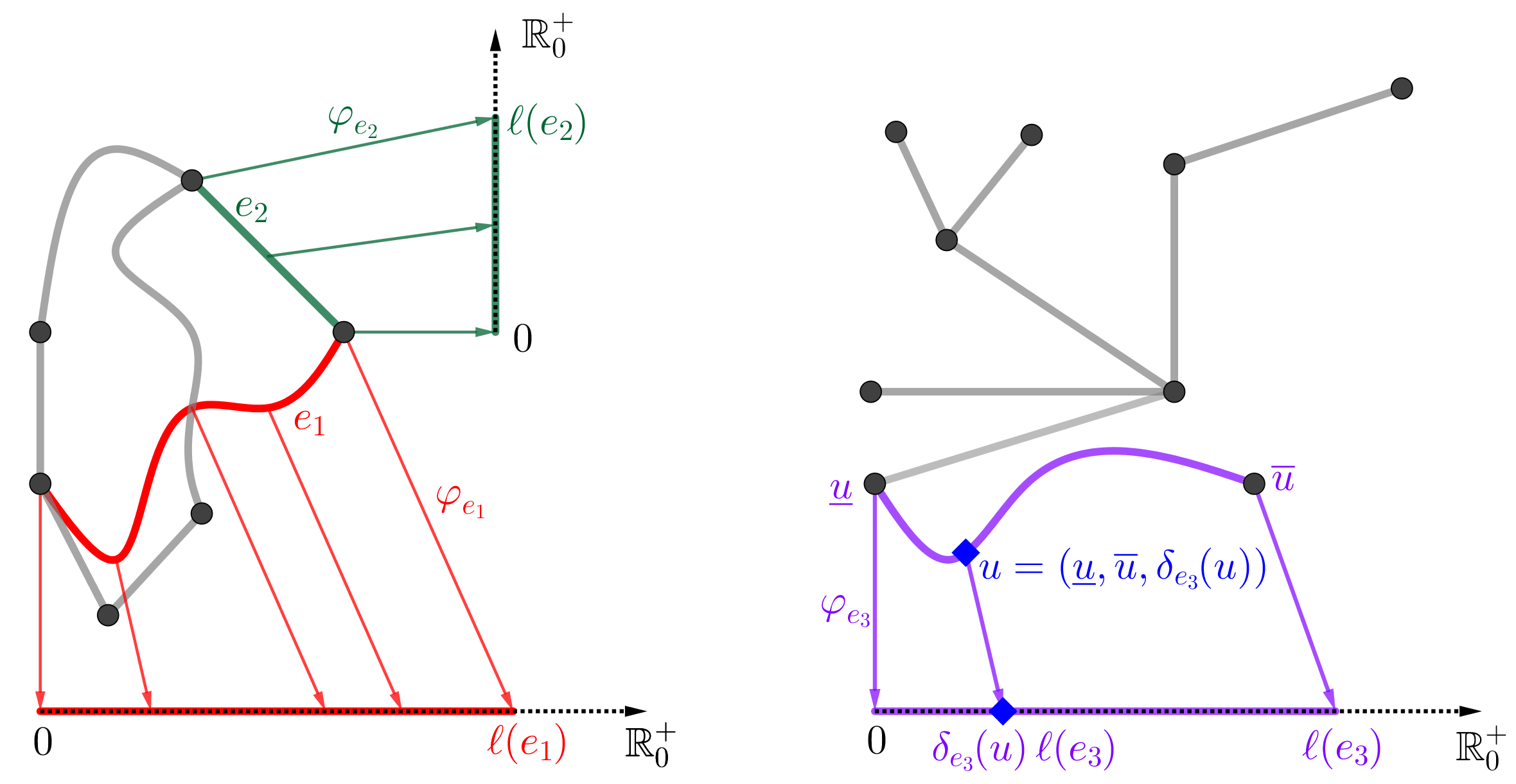}
        \caption{Left: a graph with Euclidean edges, where the bijections $\varphi_{e_1}$ and $\varphi_{e_2}$ have been highlighted. Right: a Euclidean tree {with $5$ leaves} where the role of $\ul u$, $\ol u$ and $\delta_{e_3}(u)$ have been stressed. Adaptation of \citet[Figure 1]{filosi2023temporally}.}
        \label{fig:graphEE_eucTree}
    \end{figure}

 
    %\emi{Tobia: here should go the definition of Euclidean tree.} \tobi{Please someone check the definition and the figure.} \xave{OK for me}
    A graph with Euclidean edges is called an \emph{Euclidean tree} if it has a tree structure, \ie given any two points $u_1,u_2 \in \mathcal{G}$, there is exactly one path that connects $u_1$ to $u_2$. {A vertex of a Euclidean tree is called a \emph{leaf} if it is connected to only one other vertex.} Figure \ref{fig:graphEE_eucTree} aims to clarify the notions introduced so far.
    
\subsection{Classes of matrices used in the paper}
    An $n\times n$ real matrix $\latMatr A$ is said to be \emph{positive semidefinite} if and only if $\dVec{c}^\top \latMatr{A} \dVec{c} \geq 0$ for all $\dVec{c} \in \mathbb{R}^n$, and \emph{conditionally negative semidefinite} if and only if  $\dVec{c}^\top \latMatr{A} \dVec{c} \leq 0$ for all $\dVec{c} \in \mathbb{R}^n$ such that $\onen^\top \dVec{c} = 0$.

    An $n\times n$ matrix $\latMatr L$ is called \emph{quasi-Laplacian} if and only if it is symmetric, positive semidefinite, and has exactly one null eigenvalue, with corresponding eigenvector $\onen$. A quasi-Laplacian matrix $\latMatr L$ is called \emph{Laplacian} if and only if it has non-positive off-diagonal entries.
    Notice that Laplacian matrices arise naturally from weighted graphs: there is a bijection between the set of (positively) weighted, simple and undirected graphs and the set of Laplacian matrices \cite[Subsection 2.1]{devriendt_effective_2022}. Henceforth, we will consider the Laplacian matrix of a given graph with Euclidean edges, with entries defined as \citet[Equation (1)]{devriendt_effective_2022}: 
    {\begin{equation}
        \label{eq:laplacianofagraph}
        L_{ij} = 
        \begin{dcases}
            \sum_{v\in V} w((v, v_i)) &\text{ if $i=j$}\\
            -w((v_i,v_j))  &\text{ if $i \neq j$,}
        \end{dcases}
    \end{equation}
    with $w((v_i,v_j))$ the weight of the edge joining vertices $v_i$ and $v_j$ if these vertices are connected, $0$ otherwise.} \par %\tobi{Beware that the elements of the laplacian are formed with the weight (\ie the inverse of the length): I fixed this.}  \par
    Quasi-Laplacian matrices have been introduced as a generalisation of Laplacian matrices, as they will often appear in this manuscript. Their interpretation is slightly blurred, since, thinking about the bijection between Laplacian matrices and weighted graphs, it seems natural to interpret quasi-Laplacian matrices as Laplacian matrices of graphs with possibly negative weights. Nevertheless, as one may expect, quasi-Laplacian matrices share several properties of Laplacian matrices formally stated below.
    \begin{proposition}
        \label{prop:quasiLaplacianProperties}
        Let $\latMatr L$ be a quasi-Laplacian matrix partitioned as follows: 
        \begin{equation}
            \label{eq:partitionQuasiLapl}
            {\latMatr{L}}=\begin{bmatrix}
                {\latMatr{A}} & {\latMatr{B}}\\
                {\latMatr{B}}^\top & {\latMatr{C}}
            \end{bmatrix},
        \end{equation}
        where ${\latMatr{A}}$ and ${\latMatr{C}}$ are square matrices. Then the following statements hold.
        \begin{enumerate}
            \item ${\latMatr{B}}\ne \grMatr 0$.
            \item ${\latMatr{A}}$ and ${\latMatr{C}}$ are strictly positive definite.
            \item Both ${\latMatr{L}}/{\latMatr{C}}:={\latMatr{A}}-{\latMatr{B}} {\latMatr{C}}^{-1} {\latMatr{B}}^\top$ and ${\latMatr{L}}/{\latMatr{A}}:={\latMatr{C}} - {\latMatr{B}}^\top {\latMatr{A}}^{-1} {\latMatr{B}}$ are quasi-Laplacian matrices, \ie, the set of quasi-Laplacian matrices is closed under the Schur complement operation.
            \item ${\latMatr{L}}^-$ is a quasi-Laplacian matrix.
        \end{enumerate}
    \end{proposition}

\subsection{Metrics used in the paper}
    The work by \cite{anderes_isotropic_2020} shows that there is no unique way to define a metric over a graph with Euclidean edges. While the geodesic distance is physically intuitive, \cite{anderes_isotropic_2020} prove that it has very limited use in terms of available covariance functions. Further, when considering the geodesic distance not all graphs with Euclidean edges become {\em permissible} { \citep{anderes_isotropic_2020}}, and a collection of technical restrictions is required. Alternatively, one can use the resistance metric, being a generalisation of the electric distance used for electric circuits \citep{Klein}. To provide a description of the resistance metric, some further background is needed.
    Let $X$ be a set and $Z:X\to \R$ a square-integrable random {field}. The \emph{variogram} of $Z$ is defined via
    \begin{align}
        \gamma_{{Z}}: X\times X &\to \R \nonumber\\
        (x_1,x_2)&\mapsto \gamma_{{Z}} (x_1,x_2):=\Var\round{{Z}(x_1)-{Z}(x_2)}.\nonumber
    \end{align}
    Clearly, a variogram is always symmetric and non-negative valued. In addition, it is conditionally negative semidefinite, \ie it satisfies \citep{chiles2012geostatistics}
    \begin{equation}
        \sum_{i,j=1}^n c_i c_j \gamma(x_i,x_j) \leq 0
    \end{equation}
    whenever $\sum_{i=1}^n c_i = 0$, for any $x_1,\dots,x_n\in X$. \par 
    \cite{anderes_isotropic_2020} define the resistance metric, denoting $d_R$ throughout, as the variogram of a  {random field} $Z$ on the graph with Euclidean edges $\mathcal G$, that is constructed {\em ad hoc} in order to be as much as possibly parenthetical to a Brownian bridge on the graph. Technicalities are deferred to their paper. \par
    The concept of variogram can be generalised to the case of a zero-mean square-integrable vector-valued  random field $\rVec{Z}$ on a non-empty set $X$, through the following lines. 
    %\emi{Let $X$ be a set and let $\rVec{Z}:X\to \R^p$ be a zero-mean square-integrable random  field. THIS STATEMENT IS REDUNDANT and NEEDS TO GO OUT.} 
    The \emph{pseudo-variogram} of $\rVec{Z}$ is defined through the identity
    \begin{align}
        \label{eq:pseudoVariogramDef}
        \grMatr \Gamma_{\rVec{Z}}:X\times X &\to {[0,+\infty)}^{p\times p}\nonumber \\
        (x_1,x_2) &\mapsto \grMatr \Gamma_{\rVec{Z}}(x_1,x_2):=\Big [ \Var\round{Z_i(x_1)-Z_j(x_2)} \Big ]_{i,j=1}^p. 
    \end{align} 
    Straightforward calculations allow to rewrite the above as
    \begin{equation}
        \label{eq:pseudoVariogramMatricial}
        \grMatr \Gamma_{\boldsymbol{Z}}(x_1,x_2) = \round{\diag \Var \boldsymbol{Z}(x_1)} \,\onep^\top + \onep\, \round{\diag \Var \boldsymbol{Z}(x_2)}^\top - 2\, \Cov{\boldsymbol{Z}(x_1)}{\boldsymbol{Z}(x_2)},
    \end{equation}
    where $\Var \boldsymbol{Z}(x_1)$ is the $p\times p$  {collocated} variance-covariance matrix of $\boldsymbol{Z}(x_1)$, $\diag \Var \boldsymbol{Z}(x_1)$ is the $p\times 1$ vector containing the main diagonal of $\Var \boldsymbol{Z}(x_1)$ and $\Cov{\boldsymbol{Z}(x_1)}{\boldsymbol{Z}(x_2)}:=\Exp{\boldsymbol{Z}(x_1)\boldsymbol{Z}(x_2)^\top}$.
    The pseudo-variogram satisfies the following properties \citep{dorr2023characterization}:
    \begin{itemize}
        \item {for any $x_1,x_2 \in X$,} 
        \begin{equation}
            \label{sympseudo}
            \grMatr \Gamma_{{\rVec{Z}}}(x_1,x_2)=\grMatr \Gamma_{{\rVec{Z}}}(x_2,x_1)^\top;
        \end{equation}
        \item {for any $x \in X$, the diagonal entries of $\grMatr \Gamma_{{\rVec{Z}}}(x,x)$ are zero;}
        \item $\grMatr \Gamma_{{\rVec{Z}}}$ is $\onep$-conditionally negative semidefinite, \ie for all {$n \in \mathbb{N}^+$, $x_1,\dots,x_n\in X$ and } $\dVec c_1,\dots,\dVec c_n\in\R^p$ such that $ {\onep}^\top \sum_{i=1}^n \dVec c_i = 0$, it holds
        \begin{equation}
            \label{eq:negativesemidef}
            \sum_{i,j=1}^n \dVec c_i^\top \grMatr \Gamma_{{\rVec{Z}}}(x_i,x_j) \dVec c_j \leq 0.
        \end{equation}
    \end{itemize}
    
    In the following, we define %\emi{\emph{define}: do not need to emphasize. You actually define it.} 
    the multivariate distance $\latMatr D:X\times X \to {[0,+\infty)}^{p\times p}$ via $\latMatr D(x_1,x_2):=\grMatr{\Gamma}_{\latMatr L}(x_1,x_2)$. It is in order to notice that at this point the matrix-valued distance $\latMatr D$ is \emph{not} symmetric, \ie in general $D_{ij}(x_1,x_2)\ne D_{ji}(x_1,x_2)$. This may pave the way for new asymmetric metrics. %\todo{comment on this}  \emi{I do not understand this comment.} \tobi{Please Alfredo add some comments on this and/or rephrase this part: you were the one stressing asymmetry properties. :)}
\subsection{Classes of functions and Schoenberg characterisations}
    \label{ssec:positiveDefiniteFunctions}
    This section follows closely \cite{zastavnyi_analog_2023}. A %symmetric \xave{\textcolor{red}{¡} HERE WE STATE A SYMMETRY CONDITION \textcolor{red}{!} } \emi{Yes, but this is for a scalar-valued function!} \tobi{I would remove this condition: there is no need for the function $k$ to be symmetric. In such a case, we solve all our problems, don't we?} \xave{I had problems in past papers because the symmetry condition was not explicit... Here we only use symmetric functions and matrix-valued functions (considering transpose and permutation of $x$ and $y$ at the same time), so I think the problem would be solved if we do as above and say: ``a symmetric function is positive semidefinite if and only if...''; same for matrix-valued functions: ``a matrix-valued function such that $K(x,y)=K(y,x)^\top$ is positive semidefinite if and only if...''} 
    function $k:X\times X \to \R$ is positive semidefinite if, for all $n\in\curly{1,2,\dots}$, ${x}_1,\dots,{x}_n \in X$ and $c_1,\dots,c_n \in \R$, one has 
    %\emi{let us choose one unique notation: $x$ or $u$. } \tobi{I agree on using only $x$ for generic sets and $u$ for graphs EE. In addition, let us use $X$ for a generic set, while use $\cal G$ for the graphs EE (As Anderes).}
    \begin{equation}
        \label{eq:positSemiDefFunction}
        \sum_{i,j=1}^n c_i c_j k({x}_i, {x}_j) \geq 0.
    \end{equation}
    In addition, $k$ is a strictly positive definite function if it is positive semidefinite and
    \begin{equation}
        \sum_{i,j=1}^n c_i c_j k(x_i,x_j) = 0 \implies c_1=\dots=c_n = 0.
    \end{equation} 
    The extension to the matrix-valued case reads as follows. A matrix-valued function $\latMatr K :X\times X \to \R^{p\times p}$ is positive semidefinite if, for all $n\in\curly{1,2,\dots}$, ${\dVec{c}}_1,\dots,{\dVec{c}}_n\in\R^p$ and ${x}_1,\dots,{x}_n\in X$, it holds:
    \begin{equation}
        \label{eq:positiveSemidef}
        \sum_{i,j=1}^n {\dVec{c}}_i^\top \latMatr{K}(x_i,x_j) {\dVec{c}}_j \geq 0.
    \end{equation}
    %\xave{(Do we need a symmetry condition for $\latMatr K$? For example, $\latMatr K(x,y) = \latMatr K(y,x)^\top$)} \tobi{I don't think so, as \cite{zastavnyi_analog_2023} never cite this condition, yet we can discuss on this.} \emi{I don't think so.} \xave{I disagree: this definition does not match the univariate one in the case $p=1$, where $k(x,y)=k(y,x)$, see above comment... Furthermore, we work with covariance kernels, so the symmetry should hold: $\latMatr K(x,y) = \latMatr K(y,x)^\top$} \tobi{See my above comment. Here are general definitions, when we specialise them to our case, we add the properties that are satisfied.}
    In addition, $\latMatr{K}$ is strictly positive definite if it is positive semidefinite and the condition
    \begin{equation}
        \sum_{i,j=1}^n {\dVec{c}}_i^\top \latMatr K ({x}_i,{x}_j) {\dVec{c}}_j = 0 
    \end{equation}
    implies that there exists a pair $ i\ne j$ such that ${x}_i={x}_j$ or ${\dVec{c}}_1=\dots={\dVec{c}}_n=0$. {A relevant result for our developments is coming from Proposition 1 in \cite{zastavnyi_analog_2023}, which is formally stated for a neater exposition.} 
    \begin{proposition}[]
        \label{prop:equivalentFormPositDef}
        Let $\latMatr{K} : X \times X \to \R^{p\times p}$. Then, {the following} are equivalent:
        \begin{itemize}
            \item $\latMatr{K}$ is positive semidefinite;
            \item for all $n\in\curly{1,2,\dots}$ and for all ${x}_1,\dots, {x}_n\in X$, the $np \times np$ block matrix $\latMatr C:=[C_{ij}]_{i,j=1}^{np}:=[\latMatr{K}({x}_i,{x}_j)]_{i,j=1}^{{n}}$ is positive semidefinite, \ie, $\dVec{c}^\top \latMatr C \dVec{c} \geq 0$ for any $\dVec{c} \in \mathbb{R}^{np}$. 
        \end{itemize}
        %\emi{NOTATION INCONSISTENCY: if you do not want to use $\R_+$, then you need to consistently use $\{1,2,\ldots\}$ instead of $\mathbb{N}^+$.}   \tobi{Fixed.}
    \end{proposition}
        
    A function $\psi:[0,+\infty)\to {[0,+\infty)}$ is called \emph{completely monotone} if it is continuous on $[0,+\infty)$, infinitely differentiable on $(0,+\infty)$ and for each $i\in \N$ it holds
    \begin{equation*}
        (-1)^i\,\psi^{(i)}(x)\geq 0,
    \end{equation*}
    where $\psi^{(i)}$ denotes the $i^\text{th}$ derivative of $\psi$ and $\psi^{(0)}:=\psi$. A function $\psi:[0,+\infty)\to {[0,+\infty)}$ is completely monotone if and only if it is the Laplace transform of a (unique) finite Borel measure $\mu$ on $[0,+\infty)$, \ie  
    \begin{equation}
        \label{eq:compMonotLaplaceTransform}
        \psi(x)=\int_0^{+\infty} \exp(-tx) \, \mu({\rm d}t), \quad x\in [0,+\infty).
    \end{equation}
%    \xave{I WOULD PREFER TO WRITE $\mu({\rm d}t)$ ($\mu$=MEASURE) RATHER THAN ${\rm d}\mu(t)$ ($\mu$=FUNCTION); SAME FOR OTHER OCCURRENCES OF MEASURES THROUGHOUT THE PAPER.} \emi{That is what I originally did, Xavier. Tobia then changed it. I agree with you.} \tobi{Fine. Fixed it I guess.}
    A function $g:[0,+\infty) \to [0,+\infty)$ is called a \emph{Bernstein function} if $g$ in infinitely differentiable on $[0,+\infty)$ and $g^{(1)}$ is a completely monotone function. \par
    {Our expository material is now completed by reporting Theorem 5 of } \cite{zastavnyi_analog_2023}.
    \begin{theorem}
        \label{theo:MultivariateSchoenberg}
        Let $\psi: [0,+\infty)\to {(0,+\infty)}$ be a {non-constant} completely monotone function, $g:[0,+\infty)\to {[0,+\infty)}$ {be} a Bernstein function and let $\grMatr \Gamma_{{\rVec{Z}}}: X\times X\to [0,+\infty)^{p\times p}$, for a non-empty set $X$. If, for all ${x_1,x_2}\in X$, $\grMatr \Gamma_{{\rVec{Z}}}^\top(x_1,x_2)=\grMatr \Gamma_{{\rVec{Z}}}(x_2,x_1)$ and $\grMatr \Gamma_{{\rVec{Z}}}$ is $\onep$-conditionally negative semidefinite, then, for all ${\xi}>0$,
        \begin{equation*}
            ({x_1,x_2})\mapsto \psi\Big ({\xi} \, g\round{\grMatr \Gamma_{{\rVec{Z}}}({x_1,x_2}}\Big ), \qquad {x_1,x_2}\in X,
        \end{equation*}
        where $\psi$ and $g$ are applied element-wise, is a positive semidefinite function on $X$.
    \end{theorem}
    % \tobi{Divide general results on the construction (put them here) and the properties of our specific metric (below, Remark \ref{prop:actuallyQuasiMetricSpace})}
\section{General Results} 
\label{sec:results}
The following notation will ease the exposition throughout:
\begin{itemize}
    %\item \tobi{we recall from Definition \ref{def:graphWithEE}} that ${\cal G} = (V, E, \curly{\varphi_e}_{e\in E})$ is a graph with Euclidean edges, \emi{this line is redundant - it has already been defined.} \xave{I think that it is a necessary recall, as we used $X$ in the preceding paragraphs and not $\mathcal G$} \tobi{What about adding this at the beginning?}     \emi{I disagree with both. This triple has been defined at the beginning of Definition \ref{def:graphWithEE}. Please check for eventual inconsistencies.} 
    %\textcolor{blue}{(WHY DO WE ABANDON THE NOTATION $\varphi_e$ AND NOW USE $w$? SEEMS NOT TO BE THE SAME CONCEPT: weight $w$ defined in Eq. \ref{eq:laplacianofagraph}, vs. mapping $\varphi_e$ defined in Def. \ref{def:graphWithEE}; see also the ``weight matrix'' in Example \ref{ex:counterExampleDistance} as an additional concept not introduced before, that we may introduce here)} \emi{I agree with Xavier. I think that at the beginning of the project with Tobia we took some decision regarding that, and I do not recall it anymore. We should fix this part.} \emi{N.B.: ${\cal G}$ has been defined at Def 1 (see my correction).} \textcolor{blue}{I suggest removing ``$G$'' and ``${\cal G}$'' from Def. 1, so that it can be introduced here and not be redundant} \tobi{I should have fixed this, please check.} \xave{I also fixed it in Proposition \ref{prop:asymptoticResults} by removing $(V,E,w)$ from this proposition}
    \item $n$ is the number of vertices of {the graph} ${\cal G}$ (cardinality of $V$),
    %\item $G = V \cup E$ is the union of the graph vertices and edges, \emi{this info is not needed as implicit to DEF1}
    \item ${\latMatr{L}}$ is the Laplacian matrix of ${\cal G}$ as per Equation (\ref{eq:laplacianofagraph}), %\emi{WHY do we still have these inconsistencies? We are using Eq. and Equation, we need to choose one of them.},
    \item ${\grMatr{\Sigma}}:={\latMatr{L}}^-$ is the Moore-Penrose inverse of ${\latMatr{L}}$,
    \item ${\grMatr{\Lambda}}$ and ${\latMatr{W}}$ are the eigenvalues and eigenvectors matrices of ${\grMatr{\Sigma}}$, \ie ${\grMatr{\Sigma}} = {\latMatr{W}} {\grMatr{\Lambda}} {\latMatr{W}}^\top$,
    %\item ${\latMatr{I}_n}$ is the $n\times n$ identity matrix,
    %\item $\latMatr{J}_n := \diag{\begin{bmatrix}
    %0, 0, \dots, 0, 1
    %\end{bmatrix}^\top}$ ($n \times n$ diagonal matrix),
    \item for two matrices ${\latMatr{A}},{\latMatr{B}}$ {of the same size} and for a positive integer $p$, we define
    \begin{equation}
        \label{eq:notationM}
        \Mh{p}{{\latMatr{A}}}{{\latMatr{B}}}:=\Mv{p}{{\latMatr{A}}}{{\latMatr{B}}}:={\latMatr{I}}_p \otimes {\latMatr{A}} + (\onepp - {\latMatr{I}}_p) \otimes {\latMatr{B}} = \begin{bmatrix}
            {\latMatr{A}} & {\latMatr{B}} & \dots & {\latMatr{B}}\\
            {\latMatr{B}} & {\latMatr{A}} & \dots & {\latMatr{B}}\\
            \vdots & \vdots & \ddots & \vdots\\
            {\latMatr{B}} & {\latMatr{B}} & \dots & {\latMatr{A}}
        \end{bmatrix}.
    \end{equation}
\end{itemize}

The construction provided in this section extends, to the matrix-valued case, that in \cite{anderes_isotropic_2020}. Specifically, we shall consider a vector-valued random field, $\rVec{Z}$, that is obtained through the identity 
\begin{equation}
    \label{eq:multivariate_random_field_Z}
    \rVec{Z}(u) := \rVec{Z}_{V}(u) + \rVec{Z}_{E}(u), \qquad u \in \mathcal G.
\end{equation}
%where $V$ and $E$ are the vertices and the edges, respectively, of a graph %{${\cal G} = \left( V,E,\{\varphi_e\}_{e\in E} \right)$} 
%with Euclidean edges. 
In turn, the {random fields} $\rVec{Z}_{V}$ and $\rVec{Z}_{E}$ obeys to different construction principles. Specifically, $\rVec{Z}_V$ is a multivariate Gaussian random vector on the vertices and is linearly interpolated on the edges, whilst $\rVec{Z}_E$ is zero on the vertices and adds some variability on the edges. We note that \cite{anderes_isotropic_2020} consider this construction for the case $p=1$ and this entails the use of a specific class of variograms that are then used to build the resistance metric, $d_R$. It is not surprising that our construction will rely on the pseudo-variogram $\grMatr \Gamma_{\rVec{Z}}$ {for the generalisation to the $p$-variate setting with $p>1$}. While these constructions are mathematically involved, a simple description is provided here, with technicalities and proofs deferred to Appendix \ref{A_proofs}. \par

%The following notation will ease the exposition throughout. For two given matrices ${\latMatr{A}},{\latMatr{B}}$ with dimension $n\times n$ and for a positive integer $p$, we define
%\begin{equation}
    %\label{eq:notationM}
    %\Mh{p}{{\latMatr{A}}}{{\latMatr{B}}}:=\Mv{p}{{\latMatr{A}}}{{\latMatr{B}}}:={\latMatr{I}}_p \otimes {\latMatr{A}} + (\onepp - {\latMatr{I}}_p) \otimes {\latMatr{B}} = \begin{bmatrix}
        %{\latMatr{A}} & {\latMatr{B}} & \dots & {\latMatr{B}}\\
        %{\latMatr{B}} & {\latMatr{A}} & \dots & {\latMatr{B}}\\
        %\vdots & \vdots & \ddots & \vdots\\
        %{\latMatr{B}} & {\latMatr{B}} & \dots & {\latMatr{A}}
    %\end{bmatrix}.
%\end{equation}

\subsection{Construction for the Vertices} 
    \label{ssec:constructionZ_V}
    This section provides a construction for the random field $\rVec{Z}_V$ as per the identity (\ref{eq:multivariate_random_field_Z}). While the search for metrics is probably unlimited in terms of alternatives, our spectrum is suitably restricted by providing two \emph{desiderata}, denoted $\mathbb{D}$ throughout.
    \begin{description}
        \item[$\mathbb{D}1$:] The construction (\ref{eq:multivariate_random_field_Z}) needs to provide a pseudo-variogram $\grMatr \Gamma_{\rVec{Z}}$ for which the diagonal entries should obey to the variogram construction as in \cite{anderes_isotropic_2020}. 
        \item[$\mathbb{D}2$:] The multivariate  {random field} $\rVec{Z}_V$ on $V$ should enjoy a homogeneous conditional independence structure. More precisely, we set the block-precision matrix ${\grMatr{\Theta}}$ of $\rVec{Z}_V$ to have the following structure:
        \begin{equation}
            \label{eq:defTheta}
            {\grMatr{\Theta}} = \Mv{p}{\latMatr{Q}}{-\alpha {\latMatr{I}_n}},
        \end{equation}
        where $\alpha$ is a positive real value and $\latMatr{Q}$ is a suitable $n\times n$ matrix that will be defined next. This means that, for $i\ne j$ and $v_1\ne v_2$, ${Z}_{V,i}(v_1)$ and $Z_{V,j}(v_2)$ are conditionally independent given everything else.
    \end{description}
    
    Define the $n\times n$ matrices:
    \begin{align}
        \label{eq:defQ}
        \latMatr{Q}&:=\frac{1}{2} \round{{\latMatr{L}}+\alpha(p-2){\latMatr{I}_n}+\sqrt{\latMatr{L}^2-2\alpha (p-2){\latMatr{L}}+\alpha^2 p^2 {\latMatr{I}_n}}}\\
        \label{eq:defQEigen}
        &=\frac{1}{2} {\latMatr{W}} \round{{\grMatr{\Lambda}}^-+\alpha(p-2){\latMatr{I}_n}+\sqrt{({\grMatr{\Lambda}}^{-})^2-2\alpha(p-2){\grMatr{\Lambda}}^-+\alpha^2 p^2 {\latMatr{I}_n}}} {\latMatr{W}}^\top,\\
        \label{eq:defX}
        {\latMatr{X}}&:=\frac{1}{\alpha(p-1)}\round{{\grMatr{\Sigma}} \latMatr{Q} - {\latMatr{I}_n}}\\
        \label{eq:defXEigen}
        &=\frac{1}{2\alpha (p-1)} {\latMatr{W}} \round{\alpha(p-2){\grMatr{\Lambda}} + \sqrt{{\latMatr{I}_n - \latMatr{J}_n} - 2 \alpha (p-2){\grMatr{\Lambda}} + \alpha^2 p^2 {\grMatr{\Lambda}}^2} -{\latMatr{I}_n}-\latMatr{J}_n} {\latMatr{W}}^\top.
    \end{align}
    Notice that in Equations (\ref{eq:defQEigen}) and (\ref{eq:defXEigen}) the principal square root coincides with the element-wise square root, being its arguments diagonal matrices. {These equations have been established from Equations (\ref{eq:defQ}) and (\ref{eq:defX}) by using the fact that $\latMatr{W}$ is an orthogonal matrix and that $\grMatr{\Lambda} \grMatr{\Lambda}^- = \latMatr{I}_n - \latMatr{J}_n$ insofar as the first $p-1$ diagonal entries of $\grMatr{\Lambda}$ are positive and the last diagonal entry is zero.} 
    In addition, we define the two positive constants
    \begin{equation*}
        %\label{eq:def_k1_k2_constants}
        k_1 := \frac{p-1}{\alpha n p^2} \qquad k_2 := \frac{p^2-p+1}{\alpha n p^2(p-1)},
    \end{equation*} 
    and the $n \times n$ matrices 
    \begin{equation*}
        \Tilde {\grMatr{\Sigma}} := {\grMatr{\Sigma}} + k_1 \onenn \qquad
        \Tilde {\latMatr{X}} := {\latMatr{X}} + k_2 \onenn.
    \end{equation*}
    Although these definitions may appear {peculiar}, they are the only solution to the \emph{desiderata} $\mathbb D$ previously stated and are crucial to define the random field $\rVec{Z}_V$. We start with the vertices, where  $\rVec{Z}_V$ is assumed a zero-mean $p$-variate Gaussian random vector $\rVec{Z}_V\big|_V:V\to \R^p$ having covariance matrix-valued function
    \begin{equation}
        \label{eq:covZ_V_onV}
        \Cov{\rVec{Z}_V(v_1)}{\rVec{Z}_V(v_2)}:=\Mv{p}{\Tilde {\grMatr{\Sigma}}[v_1, v_2]}{\Tilde {\latMatr{X}}[v_1, v_2]} \in \R^{p \times p}.
    \end{equation}
    Checking the positive semidefiniteness of  (\ref{eq:covZ_V_onV}) %needs technical work (see the Appendix) 
    relies on a result of independent interest which implies that both $\grMatr \Theta$ and $\grMatr \Theta^-$ are positive semidefinite matrices. We state this formally below.

    \begin{proposition}
	\label{prop:ThetaQuasiLaplacian}
         The matrix $\grMatr \Theta$ defined in Equation (\ref{eq:defTheta}) and its Moore-Penrose inverse $\grMatr \Theta^-$ are quasi-Laplacian.
    \end{proposition}

    \begin{proposition}
        \label{prop:covZ_V_isACov}
        The function defined in (\ref{eq:covZ_V_onV}) is positive semidefinite.
    \end{proposition}

    Once $\rVec{Z}_V$ is defined over the vertices, the corresponding values over the edges are attained through a sheer linear interpolation. Namely, for $u=(\ul u, \ol u,\delta(u))\in \mathcal G$, we set 
    \begin{equation}
        \label{eq:interpZv}
        \rVec{Z}_V(u):=(1-\delta(u)) \rVec{Z}_V(\ul u) + \delta(u) \rVec{Z}_V(\ol u).
    \end{equation}
    %\xave{I suggest changing ``Lemma'' to ``Proposition'': we are in the main part of the paper, not in an appendix, so the presence of a lemma seems strange.} \tobi{Agreed. Done.}
    
\subsection{Construction for the Edges}     
    \label{ssec:constructionZ_E}
    We define $\rVec{Z}_E$ as a $p$-variate zero-mean Gaussian random field independent of $\rVec{Z}_V$ whose covariance matrix-valued function is:
    \begin{equation}
        \label{eq:covZ_E}
        \Cov{\rVec{Z}_E(u_1)}{\rVec{Z}_E(u_2)}:=\indOne\round{e_1=e_2}\ell(e_1)\round{\delta_1 \wedge \delta_2 - \delta_1 \delta_2}\Mv{p}{1}{\beta} \in {[0,+\infty)}^{p\times p},
    \end{equation}
    for $e_i:=e(u_i)$, $\delta_i:={\delta}(u_i)$, $i\in\curly{1,2}$ and where $\beta\in\square{0,1}$ is a fixed parameter. Notice that $\rVec{Z}_E\big|_{e_1}$ is independent from $\rVec{Z}_E\big|_{e_2}$ whenever $e_1 \ne e_2$, in addition $\ell(e_1)\round{\delta_1 \wedge \delta_2 - \delta_1 \delta_2}$ is the covariance function of a standard Brownian bridge on $[0,\ell(e_1)]$, whilst $\Mh{p}{1}{\beta}$ is positive semidefinite for $\beta\in[0,1]$. This ensures that  (\ref{eq:covZ_E}) is positive semidefinite, as the {element-wise}  product of positive semidefinite functions is itself positive semidefinite.
    
\subsection{Compendium}
    \label{ssec:compendium}
    The following result comes directly from the above definitions.
    \begin{proposition}
        \label{prop:explicit_K_ZV_and_K_ZE}
        Let ${\latMatr{K}}_{\rVec{Z}_V}$, ${\latMatr{K}}_{\rVec{Z}_E}$ be the covariance mappings associated respectively with the {random fields} $\rVec{Z}_V$ and $\rVec{Z}_E$. Then, it is true that 
        \begin{equation}
            {\latMatr{K}}_{\rVec{Z}_V}(u_1,u_2)=\Mv{p}{\dVec{\delta}_{{1,V}}^\top \Tilde {\grMatr{\Sigma}}\dVec{\delta}_{{2,V}}}{\dVec \delta_{{1,V}}^\top \Tilde {\latMatr{X}}\dVec \delta_{{2,V}}} \in \R^{p\times p},
        \end{equation}
        for $u_i=(\ul u_i, \ol u_i, \delta_i)\in \mathcal G$ ($i\in\curly{1,2}$), and where $\dVec{\delta}_{{i,V}}$ is the $n$-dimensional vector whose entries are 
        \begin{equation}
            \label{eq:delta_i}
            {\delta}_{{i,{V}}}:=\begin{cases}
                1-\delta_i \quad &\text{if }v=\ul u_i\\
                \delta_i \quad &\text{if }v=\ol u_i\\
                0 \quad &\text{otherwise}  
            \end{cases}, \quad  v\in V.
        \end{equation}
        Further, it is true that 
        \begin{equation*}
            {\latMatr{K}}_{\rVec{Z}_E}(u_1,u_2)=\Mv{p}{\indOne\round{e_1=e_2}\ell(e_1)\round{\delta_1 \wedge \delta_2 - \delta_1 \delta_2}}{\beta \indOne\round{e_1=e_2}\ell(e_1)\round{\delta_1 \wedge \delta_2 - \delta_1 \delta_2}}.
        \end{equation*}
    \end{proposition}
    As a direct implication of Proposition \ref{prop:explicit_K_ZV_and_K_ZE} in concert with the fact that $\rVec{Z}_V$ and $\rVec{Z}_E$ are independent, we obtain the following.
    \begin{corollary}
        \label{cor:explicitExprK_Z}
        It is true that
        \begin{align}
            \label{eq:covZ_full}
            &{\latMatr{K}}_{\rVec{Z}}(u_1,u_2)={\latMatr{K}}_{\rVec{Z}_V}(u_1,u_2)+{\latMatr{K}}_{\rVec{Z}_E}(u_1,u_2)\nonumber\\
            &=\Mv{p}{\boldsymbol{\delta}_{{1,V}}^\top \Tilde {\grMatr{\Sigma}}\dVec{\delta}_{{2,V}} + \indOne\round{e_1=e_2}\ell(e_1)\round{\delta_1 \wedge \delta_2 - \delta_1 \delta_2}}{\dVec \delta_{{1,V}}^\top \Tilde {\latMatr{X}}\dVec \delta_{{2,V}} + \beta \indOne\round{e_1=e_2}\ell(e_1)\round{\delta_1 \wedge \delta_2 - \delta_1 \delta_2}}.
        \end{align} 
    \end{corollary}

    We are now able to derive the matrix-valued metric associated with the {random field} $\rVec{Z}$ in (\ref{eq:multivariate_random_field_Z}).   To do so, we let  
    $\latMatr{D}: \mathcal G \times \mathcal G \to {[0,+\infty)}^{p \times p}$ be the pseudo-variogram of $\rVec{Z}$, \ie
    \begin{equation}
        \label{eq:definitionOfDistance}     
        \latMatr{D}(u_1,u_2):=\grMatr \Gamma_{\rVec{Z}}(u_1,u_2)=\square{\Var\round{Z_i(u_1)-Z_j(u_2)}}_{i,j=1}^p.
    \end{equation}
    A straightforward application of Equation (\ref{eq:pseudoVariogramMatricial}) in concert with tedious calculations provides the following fact.
    \begin{proposition}
        \label{prop:explicitWritingMultivMetric}
        It is true that
        \begin{align}
            \label{eq:explicitWritingMultivMetric}
            \latMatr{D}(u_1,&u_2)=\Mv{p}{\round{\boldsymbol{\delta}_{ {1,V}}-\boldsymbol{\delta}_{ {2,V}}}^\top \Tilde {\grMatr{\Sigma}} \round{\boldsymbol{\delta}_{ {1,V}}-\boldsymbol{\delta}_{ {2,V}}}}{\boldsymbol{\delta}_{ {1,V}}^\top \Tilde {\grMatr{\Sigma}} \boldsymbol{\delta}_{ {1,V}} + \boldsymbol{\delta}_{ {2,V}}^\top \Tilde {\grMatr{\Sigma}} \boldsymbol{\delta}_{ {2,V}}-2\boldsymbol{\delta}_{ {1,V}}^\top \Tilde {\latMatr{X}} \boldsymbol{\delta}_{ {2,V}}}\nonumber\\
            &+\Mv{p}{\ell(e_1) \delta_1(1-\delta_1)+\ell(e_2)\delta_2(1-\delta_2) - 2\indOne \round{e_1=e_2}\ell(e_1)\round{\delta_1 \wedge \delta_2 - \delta_1 \delta_2}}{\ell(e_1) \delta_1(1-\delta_1)+\ell(e_2)\delta_2(1-\delta_2) - 2\beta \indOne\round{e_1=e_2} \ell(e_1)\round{\delta_1 \wedge \delta_2 - \delta_1 \delta_2}},
        \end{align}
        for $u_1,u_2 \in \mathcal G$, with notation as in Proposition \ref{prop:explicit_K_ZV_and_K_ZE}.
    \end{proposition}    
    Notice that the former line coincides with the contribution of the  {random field} $\boldsymbol{Z}_V$, whilst the latter is the contribution of $\boldsymbol{Z}_E$. 

    We describe below some properties of the metric $\latMatr{D}$ provided through our construction. 

    \begin{proposition}
        \label{prop:actuallyQuasiMetricSpace}
        The matrix-valued expression (\ref{eq:definitionOfDistance}) satisfies properties similar to the (real-valued) quasi-metrics, \ie it satisfies:
        \begin{enumerate}
            \item $D_{ij}(u_1,u_2)\geq 0$;
            \item $D_{ij}(u_1,u_2)=0 \iff u_1=u_2$ and $i=j$;
            \item $\latMatr{D}(u_1,u_2) = \latMatr{D}(u_1,u_2)^\top = \latMatr{D}(u_2,u_1) = \latMatr{D}(u_2,u_1)^\top$ for any $u_1, u_2 \in \mathcal G$;
            \item The mapping $(u_1,u_2) \mapsto \latMatr{D}(u_1,u_0) \latMatr{A} + \latMatr{A}^\top \latMatr{D}(u_0,u_2) - \latMatr{D}(u_1,u_2) - \latMatr{A}^\top \latMatr{D}(u_0,u_0) \latMatr{A}$ is positive semidefinite for any $u_0 \in \mathcal G$ and any $p\times p$ matrix $\latMatr{A}$ such that $\latMatr{A}^\top \onep = \onep$;
            \item {The matrices $\latMatr{D}(u_1,u_2)-\latMatr{D}(u_0,u_0)$ and  $2\latMatr{D}(u_1,u_0)+2\latMatr{D}(u_2,u_0)-\latMatr{D}(u_1,u_2)-3\latMatr{D}(u_0,u_0)$ are symmetric and positive semidefinite for any $u_0,u_1,u_2 \in \mathcal G$.}
        \end{enumerate}
    \end{proposition}
    {Note the resemblance between the fourth property in Proposition \ref{prop:actuallyQuasiMetricSpace} and the second statement in \citet[Proposition 4]{zastavnyi_analog_2023}.}
       
    \begin{proposition}
        \label{prop:ourDistProperties}
        Let ${\cal G}$ be a graph with Euclidean edges and let $\latMatr{D}: \mathcal G\times \mathcal G \to {[0,+\infty)}^{p\times p}$ as defined in (\ref{eq:definitionOfDistance}) and having elements $D_{ij}$, $i,j=1,\ldots, p$. Then, the following properties hold: 
        \begin{enumerate}
            \item for all $i\in\curly{1,\dots,p}$, the distance $D_{ii}$ coincides with the resistance distance defined by \cite{anderes_isotropic_2020}, \ie $D_{ii}(u_1,u_2)=d_{R}(u_1,u_2)$;
            \item the distance $\latMatr{D}$ is matrix-homogeneous, \ie for every $ i\ne j$ and $\,i'\ne j'$, we have $D_{ij}(u_1,u_2)=D_{i'j'}(u_1,u_2)$ and $D_{ii}(u_1,u_2)=D_{i'i'}(u_1,u_2)$.  %\alfr{this is different from the "matrix-homogeneous" condition mentioned in (c), at the beginning of page 4} \tobi{Now it should be ok.} \emi{Now it is okay.}
        \end{enumerate}
    \end{proposition}
    Albeit counter-intuitive, the distance within the same variable is not always less than the same distance among variables. In formulae, $D_{ii}(u_1,u_2)\not\leq D_{ij}(u_1,u_2)$, in general. Appendix \ref{A_B} shows this fact through a worked example.
    
    %The following example shows this fact. \emi{I would put this example in the Appendix and write: {\em This fact is proved in the Appendix through an example.} In fact, this example is covering almost two pages and breaks the narrative while distracting the reader from the flow.} \xave{I agree. Maybe an Appendix 2?} \emi{Yes, I agree.}
    
    The following result illustrates some %\emi{asymptotic: we should clarify that our asymptotics here is w.r.t. $n$ and not to $p$, which is a different kind of asymptotics.}     
    asymptotic properties of the proposed metric. %\xave{PROPOSITION \ref{prop:asymptoticResults} NOT PROVEN IN APPENDIX?} \tobi{Not yet, right.} \emi{Looking forward to studying the proof.} \tobi{Done now.}
    \begin{proposition}
        \label{prop:asymptoticResults}
        Let ${\cal G}$ be a graph with Euclidean edges and let $\latMatr{Q}$ and ${\latMatr{X}}$ as defined in Equations (\ref{eq:defQ})-(\ref{eq:defXEigen}). Then, the asymptotic results {indicated in Table \ref{tab:asympt}} hold.\par %The expressions for $D_{ij}$ ($5^\text{th}$ column) hold for $i\ne j$ (recall that $D_{ii}(u_1,u_2)=d_R(u_1,u_2)$). In addition, concerning 
        %In the $4^\text{th}$ row (\emi{OF WHAT?? Shouldn't all this description be in the caption??}), $\tau \geq \lambda_{n-1}^{-1}$, with $\lambda_{n-1}$ being the smallest positive eigenvalue of $\grMatr{\Sigma}$, while in the last row, $0<\tau \leq \lambda_1^{-1}$. Lastly, whenever in the last column the term $d_R(u_1,u_2)$ appears, the result in that column holds only for $\beta \to 1$ (the other results hold for every $\beta\in[0,1]$).
        %\begingroup
        \begin{table}
            \caption{Asymptotic results with respect to $p$ and $\alpha$. In the $4^\text{th}$ row, $0<\tau \leq \lambda_{n-1}$, with $\lambda_{n-1}$ being the smallest positive eigenvalue of $\grMatr{\Sigma}$, while in the last row, $\tau \geq \lambda_1$. In the last column, the results involving $d_R(u_1,u_2)$ hold only for $\beta \to 1$ (the other result holds for every $\beta\in[0,1]$). \tobi{Emilio, per essere formali qui dovremmo modificare le prime due righe della colonna Q e metterci il limite: sei d'accordo o lasciamo così?}}
            \renewcommand{\arraystretch}{1.2}
            \begin{center}
                \begin{tabular}{|c|c|c|c|c|}
                    \hline
                    $p$ & $\alpha$ & $\latMatr{Q}$ & ${\latMatr{X}}$ & $D_{ij}(u_1,u_2)$ for $i\ne j$\\
                    \hline
                    $+\infty$ & fixed & $\alpha(p-1){\latMatr{I}_n}$ & ${\grMatr{\Sigma}}$ & $d_R(u_1,u_2)$ \\
                    fixed & $+\infty$ & $\alpha(p-1){\latMatr{I}_n}$ & ${\grMatr{\Sigma}}$ & $d_R(u_1,u_2)$ \\
                    fixed & $0^+$ & ${\latMatr{L}}$ & $-\frac{1}{\alpha n (p-1)}\onenn$ & $\frac{2}{\alpha n p}$\\
                    $+\infty$ & $\frac{1}{p\tau}$ & $\frac{\latMatr I_n}{\tau}$ & $\grMatr{\Sigma} - \tau \latMatr I_n$ & $d_R(u_1,u_2) + 2\tau \dVec \delta_{1,V}^\top \dVec \delta_{2,V}$\\
                    $+\infty$ & $\frac{1}{p\tau}$ & $\latMatr L + \frac{1}{n\tau} \onenn $ & $-\frac{\tau}{n}\onenn $ & $d_R(u_1,u_2)+\frac{2\tau}{n} + 2\dVec{\delta}_{1,V}^\top \latMatr X \dVec{\delta}_{2,V}$\\
                    \hline
                \end{tabular}
            \end{center}
            %\endgroup 
            \label{tab:asympt}
        \end{table}
    \end{proposition}
    
    \begin{theorem}
        \label{theo:covFunctionBuilding}
        Let ${\cal G}$ be a graph with Euclidean edges and $\latMatr{D}: \mathcal G\times \mathcal G \to {[0,+\infty)}^{p\times p}$ the distance defined at (\ref{eq:definitionOfDistance}). In addition, let $\psi:[0,+\infty)\to (0,+\infty)$ be a non-constant completely monotone function and $g: [0,+\infty)\to [0,+\infty)$ a Bernstein function. Then, for every $ {\xi}>0$, the mapping  {$\latMatr{K}$ defined as}
        \begin{equation}
            \label{eq:validCovarianceFunction}
            (u_1,u_2)\mapsto  {\latMatr{K}(u_1,u_2)} = \psi\round{ {\xi}\,g\round{\latMatr{D}(u_1,u_2)}},
        \end{equation}
        where $\psi$ and $g$ are applied element-wise, is a valid covariance function. If, in addition, $\psi(0)=1$, then (\ref{eq:validCovarianceFunction}) is a valid correlation function.
    \end{theorem}
    One of the main drawbacks of this construction is that the distance is homogeneous (see Proposition \ref{prop:ourDistProperties}). Consequently, {not only do} the resulting {matrix-valued covariance function (\ref{eq:validCovarianceFunction}) assigns the same marginal covariance for all the variables, \ie $\psi( {\xi} g( {D_{ii}}(u_1,u_2))) = \psi( {\xi} g( {d_R}(u_1,u_2)))$ for any $u_1, u_2\in \mathcal G$ and any $i \in \{1,\ldots,p\}$}, {but also the cross-covariance between two different variables is independent from the variables themselves, \ie $\psi(\xi\, g(D_{ij}(u_1,u_2)))$ does not depend on $i,j$ as long as $i\ne j$.} %\xave{It seems to be stronger: part 1 of Proposition \ref{prop:ourDistProperties} implies that the marginal covariances (diagonal entries) are the same, doesn't it?} \tobi{Right, thx. Is it better now?} \emi{Check on the chain of comments. Fixed.}. 
    Such a restriction can be easily circumvented. Let $[\rho_{ij}]_{i,j=1}^p$ be a collocated correlation matrix, {\em e.g.}, $\rho_{ii}=1$ and $-1\leq \rho_{ij } \le 1$ when $i \ne j$, and the matrix being symmetric and positive semidefinite. Let $\dVec{\sigma} = (\sigma_1 ,\ldots, \sigma_p )^{\top}$ with $\sigma_i>0$ for $i=1,\ldots,p$. Let ${\latMatr{K}}$ be the mapping in (\ref{eq:validCovarianceFunction}) and define the mapping $\widetilde{{\latMatr{K}}} : \mathcal G \times \mathcal G \to \R^{p \times p}$ with entries %$\widetilde{K}_{ij}$ defined through
    $$ \widetilde{{{K}}}_{ij}(u_1,u_2) := \sigma_i \sigma_j \rho_{ij} {{K}}_{ij}(u_1,u_2), \quad i, j = 1,\ldots,p, \quad u_1,u_2 \in \mathcal G.  $$
    Then, by straight application of the Schur product theorem, one gets that $\widetilde{{\latMatr{K}}}$ is a positive semidefinite matrix-valued function. {However, all the marginal covariances (diagonal entries of $\widetilde{{\latMatr{K}}}$) are proportional, and the same happens for all the cross-covariances (off-diagonal entries of $\widetilde{{\latMatr{K}}}$), which is still restrictive.} The results following throughout allow to circumvent such a limitation.
    
\subsection{Generalisations}
    An interesting generalisation is provided below. 
    \begin{theorem} 
        \label{theo:la_raja}
        {Let $\theta \in (0,1]$.} Let $ {\latMatr{F}}: [0,+\infty) \to {\R^{p \times p}}$ be a Borel measure such that, for any fixed $\xi \geq 0$, the matrix $ {{\latMatr{F}} ({\rm d}\xi)}$ 
        is symmetric and positive semidefinite. {Let $\latMatr D$ as in Equation (\ref{eq:definitionOfDistance})} and let $\grMatr \Psi(\cdot): [0,+\infty) \to \mathbb{R}^{p \times p}$ be a continuous mapping determined through
        \begin{equation}
            \label{eq:mapping}
            \grMatr \Psi(\latMatr{D}) := \Big[\int_{0}^{+\infty} \exp(- \xi {{D}_{ij}^\theta}) \,{{F}_{ij}} ({\rm d}\xi)\Big]_{i,j=1}^p,
        \end{equation}
        where ${{F}}_{ij} ({\rm d} \cdot)$ and $D_{ij}$ are the $(i,j)$-th entries of ${\latMatr{F}}({\rm d} \cdot)$ and $\latMatr{D}$, respectively.
        Then, $\grMatr \Psi$ is the covariance function of a Gaussian random field $\rVec Z$ constructed according to the addition principle (\ref{eq:multivariate_random_field_Z}). %<\tobi{Question: is it possible to achieve negative covariance values with this? That is: is $dF_{ij}$ non-negative? If not so, I'd change the codomain of $\Psi$.} \emi{I would not get into that. $\latMatr{F}$ is a signed measure (I do not use this terminology intentionally in the paper). My intuition is that the answer to your question is NO, because we are using an embedding argument into an Hilbert space.}
    \end{theorem}

    As a  {first} application of this result, we adapt a construction that has been originally provided by \cite{porcu2018shkarofsky} in Euclidean spaces and working with the classical Euclidean norm. The {Shkarofsky-Gneiting} family of functions ${\cal SG}_ {\alpha, \beta, \nu}(\cdot) $ is defined on $[0,+\infty)$ through 
    \begin{equation}
        \label{eq:gen_hyper_cov}
        {\cal SG}_ {\alpha, \beta, \nu}(t) = \left  ( 1+ \frac{t}{\beta } \right )^{- \nu/2} \frac{\mathcal{K}_{\nu} \left ( \frac{1}{\sqrt{\alpha}} \sqrt{\beta + t}\right ) }{\mathcal{K}_{\nu} \left ( \sqrt{\frac{\beta}{\alpha}}\right )}, \qquad t \ge 0,
    \end{equation}
    where $\mathcal{K}_{\nu}$ stands for the modified Bessel function of the second kind. 
    Arguments therein show that ${\cal SG}$ is the Laplace transform of a positive and bounded function. The family is very interesting as it admits the special limit cases
    \begin{equation*}
        \label{eq:matern_cov}
        {\cal SG}_ {\alpha, 0, -\nu}(t) = {\cal M}_{\alpha,\nu}(t):= \frac{2^{1-\nu}}{\Gamma(\nu)} \left ( {\frac{t}{\alpha}} \right )^{\nu/2} \mathcal{K}_{\nu} \left ( \sqrt{\frac{t}{\alpha}} \right ), \qquad t \ge 0, 
    \end{equation*}
    for $\alpha, \nu>0$, and 
    \begin{equation*}
        \label{eq:cauchy_cov} 
        {\cal SG}_ {0, \beta, \nu}(t) = {\cal C}_{\beta,\nu}(t):= \left ( 1+ \frac{t}{\beta} \right )^{-\nu}, \qquad t \ge 0, 
    \end{equation*}
    for $\beta, \nu >0$, where ${\cal M}_{\alpha,\nu}(t) = \frac{2^{1+\nu}}{\Gamma(-\nu)} \left ( \frac{\sqrt{t}}{\sqrt{\alpha}} \right )^{-\nu} \mathcal{K}_{-\nu} \left ( \frac{\sqrt{t}}{\sqrt{\alpha}} \right )$, $\alpha>0, \nu>0$ and where ${\cal C}_{\alpha,a, \nu}(t)=\left ( 1+ \frac{t^a}{\alpha} \right )^{-\nu}$, $\alpha>0$, $\nu>0$.
    %Further,
    %\begin{equation}
        %\label{eq:nanc-ninc}  
        %{\cal C}_{\beta,a,\nu}(t) = {\cal SG}_ {0, \beta, \nu}(t^a), \qquad t \ge 0. 
    %\end{equation}
    We {consider the case $p=2$ and} define the bivariate Shkarofsky-Gneiting mapping $\boldsymbol{{\cal SG}}_2 : [0,+\infty)^{{2 \times 2}} \to \R^{2 \times 2}$ through 
    \begin{equation}
        \label{eq:elpatron_model} 
        \boldsymbol{{\cal SG}}_2  (\latMatr{D}) = \begin{bmatrix} 
            \sigma_1^2 \;{\cal SG}_{\alpha_{1},\beta_{1},\nu_{1}}({{D}_{11}^{\theta}}) & \sigma_1 \sigma_2 \rho  \;{\cal SG}_{\alpha_{12},\beta_{12},\nu_{12}}({{D}_{12}^{\theta}}) \\
            & \\
            \sigma_1 \sigma_2 {\rho} \;{\cal SG}_{\alpha_{12},\beta_{12},\nu_{12}}({{D}_{12}^{\theta}}) & \sigma_2^2 \;{\cal SG}_{\alpha_{2},\beta_{2},\nu_{2}}({{D}_{22}^{\theta}})
        \end{bmatrix}, 
    \end{equation}
    %with  ${\cal SG}$ being the Shkarofsky-Gneiting class of functions defined through Equation (\ref{eq:gen_hyper_cov}).  Here, 
    where {$\theta \in (0,1]$, $(\sigma_{1}, \sigma_{2}, \alpha_{1}, \alpha_{2}, \alpha_{12}, \beta_{1}, \beta_{2}, \beta_{12}) \in (0,+\infty)^8$ and $(\rho, \nu_1, \nu_2, \nu_{12}) \in \R^4$}. % are marginal variances, and $\rho$ is a %\tobi{what do you mean with collocated?} \emi{see previous section} 
    %collocated correlation coefficient.
    %\alfr{[[Should the entries of the matrix above be evaluated at the different $D_{ij}^\theta$?]]} \emi{I agree!}
    The following result is a straightforward combination of Theorem \ref{theo:la_raja} with the arguments in Theorem 1 of \cite{porcu2018shkarofsky}, so that the proof is omitted.
    \begin{theorem} 
        \label{theo:el_patron}
        %Let $i=1,2$. Let $\alpha_{i}, \alpha_{12}, \beta_{i},\beta_{12}>0$. Let $\nu_i, \nu_{12} \in \R$.
        Let $\boldsymbol{{\cal SG}}_2$ be the matrix-valued function defined by Equation (\ref{eq:elpatron_model}). % with ${\cal SG}$ being the class of functions defined by Equation (\ref{eq:gen_hyper_cov}). 
        If either
        \begin{description}
            \item[A. {\sc Parsimonious} $\boldsymbol{{\cal SG}}_2$: ] $\alpha_{12}= \frac{1}{2} \frac{\alpha_1 \alpha_2}{(\alpha_1+ \alpha_2)}$, $\beta_{12}= \frac{1}{2}(\beta_1 + \beta_2)$, $\nu_{12}= \frac{1}{2} \left ( \nu_1 + \nu_ 2 \right )$, and $$ | \rho | \le   \sqrt{ \frac{A(\alpha_1,\beta_1,\nu_1) A(\alpha_2,\beta_2,\nu_2)}{A^2(\alpha_{12},\beta_{12},\nu_{12})}},$$ with $A(\cdot,\cdot,\cdot)$ defined as 
            \begin{equation*}
                \label{eq:normalizing}
                A(\alpha,\beta,\nu) = \frac{2^{\nu-1} \left ( \alpha \beta \right )^{\nu/2}}{\mathcal{K}_{\nu} \left ( \sqrt{\frac{\beta}{\alpha}} \right ) },
            \end{equation*} 
            or 
            \item[B. {\sc Full} $\boldsymbol{{\cal SG}}_2$: ] $\alpha_{12} < \frac{1}{2} \frac{\alpha_1 \alpha_2}{(\alpha_1+ \alpha_2)}$, $\beta_{12} > \frac{1}{2}(\beta_1 + \beta_2)$, $\nu_{12} \ne \frac{1}{2} \left ( \nu_1 + \nu_ 2 \right )$, and 
            \begin{equation*}
                \abs{\rho} \le  \sqrt{ \frac{A(\alpha_1,\beta_1,\nu_1) A(\alpha_2,\beta_2,\nu_2)}{A^2(\alpha_{12},\beta_{12},\nu_{12})}  B \left (\frac{1}{\alpha_1}+\frac{1}{\alpha_2}- \frac{1}{2 \alpha_{12}}, \beta_1+\beta_2- 2 \beta_{12}, \nu_1+ \nu_2 - 2 \nu_{12}  \right )  },
            \end{equation*}
            with $B$ defined as 
            \begin{equation*}
                \label{eq:bound2} 
                B(\alpha,\beta, \nu) = \left ( \frac{\left ( \nu - \sqrt{\alpha\beta+\nu^2} \right )}{2 \beta} \right )^{\nu} \exp \left ( \frac{\alpha \beta + \nu (-\nu- \sqrt{\alpha \beta + \nu^2})}{- \nu + \sqrt{\alpha \beta + \nu^2}}   \right ),
            \end{equation*} 
        \end{description}
        then $\boldsymbol{{\cal SG}}_2$ is positive semidefinite on any graph with Euclidean edges.
    \end{theorem}

     {As a second application to a full multivariate setting ($p \geq 2$), consider the mappings $\boldsymbol{{\cal M}}_p$, $\boldsymbol{{\cal C}}_p$ and $\boldsymbol{{\cal SG}}_p$ defined on $[0,+\infty)^{p \times p}$ by
    \begin{equation}
    \label{eq:dios}
        \boldsymbol{{\cal M}}_p(\latMatr{D}) = \Big[\sigma_{ij} \, {\cal M}_{\alpha_{ij},\nu_{ij}}\left(D_{ij}^\theta \right) \Big]_{i,j=1}^p,
    \end{equation}
    \begin{equation}
    \label{eq:dios2}
        \boldsymbol{{\cal C}}_p(\latMatr{D}) = \Big[\sigma_{ij} \, {\cal C}_{\alpha_{ij},\nu}\left(D_{ij}^\theta \right) \Big]_{i,j=1}^p,
    \end{equation}
    \begin{equation}
        \label{eq:dios3}
        \boldsymbol{{\cal SG}}_p(\latMatr{D}) = \Big[\sigma_{ij} \, {\cal SG}_{\alpha_{ij},\beta_{ij},\nu}\left(D_{ij}^\theta \right) \Big]_{i,j=1}^p,
    \end{equation}
    with $[\sigma_{ij}]_{i,j=1}^p$, $[\nu_{ij}]_{i,j=1}^p$, $[\alpha_{ij}]_{i,j=1}^p$ and $[\beta_{ij}]_{i,j=1}^p$ being real symmetric matrices, $\nu$ a real number {and $\theta \in (0,1]$}.} Then {the following statement holds.}
    \begin{theorem} 
        \label{theo:dios}
        Let $\boldsymbol{{\cal M}}_p$ defined by Equation (\ref{eq:dios}). 
        If either
        \begin{description}
            \item[A.] $[\nu_{ij}]_{i,j=1}^p$             and $\big[\nu_{ij}\alpha_{ij}\big]_{i,j=1}^p$ are conditionally negative semidefinite 
            and $\Big[\frac{\sigma_{ij} \nu_{ij}^{\nu_{ij}} \exp(-\nu_{ij})}{ \Gamma(\nu_{ij})}\Big]_{i,j=1}^p$ is positive semidefinite or 
            \item[B.] there exists $\beta>0$ such that $[\nu_{ij}]_{i,j=1}^p$             and $[\alpha_{ij}^{-1}-\beta \nu_{ij}]_{i,j=1}^p$ are conditionally negative semidefinite 
            and $\big[\sigma_{ij} (\beta \alpha_{ij} )^{-\nu_{ij}} \exp(-\nu_{ij}) \big]_{i,j=1}^p$ is positive semidefinite,
        \end{description}
        then $\boldsymbol{{\cal M}}_p$ is positive semidefinite on any graph with Euclidean edges.
    \end{theorem}

%\emi{the notation ${\cal MM}$ is suboptimal. What about $\boldsymbol{{\cal M}}$?}

{\begin{theorem} 
        \label{theo:dios2}
        Let $\boldsymbol{{\cal C}}_p$ defined by Equation (\ref{eq:dios2}) with $\nu > 0$. 
        If $[\beta_{ij}]_{i,j=1}^p$ is conditionally negative semidefinite and $\big[\sigma_{ij} \beta_{ij}^{\nu} \big]_{i,j=1}^p$ is positive semidefinite, then $\boldsymbol{{\cal C}}_p$ is positive semidefinite on any graph with Euclidean edges.
\end{theorem}}

{\begin{theorem} 
        \label{theo:dios3}
        Let $\boldsymbol{{\cal SG}}_p$ defined by Equation (\ref{eq:dios3}). 
        If $[\alpha_{ij}^{-1}]_{i,j=1}^p$ and $[\beta_{ij}]_{i,j=1}^p$ are conditionally negative semidefinite and $\big[\sigma_{ij} (\alpha_{ij} \beta_{ij})^{\nu/2}/{\cal K}_\nu\left(\sqrt{\frac{\alpha_{ij}}{ \beta_{ij}}}\right) \big]_{i,j=1}^p$ is positive semidefinite, then $\boldsymbol{{\cal SG}}_p$ is positive semidefinite on any graph with Euclidean edges.
\end{theorem}}

\section{Euclidean Trees}
\label{sec:linearNetworks}
%\tobi{Here as well, can you pay attention to the domains and codomain of the functions? Try to avoid $\R$ if the actual domain/codomain is \eg $[0,+\infty)$.}
Euclidean trees with a number $m$ of leaves can be embedded on the $m^\prime$-dimensional Euclidean space {with $m^\prime:=\lceil m/2 \rceil$} %\xave{I THINK IT SHOULD BE $\lceil m/2 \rceil$ HERE AND IN WHAT FOLLOWS}
\citep{anderes_isotropic_2020} endowed with the $\ell_1$ distance: for a positive integer ${m^\prime}$ and two points $x,y \in \R^{{m^\prime}}$, the $\ell_1$ or Manhattan-city block distance $d_{M}$ can be defined as 
$d_{M}(x,y) = \sum_{k=1}^{{m^\prime}} |x_k-y_k|$. According to this embedding, for every ${{m^\prime}} \ge 2$, every matrix-valued covariance function depending on $d_M$ can be used as an isotropic covariance function for a {Euclidean tree} %\xave{IS LINEAR NETWORK A SYNONYM OF TREE?} \emi{No. My statement is accurate, but might cause confusion for unexperienced reader.}. 
The surprising fact is that, to the best of our knowledge, the literature has no models of multivariate covariance functions depending on such a metric. %\xave{(What about the exponential covariance? maybe instead of "no models" we can say "few models")} \emi{Do you have a reference for that? Clearly, for $m=1$, all norms are isometric. But a part from this redundancy, I have seen no paper about it.} \xave{THIS REFERENCE: https://link.springer.com/article/10.1007/s11004-006-9055-7. This implies that Loewner completely monotone functions are valid covariances on any Euclidean tree, doesn't it? This could be an extra theorem of this paper!} \emi{These papers are for scalar fields. \cite{anderes_isotropic_2020} mentions alternative facts about them. I would not get into this, we have so much material.}  . 
There is a fact that is even more surprising. Call $\Phi_{m}^p$ the class of continuous mappings ${\mathbf{K}}: [0,+\infty) \to \mathbb{R}^{p \times p}$ such that the composition ${\mathbf{K}}(d_{M}(\cdot,\cdot))$ is positive {semi}definite on a Euclidean tree with $m$ leaves. A characterisation of this class is {missing in the literature} %elusive \tobi{what do you mean with elusive?} \emi{elusive means evasive, or hard to understand},
and we provide it below. 

%\alfr{Emilio, when you consider a tree with $m$ leaves, the results below should be written with $\lfloor m/2 \rfloor$ instead of $m$. Probably, it is easier to modify the text above: $\kappa$ leaves and $m=\lfloor \kappa/2 \rfloor$.} \emi{I do not understand. Can you modify the text to see what you mean? } %\xave{Alfredo seems right: Cambanis et al. (theorem 3.1) obtain the representation (\ref{eq:omega}) where $m$ is the dimension of the space endowed with the Manhattan distance ($m'$ in our notation), not the number of leaves. I made the corrections in the following.}

\begin{proposition} 
    \label{prop:p-Schoen}
    Let $m,p$ be two positive integers. {Let $m^\prime:=\lceil m/2 \rceil$.} Let ${\latMatr{K}}: [0,+\infty) \to \mathbb{R}^{p \times p}$ be continuous with ${\latMatr{K}}(0)= \mathbf{1}_{p \times p}$. Then, ${\latMatr{K}}$ belongs to the class $\Phi_m^p$ if and only if 
    \begin{equation}
        \label{eq:alpha-1} {\latMatr{K}}(t) = \int_0^{+\infty} \omega_{{m^\prime}} (rt ) \,{\latMatr{F}}({\rm d}r), \qquad t \ge 0,
    \end{equation}
    where ${\latMatr{F}}$ is a {bounded} {symmetric} matrix-valued measure such that ${\mathbf{F}} ({\rm d}r)$ is positive {semi}definite for all $r >0$, and where 
    \begin{equation}
        \label{eq:omega}
        \omega_{{m^\prime}}(t) = \frac{\Gamma({m^\prime}/2)}{\sqrt{\pi} \Gamma(({m^\prime}-1)/2)} \int_{1}^{+\infty} \Omega_{{m^\prime}} (\sqrt{v}t) v^{-{m^\prime}/2} (v-1)^{({m^\prime}-3)/2} {\rm d} v,  
    \end{equation}
   with $\Omega_{{m^\prime}}$ being the characteristic function of a random vector that is uniformly distributed over the {unit sphere} $\mathbb{S}^{{m^\prime}-1}$ embedded in $\mathbb{R}^{{m^\prime}}$.
\end{proposition}
Some comments are in order. For $p=1$, the result has been proved by \cite{cambanis1983symmetric}. The following strict inclusion relations hold: 
$$ \Phi_{1}^p \supset \Phi_2^p \supset \cdots \supset \bigcap_{m} \Phi_m^p =: \Phi_{\infty}^p.    $$
A convergence argument from \cite{schoenberg} applies \emph{mutatis mutandis} to assert (no proof needed) that ${\mathbf{K}} \in \Phi_{\infty}^p$ if and only if 
\begin{equation}
    \label{eq:alpha-2} {\mathbf{K}}(t) = \int_0^{+\infty} \exp(-rt) {\mathbf{F}}({\rm d}r), \qquad t \ge 0,
\end{equation}
with ${\mathbf{F}}$ as in Proposition \ref{prop:p-Schoen}. Closed-form expressions for the inner kernel $\omega_{{m^\prime}}$ have been available thanks to  \cite{gneiting1998onalpha}. In particular, we have 
\begin{equation*}
    \omega_2(t) = - \frac{2}{\pi} {\rm si}(t) \quad \text{and} \quad \omega_3(t) = \frac{1}{2} \Bigg (  \frac{\sin t}{t} + \cos t + t \; {\rm si}(t)  \Bigg ), \quad t \ge 0, 
\end{equation*}
where ${\rm si}$ denotes the sine integral function. %(\emi{Please, use the defibrillator to read this part}). 

There exists a wealth of multivariate covariance  models that depend on the Euclidean distance, denoted $\|\cdot\|_{m}$, in $\mathbb{R}^m$. Call $\Psi_m^p$ the class of {continuous} mappings $\latMatr H: [0,+\infty) \to \R^{p \times p}$ with $\latMatr H(0)= \onepp$ and $\latMatr H(\|\cdot\|_m)$ being positive {semi}definite in $\R^m$. The following result proves something extremely useful even for the scalar-case, that apparently was overlooked by \cite{anderes_isotropic_2020}. 

\begin{proposition}
    \label{prop:PleaseEmilioChooseWiseLabels} 
    Let $m,p$ be positive integers. {Let $m^\prime:=\lceil m/2 \rceil$.} Let $\latMatr H$ be a member of the class $\Psi_{2{m^\prime}-1}^p$. Let the mapping ${\latMatr{K}}$ be defined through
    \begin{equation}
        \label{eq:K} 
        {\latMatr{K}}(t) = I^{({m^\prime}-1)} \latMatr H(t), \qquad t \ge 0,
    \end{equation}
    where $I^{({m^\prime}-1)}$ is the $({m^\prime}-1)$ iterated application of the operator $I$, defined as 
    \begin{equation*}  
        (If)(t) := \frac{\int_{t}^{+\infty} f(u) {\rm d} u}{\int_{0}^{+\infty} f(u) {\rm d} u}.
    \end{equation*}
    Then, ${\mathbf{K}}$ belongs to the class $\Phi_m^p$.
\end{proposition}
Propositions \ref{prop:p-Schoen} and \ref{prop:PleaseEmilioChooseWiseLabels} provide a recipe to build new members of the class $\Phi_m^p$. We start by {univariate models, as described in Table \ref{tab:annal-of-stat-2}, where} %describing the elements of Table \ref{tab:annal-of-stat-2}. 
we use the shortcuts $\Psi_m$ and $\Phi_m$ for $\Psi_m^1$ and $\Phi_m^1$, respectively. 

\begingroup
   \renewcommand{\arraystretch}{1.2}
    \begin{center}
        \begin{table}[ht]
            \begin{tabular}{|c|c|c|}
            \hline 
               Parameter restriction  & Member of $\Psi_{2{m^\prime}-1}$ & Member of $\Phi_m$ \\
                 \hline
                 $\nu \ge {m^\prime}$ & $\psi_{\nu}(t) = \left ( 1-t \right )_{+}^{\nu}$ & $ {\phi}_{\nu}(t) = \left ( 1-t \right )_{+}^{\nu+{m^\prime}-1}$ \\
                 $\nu \ge {m^\prime+1}$ & $\psi_{\nu}(t) = \left ( 1-t \right )_{+}^{\nu} \left (1+\nu t \right )$ & ${\phi}_{\nu}(t)= \frac{1}{{m^\prime}} \left (1-t \right )_{+}^{\nu+{m^\prime}-1}  \left ({m^\prime}+\nu t \right )$ \\
                 \hline        
            \end{tabular}
            \vspace{0.1cm}
            \caption{Elements of the class $\Psi_{2{m^\prime}-1}$ and corresponding element of the class $\Phi_{m}$ after iterative application of the operator $I$ as defined through Proposition \ref{prop:PleaseEmilioChooseWiseLabels}. The first column provides the parametric restriction ensuring the corresponding element in the second column to be in $\Psi_{2{m^\prime}-1}$.}
            \label{tab:annal-of-stat-2}
        \end{table}
    \end{center}
\endgroup

We can now use Proposition \ref{prop:PleaseEmilioChooseWiseLabels} in concert with Table \ref{tab:annal-of-stat-2} and the following Lemma, for which a proof is not provided as it follows the same path as Theorem \ref{theo:la_raja}. 
\begin{lemma} 
    \label{lem:tree}
    Let $m,p$ be positive integers and ${\phi}:[0,+\infty)\to \R $ be a member of the class $\Phi_m$. Let ${\latMatr{F}}$ be a matrix-valued measure as in Proposition \ref{prop:p-Schoen}. Then, the scale mixture 
    \begin{equation}
        \label{eq:gen_tree} 
        {\latMatr{K}}(t) = \int_{0}^{+\infty} {\phi}(rt) \,  {\latMatr{F}}({\rm d}r), \qquad t \ge 0, 
    \end{equation}
    provides an element of the class $\Phi_m^p$.
\end{lemma}
Using the arguments in Theorem 1 of \cite{daley2015classes} in concert with Lemma \ref{lem:tree}, one can prove that the multivariate model ${\latMatr{K}}(t) = \left [ K_{ij}(t)\right ]_{i,j=1}^p$ with 
\begin{equation*}
    K_{ij}(t) = \sigma_i \sigma_{ij} \rho_{ij} \phi_{\nu_{ij}}\round{ \frac{t}{b_{ij}}}, \qquad t \ge 0,
\end{equation*}     
is a member of the class $\Phi_{m}^p$ provided the constraints in Theorem 1 of \cite{daley2015classes} hold. Here, ${\phi}_{\nu}$ can be any of the entries in Table \ref{tab:annal-of-stat-2}.\par
We finish this section by providing a direct construction that is based on the kernel $\omega_m $ as being defined through Equation (\ref{eq:omega}).

\begin{proposition}
    \label{prop:infarto_miocardio} 
    Let $m,p$ be positive integers. Let $m^\prime:=\lceil m/2 \rceil$. Let $[\rho_{ij}]_{i,j=1}^p$, $[b_{ij}]_{i,j=1}^p$ and $[\nu_{ij}]_{i,j=1}^p$ be real symmetric matrices such that $b_{ij}>0$ and $\nu_{ij}>2{m^\prime}-1$ for $i,j=1,\ldots,p$. 
    For $r > 0$, define $\latMatr{A}(r)$ with entries
    \begin{equation}
        \label{eq:AA}
        A_{ij}(r) = \frac{\rho_{ij} \Gamma({\frac{\nu_{ij}}{2}})}{\Gamma({\frac{\nu_{ij}}{2}}-{m^\prime}+\frac{1}{2}) \Gamma({{m^\prime}-\frac{1}{2}})} \left ( 1- b_{ij}^2 r^2 \right )_+^{{\frac{\nu_{ij}}{2}}-{{m^\prime}-\frac{1}{2}}}, \quad i,j=1,\ldots,p.
    \end{equation}
    Then, the mapping $\latMatr{K}: [0,+\infty) \to \R^{p \times p} $ with entries 
    \begin{equation}
        \label{eq:Omega-d}
        K_{ij}(x) = \frac{\rho_{ij}}{b_{ij}^{{2{m^\prime}-1}}} \, \omega_{\nu_{ij}} \left ( \frac{x}{b_{ij}}\right ), \qquad x\ge 0, \quad i,j=1,\ldots, p,
    \end{equation}
    belongs to $\Phi_m^p$ provided the matrix $\latMatr{A}(r)$ as defined in (\ref{eq:AA}) is positive semidefinite for any $r > 0$. 
\end{proposition}
Sufficient conditions ensuring $\latMatr{A}(r)$ to be positive semidefinite for every fixed $r \ge 0$ are provided by \citet[Proposition 1]{emery2023schoenberg}. The kernel (\ref{eq:Omega-d}) is especially flexible as it can be used to generate new families of multivariate kernels by using the scale mixture approaches highlighted above. \par
We conclude this section with an illustration to the bivariate setting ($p=2$) for a Euclidean tree with $m=4$ leaves. The entries of the covariance function are given by $K_{11}(d_R(u,v)) = (1- d_R(u,v)/4)_+^3$, $K_{22}(d_R(u,v)) = (1- d_R(u,v))_+^3$, and $K_{12}(d_R(u,v)) = 0.6 \times (1- d_R(u,v)/2.5)_+^3$, for all $u,v\in\mathcal{G}$. Here, the collocated correlation coefficient between the two random field components is $0.6$. \par
\begin{figure}
    \centering
    \includegraphics[width=0.98\textwidth]{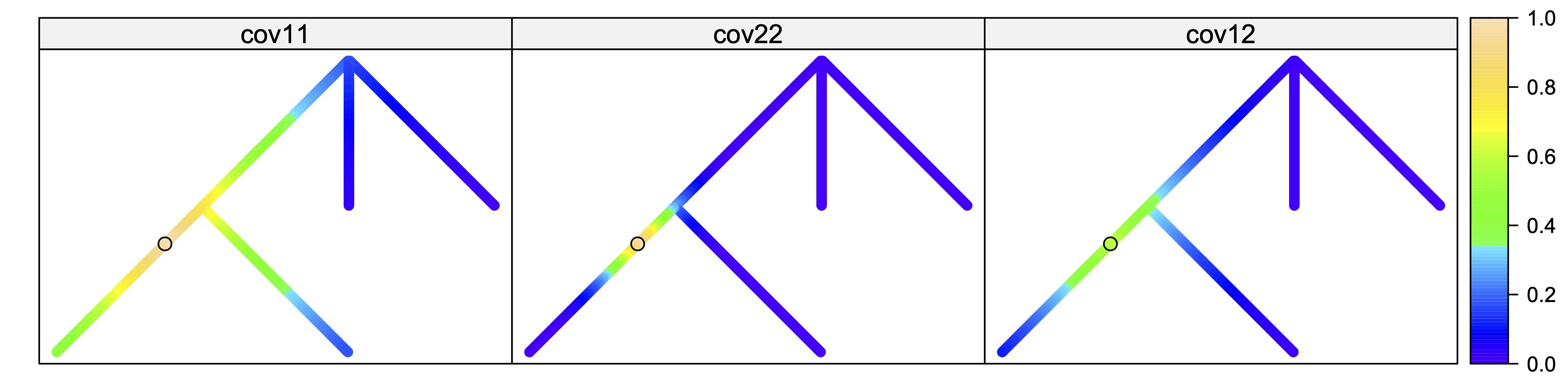}
    \caption{Marginal and cross covariance functions in terms of the resistance metric, over a tree with $m=4$ leaves, measured from a reference location (black circle).}
    \label{fig:cov_plot}
\end{figure}
\begin{figure}
    \centering
    \includegraphics[width=0.98\textwidth]{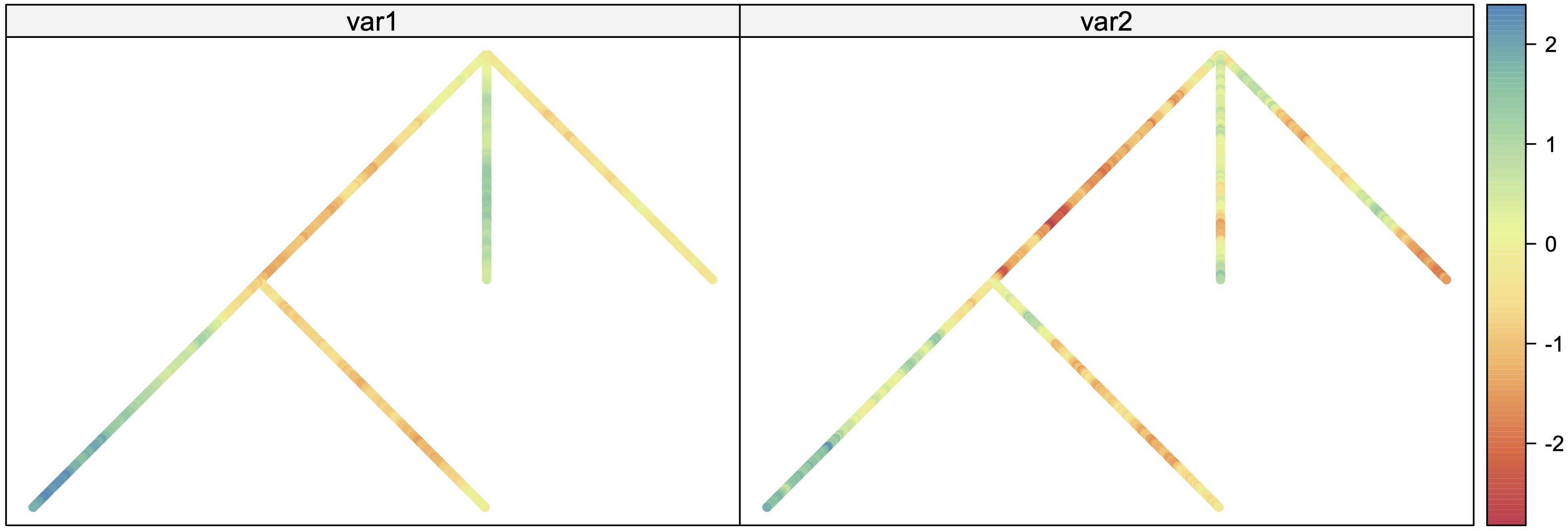}
    \caption{Realisation of a positively correlated bivariate random field on a tree with $m=4$ leaves, with a covariance function of Askey type.}     
    %\xave{SOME TEXT TO INTRODUCE THIS PICTURE WOULD BE WELCOME} \alfr{Let us first decide if we need these kind of figures in the manuscript. If the answer is Yes, I can add a few examples with the corresponding description. It was simulated through Cholesky.} \emi{I like the figure.} \tobi{Alfredo, if we decide to keep this figure, I'd like to make the lines a bit larger, so that it is easier to see the colours. Is it possible? What do you all think?} \xave{I agree} \emi{rather than larger, I think it would be nicer to have thicker lines.}
    \label{fig:simulation}
\end{figure}

Figure \ref{fig:cov_plot} displays the decay of this covariance function in terms of the resistance metric measured from a reference location. It is evident that the first component of the random field has a larger correlation range, as reflected in the left panel: high (low) values of this variable tend to be surrounded by a wider extent of high (low) values.
Figure \ref{fig:simulation} shows a realisation of a bivariate Gaussian random field over $750$ sites on the tree. The Cholesky decomposition of the covariance matrix was used to perform this simulation. \par
%As described above, the first component of the random field has a larger correlation range, as reflected in the left panel: high (low) values of this variable tend to be surrounded by a wider extent of high (low) values. 

\section{Conclusions} 
\label{sec:conclusions}

This paper has shown the intricacies related to the construction to multivariate kernels over generalised frameworks that are represented through a broad class of metric graphs. One might argue that other classes of metric graphs are preferable in certain situations. However, to our knowledge, the only alternative has been proposed by \cite{bolin2022gaussian} and coupled with the substantially different approach of stochastic partial differential equations. It is extremely challenging to attempt for a comparison between these approaches. One important advantage provided by \cite{bolin2022gaussian} for scalar-valued random fields is that an example of a once mean-square differentiable random field is provided. Yet, this is limited to a specific example, while the approach by \cite{anderes_isotropic_2020} --- hence, our approach --- allows to embrace a wealth of examples. Being the first contribution related to vector-valued random fields on metric graphs,  this paper does not have a competitor to compare with. \par
A next intuitive step to this research is represented by vector-valued space-time random fields. The machinery provided in this paper allows for a fairly general building block to start with. Space-time models might also be the building blocks to nonstationary models for networks. \par
The impact of this research is apparent. The works by \cite{baddeley2017stationary} and \cite{moradi} are a clear indication of the importance of this work to modelling point processes over networks. There is a fervent activity related to vector-valued processes in the machine learning community, with the reader referred to \cite{borovitskiy2022isotropic} for details.

%\begin{acks}[Acknowledgements]
\section*{Acknowledgements}
    The authors are grateful to Valeria Simoncini for her support regarding the solutions of matrix quadratic equations. X. Emery acknowledges the support of the National Agency for Research and Development of Chile, through grants ANID Fondecyt 1210050 and ANID PIA AFB230001. Emilio Porcu is grateful to Horst Simon for enlightening discussion about multivariate isometric embeddings.    
%\end{acks}

%\bibliographystyle{authoryear} % Style BST file (imsart-number.bst or imsart-nameyear.bst)
\bibliography{Bib}       % Bibliography file (usually '*.bib')

\newpage
\appendix
\section{Technical Results and Proofs} 
\label{A_proofs}
This part of the Appendix contains proofs and technical results. To help the reader, Figure \ref{fig:resultsRoadMap} provides the road-map of the main results of this manuscript.
\begin{lemma}
        \label{lem:MpMp_expr}
    	Let ${\latMatr{A}},{\latMatr{B}},{\latMatr{C}},{\latMatr{D}}$ be $n \times n$ matrices and let $p\in\curly{1,2,\dots}$. Then
    		\begin{equation}
    			\Mv{p}{{\latMatr{A}}}{{\latMatr{B}}} \Mv{p}{{\latMatr{C}}}{{\latMatr{D}}} = \Mv{p}{{\latMatr{A}}{\latMatr{C}}+(p-1){\latMatr{B}}{\latMatr{D}}}{{\latMatr{A}}{\latMatr{D}}+{\latMatr{B}}{\latMatr{C}}+(p-2){\latMatr{B}}{\latMatr{D}}}.
    		\end{equation}
    \end{lemma}
    \begin{proof}[Proof of Lemma \ref{lem:MpMp_expr}]
        %\xave{(Xavier: OK with the proof)} 
        By the standard properties of the Kronecker product,
        \begin{align*}
            &\round{\latMatr{I}_p \otimes {\latMatr{A}} + (\onepp - \latMatr{I}_p) \otimes {\latMatr{B}}}\round{\latMatr{I}_p \otimes {\latMatr{C}} + (\onepp - \latMatr{I}_p) \otimes {\latMatr{D}}} = \\
            & = (\latMatr{I}_p \otimes {\latMatr{A}}) (\latMatr{I}_p \otimes {\latMatr{C}}) + (\latMatr{I}_p \otimes {\latMatr{A}}) \round{(\onepp - \latMatr{I}_p) \otimes {\latMatr{D}}} \\
            &\quad  + \round{(\onepp-\latMatr{I}_p) \otimes {\latMatr{B}}} \round{\latMatr{I}_p \otimes {\latMatr{C}}} + \round{(\onepp-\latMatr{I}_p) \otimes {\latMatr{B}}} \round{(\onepp - \latMatr{I}_p) \otimes {\latMatr{D}}} \\
            &= (\latMatr{I}_p \latMatr{I}_p) \otimes ({\latMatr{A}}{\latMatr{C}}) + (\latMatr{I}_p(\onepp-\latMatr{I}_p))\otimes ({\latMatr{A}}{\latMatr{D}}) \\
            &\quad  + ((\onepp-\latMatr{I}_p) \latMatr{I}_p) \otimes ({\latMatr{B}}{\latMatr{C}}) + ((\onepp-\latMatr{I}_p)(\onepp-\latMatr{I}_p)) \otimes ({\latMatr{B}}{\latMatr{D}})\\
            &= \latMatr{I}_p \otimes ({\latMatr{A}}{\latMatr{C}})+(\onepp-\latMatr{I}_p) \otimes ({\latMatr{A}}{\latMatr{D}}+{\latMatr{B}}{\latMatr{C}}) + (p \onepp - 2\cdot \onepp+\latMatr{I}_p)\otimes ({\latMatr{B}}{\latMatr{D}})\\
            &=\latMatr{I}_p \otimes ({\latMatr{A}}{\latMatr{C}}) + (\onepp-\latMatr{I}_p) \otimes ({\latMatr{A}}{\latMatr{D}}+{\latMatr{B}}{\latMatr{C}})\\
            &\quad  + \big((p-2)(\onepp-\latMatr{I}_p) + (p-1) \latMatr{I}_p\big) \otimes ({\latMatr{B}}{\latMatr{D}}) \\
            &= \latMatr{I}_p \otimes ({\latMatr{A}}{\latMatr{C}} + (p-1){\latMatr{B}}{\latMatr{D}}) + (\onepp-\latMatr{I}_p) \otimes ({\latMatr{A}}{\latMatr{D}}+{\latMatr{B}}{\latMatr{C}}+(p-2){\latMatr{B}}{\latMatr{D}}).
        \end{align*}
    \end{proof}
    
    \begin{proof}[Proof of Proposition \ref{prop:quasiLaplacianProperties}] 
        Throughout this proof, we assume that $\latMatr{L}\in \R^{n \times n}$, ${\latMatr{A}}\in \R^{n_1 \times n_1}$ and ${\latMatr{C}}\in \R^{n_2 \times n_2}$, with $n_1,n_2>0$ and $n=n_1+n_2$.\par
        \begin{enumerate}
            \item Assume ${\latMatr{B}}=0$. Then $n-1 = \rank({\latMatr{L}})=\rank({\latMatr{A}})+\rank({\latMatr{C}})$. But since $\latMatr L \onen=\dVec 0$, it follows that ${\latMatr{A}} \dVec {1}_{n_1}=\dVec 0$ and ${\latMatr{C}} \dVec 1_{n_2}= \dVec0$. This implies that $\rank({\latMatr{A}})\leq n_1-1$ and $\rank({\latMatr{C}})\leq n_2 -1$. Contradiction.
            \item We will show only that ${\latMatr{C}}$ strictly positive definite, as the proof for the other case relies on the same argument. Since $\latMatr C$ is a principal submatrix of $\latMatr L$, it is positive semidefinite. As a consequence, it is sufficient to show that it is invertible. Assume, by contradiction, that ${\latMatr{C}}$ is singular. This is equivalent to assume that $\exists \dVec y\ne \dVec 0$ such that ${\latMatr{C}} \dVec y=\dVec 0$. Since $\latMatr L$ is positive semidefinite we have, for each $\dVec x\in \R^{n_1}$ and for each $\gamma \in \R$:
            \begin{align*}
                0 &\leq \begin{bmatrix}
                    \dVec x \\
                    \gamma \dVec y
                \end{bmatrix}^\top 
                \begin{bmatrix}
                    {\latMatr{A}} & {\latMatr{B}}\\
                    {\latMatr{B}}^\top & {\latMatr{C}}
                \end{bmatrix}
                \begin{bmatrix}
                    \dVec x \\
                    \gamma \dVec y
                \end{bmatrix}\\
                &=\dVec x^\top {\latMatr{A}} \dVec x + \dVec x^\top {\latMatr{B}} (\gamma \dVec y) + (\gamma \dVec y)^\top {\latMatr{B}}^\top \dVec x + (\gamma \dVec y)^\top {\latMatr{C}} (\gamma \dVec y)=\dVec x^\top {\latMatr{A}} \dVec x + 2 \gamma \dVec x^\top {\latMatr{B}} \dVec y.
            \end{align*}
            If $\dVec x^\top {\latMatr{B}} \dVec y\ne 0$, we can find a suitable $\gamma$ that contradicts the previous inequality. Hence $\dVec x^\top {\latMatr{B}} \dVec y=0$ and therefore, being $\dVec x$ arbitrary, ${\latMatr{B}} \dVec y=\dVec 0$. Consider now the vector $\dVec z:=(0,...,0,\dVec y^\top)^\top$: we have ${\latMatr{L}}\dVec z=\dVec 0$, but clearly $\dVec z$ does not lie in the null space of $\latMatr L$, which is $\curly{\lambda \onen\,:\,\lambda\in \R}$. Contradiction.
            \item Clearly ${\latMatr{L}}/{\latMatr{C}}$ is symmetric. Furthermore, the Guttman rank additivity formula \cite[Equation 0.9.2]{zhang_schur_2005} states that $\rank ({\latMatr{L}})=\rank ({\latMatr{C}}) + \rank({\latMatr{L}}/{\latMatr{C}})$. Here we have $\rank ({\latMatr{L}})=n-1$ and $\rank ({\latMatr{C}})=n_2$, as a consequence $\rank ({\latMatr{L}}/{\latMatr{C}})=n_1-1$. Furthermore, from \citet[Theorem 1.12]{zhang_schur_2005}, being ${\latMatr{C}}$ strictly positive definite and $L$ positive semidefinite, we have that ${\latMatr{L}}/{\latMatr{C}}$ is positive semidefinite as well. The last thing that is left to be shown is that the null eigenvector of ${\latMatr{L}}/{\latMatr{C}}$ is ${\dVec 1}_{n_1}$. Since $\onen$ is the null eigenvector of $\latMatr L$, we get
            \begin{equation*}
                \begin{bmatrix}
                    {\latMatr{A}} & {\latMatr{B}}
                \end{bmatrix}\onen=\dVec 0
                \qquad \text{ and }\qquad \begin{bmatrix}
                    {\latMatr{B}}^\top & {\latMatr{C}}
                \end{bmatrix}\onen=\dVec 0,
            \end{equation*}
            videlicet ${\latMatr{A}} \dVec{1}_{n_1}+{\latMatr{B}} \dVec{1}_{n_2}=\dVec 0$ and ${\latMatr{B}}^\top \dVec{1}_{n_1}+{\latMatr{C}} \dVec{1}_{n_2}=\dVec 0$. Using these relations, the following chain of equivalences concludes the proof:
            \begin{align*}
                &({\latMatr{A}}-{\latMatr{B}} {\latMatr{C}}^{-1} {\latMatr{B}}^\top) \mathbf{1}_{n_1}= \dVec 0 \iff
                {\latMatr{A}} \mathbf{1}_{n_1} = {\latMatr{B}} {\latMatr{C}}^{-1} {\latMatr{B}}^\top \dVec{1}_{n_1}\\
                & \iff -{\latMatr{B}} \dVec{1}_{n_2}={\latMatr{B}} {\latMatr{C}}^{-1}(-{\latMatr{C}} \dVec{1}_{n_2}) \iff -{\latMatr{B}} \dVec{1}_{n_2} = -{\latMatr{B}} \dVec{1}_{n_2}.
            \end{align*}
            \item Consider the eigendecomposition of $\latMatr L$: \begin{equation}
            \label{eq:eigL1}
                {\latMatr{L}}=\latMatr W \begin{bmatrix}
                    \grMatr{\Delta} & \dVec 0\\
                    \dVec 0 & 0
                \end{bmatrix} \latMatr W^\top,
            \end{equation}
            with $\grMatr{\Delta}$ a $(n-1) \times (n-1)$ diagonal matrix with positive diagonal entries.
            Since $\latMatr W$ is an orthogonal matrix, we have 
            \begin{equation}
                \label{eq:eigL2}
                \begin{split}
                    {\latMatr{L}}^-&=\round{\latMatr W \begin{bmatrix}
                        \grMatr{\Delta} & \dVec 0\\
                        \dVec 0 & 0
                    \end{bmatrix} \latMatr W^\top}^-=(\latMatr W^\top)^-\begin{bmatrix}
                        \grMatr{\Delta} & \dVec 0\\
                        \dVec 0 & 0
                    \end{bmatrix}^- \latMatr W^-=\latMatr W\begin{bmatrix}
                        \grMatr{\Delta}^{-1} & \dVec 0\\
                        \dVec 0 & 0
                    \end{bmatrix} \latMatr W^\top.
                \end{split}
            \end{equation}
            Therefore ${\latMatr{L}}^-$ is symmetric, is positive semidefinite and with exactly one null eigenvalue, corresponding to the same null-space eigenvector of $\latMatr L$.
        \end{enumerate}
    \end{proof}

    \begin{lemma}
        \label{lem:eigenvaluesQ}
        Let $\latMatr L$ be an $n \times n$ quasi-Laplacian matrix and ${\grMatr{\Sigma}}={\latMatr{L}}^-$. Let ${\grMatr{\Sigma}} =: \latMatr W \grMatr{\Lambda} \latMatr W^\top$ be the eigendecomposition of ${\grMatr{\Sigma}}$, with $\grMatr{\Lambda} = \diag\round{
            \begin{bmatrix}
                \lambda_1,\dots,\lambda_n
            \end{bmatrix}
        }$ and $\lambda_1\geq\dots\geq\lambda_{n-1}>\lambda_n=0$. Let $\latMatr Q$ be defined as in Subsection \ref{ssec:constructionZ_V}. Then, the eigenvalues $\lambda_n(\latMatr Q)<\lambda_1(\latMatr Q)\leq \dots \leq \lambda_{n-1}(\latMatr Q)$ of $\latMatr Q$ are:
        \begin{equation*}
            \lambda_i(\latMatr Q) = \begin{dcases}
                 \frac{1}{2} \round{\lambda_i^{-1}+\sqrt{\lambda_i^{-2}-2\alpha(p-2)\lambda_i^{-1}+\alpha^2 p^2}+\alpha(p-2)}\qquad &i\ne n\\
                \alpha(p-1)\qquad &\text{$i=n$}.
            \end{dcases}
        \end{equation*}
        In particular the smallest eigenvalue of $\latMatr Q$, $\lambda_n(\latMatr Q)$, is strictly positive, thus $\latMatr Q$ is strictly positive definite. In addition, $\latMatr Q \onen=\alpha(p-1)\onen$.
    \end{lemma}
    
    \begin{proof}[Proof of Lemma \ref{lem:eigenvaluesQ}]
        %\tobi{Added this proof.} \xave{Xavier. OK with the proof}
        Recall the eigendecomposition of matrix $\latMatr Q$ given in Equation (\ref{eq:defQEigen}): it is clear that the eigenvalues of $\latMatr Q$ are the elements on the diagonal of the matrix 
        \begin{equation*}
            \frac{1}{2}\round{\grMatr{\Lambda}^-+\sqrt{(\grMatr{\Lambda}^{-})^2-2\alpha(p-2)\grMatr{\Lambda}^-+\alpha^2 p^2 \latMatr{I}_n}+\alpha(p-2)\latMatr{I}_n},
        \end{equation*}
        which are $\curly{f(\lambda^-)\;:\;\lambda\in \curly{\lambda_1,\dots,\lambda_n}}$, where  {$\lambda_1,\dots,\lambda_n$ are the diagonal entries of $\grMatr{\Lambda}$ (positive, except for $\lambda_n$ that is null) and}
        \begin{equation*}
            f(\lambda^-):=\frac{1}{2}\round{\lambda^-+\sqrt{\round{\lambda^-}^2-2\alpha(p-2)\lambda^-+\alpha^2 p^2}+\alpha(p-2)}.
        \end{equation*}
        Here, we recall that, for a scalar $\lambda\in\R$, $\lambda^-:=\lambda^{-1}$ if $\lambda \ne 0$ and $\lambda^-:=0$ if $\lambda=0$. It is straightforward to show that $\lambda^- \mapsto f(\lambda^-)$ is a strictly increasing function on $\lambda^- \geq 0$. As a consequence, since $\lambda \mapsto \lambda^-$ is strictly decreasing on $\lambda > 0$, $\lambda \mapsto f(\lambda^-)$ is strictly decreasing on $\lambda>0$. Therefore, $\lambda_1\geq\dots\geq \lambda_{n-1}$ are mapped in $\lambda_1(\latMatr Q)=f(\lambda_1^{-1})\leq\dots\leq \lambda_{n-1}(\latMatr Q)= f(\lambda_{n-1}^{-1})$. For the case $\lambda_n=0$, we have $\lambda_n(\latMatr Q)=f(\lambda_n^-)=f(0)=\frac{1}{2}\round{\alpha p + \alpha(p-2)}=\alpha(p-1)$. In addition, since $\lambda_n=0<\lambda_1(\latMatr Q)$, we have $\alpha(p-1) = \lambda_n(\latMatr Q) < \lambda_1(\latMatr Q)$. \\
        Finally, to show that $\latMatr Q \onen = \alpha(p-1) \onen$, it is sufficient to notice from the eigendecomposition of $\latMatr Q$ given in Equation (\ref{eq:defQEigen}) that $\onen$, being (proportional to) the last column of $\latMatr W$, is an eigenvector of $\latMatr Q$. And by the above result it is clear that the corresponding eigenvalue is $\lambda_n(\latMatr Q)=\alpha(p-1)$. This concludes the proof.
    \end{proof}
    
    \begin{lemma}
        \label{lem:usefulIdentities}
        Let $\latMatr L$ be an $n \times n$ quasi-Laplacian matrix,  
        %\tobi{(actually I think a quasi-Laplacian matrix is enough} \xave{Yes, but our proof uses Lemma \ref{lem:eigenvaluesQ} that considers a Laplacian matrix... so I updated both lemmas) } 
        ${\grMatr{\Sigma}}={\latMatr{L}}^-$ and $\alpha >0$. Furthermore, let $\latMatr Q$ and $\latMatr X$ as defined in Subsection \ref{ssec:constructionZ_V}. Then the following holds:
        \begin{enumerate}
            \item \label{item:lemUsefulId_Commute}${\grMatr{\Sigma}}$, $\latMatr L$, $\latMatr Q$ and $\latMatr X$ commute;
            \item \label{item:lemUsefulId_LSigma}${\grMatr{\Sigma}} {\latMatr{L}}= {\latMatr{L}}{\grMatr{\Sigma}} = \latMatr{I}_n-\frac{1}{n} \onenn$;
            \item \label{item:lemUsefulId_Q2} $\latMatr Q^2-\alpha^2 (p-1) \latMatr{I}_n - \latMatr Q {\latMatr{L}}- \alpha (p-2) \latMatr Q + \alpha (p-2) {\latMatr{L}}= \matr 0$;
            \item \label{item:lemUsefulId_SigmaQ2} ${\grMatr{\Sigma}} \latMatr Q^2 - \alpha^2 (p-1) {\grMatr{\Sigma}} - \latMatr Q - \alpha (p-2) {\grMatr{\Sigma}} \latMatr Q + \alpha (p-2) \latMatr{I}_n + \frac{\alpha}{n} \onenn=\matr 0$.
        \end{enumerate}
    \end{lemma}
    \begin{proof}[Proof of Lemma \ref{lem:usefulIdentities}]
        %\xave{(Xavier: OK with the proof)}
        \begin{enumerate}
            \item It is sufficient to notice that all the matrices share the same eigenvector matrix $\latMatr W$, \ie for any ${\latMatr{A}},{\latMatr{B}}\in\curly{{\grMatr{\Sigma}},\latMatr L,\latMatr Q,\latMatr X}$, there are two diagonal matrices $\grMatr{\Lambda}_{\latMatr{A}},\grMatr{\Lambda}_{\latMatr{B}}$ such that ${\latMatr{A}}=\latMatr W \grMatr{\Lambda}_{\latMatr{A}}\matr  W^\top$ and ${\latMatr{B}}=\latMatr W \grMatr{\Lambda}_{\latMatr{B}} \latMatr W^\top$. Therefore 
            \begin{align*}
                {\latMatr{A}}{\latMatr{B}}&=\latMatr W \grMatr{\Lambda}_{\latMatr{A}} \latMatr W^\top \latMatr W \grMatr{\Lambda}_{\latMatr{B}} \latMatr W^\top=\latMatr W \grMatr{\Lambda}_{\latMatr{A}} \grMatr{\Lambda}_{\latMatr{B}} \latMatr W^\top = \latMatr W \grMatr{\Lambda}_{\latMatr{B}} \grMatr{\Lambda}_{\latMatr{A}} \latMatr W^\top\\
                &=\latMatr W\grMatr{\Lambda}_{\latMatr{B}} \latMatr W^\top \latMatr W \grMatr{\Lambda}_{\latMatr{A}} \latMatr W^\top = {\latMatr{B}}{\latMatr{A}}.
            \end{align*}
            \item By exploiting the diagonalisations ${\grMatr{\Sigma}} = \latMatr W \grMatr{\Lambda} \latMatr W^\top$ (Lemma \ref{lem:eigenvaluesQ}) and ${\latMatr{L}}= \latMatr W \grMatr{\Lambda}^- \latMatr W^\top$ (Eqs. (\ref{eq:eigL1}) and (\ref{eq:eigL2})), we get:
            \begin{align*}
                {\grMatr{\Sigma}} {\latMatr{L}}= \latMatr W \grMatr{\Lambda} \grMatr{\Lambda}^- \latMatr W^\top = \latMatr W (\latMatr{I}_n - \latMatr{J}_n) \latMatr W^\top = \latMatr{I}_n - \latMatr W \latMatr{J}_n \latMatr W^\top. 
            \end{align*}
            Given that $\latMatr L$ is a quasi-Laplacian matrix, its eigenvector matrix $\latMatr W$ can be written as
            \begin{equation*}
                \latMatr W = \begin{bmatrix}
                    \latMatr U & c \onen
                \end{bmatrix},
            \end{equation*}
            where $\latMatr U \in \R^{n \times (n-1)}$ and $n c^2=1$. Therefore
            \begin{align}
                \label{eq:WJWtop_expr}
                \latMatr W \latMatr{J}_n \latMatr W^\top = \begin{bmatrix}
                    \latMatr U & c \onen
                \end{bmatrix} \begin{bmatrix}
                    \matr 0_{(n-1)\times(n-1)} & \matr 0_{(n-1)\times 1}\\
                    \matr 0_{1 \times (n-1)} & 1
                \end{bmatrix} \begin{bmatrix}
                    \latMatr U^\top \\
                    c \onen^\top
                \end{bmatrix}=  \frac{1}{n} \onen \onen^\top,
            \end{align}
            and this concludes the proof.
            \item From the definition (\ref{eq:defQ}) of $\latMatr Q$, one has:
            %\begin{equation*}
                %Q:=\frac{1}{2} \round{{\latMatr{L}}+\alpha(p-2)I + \sqrt{{\latMatr{L}}^2 - 2\alpha (p-2) {\latMatr{L}}+ \alpha^2 p^2 I}},
            %\end{equation*}
            %from which 
            \begin{equation*}
                2\latMatr Q - {\latMatr{L}}- \alpha (p-2) \latMatr{I}_n = \sqrt{{\latMatr{L}}^2 - 2\alpha (p-2) {\latMatr{L}}+ \alpha^2 p^2 \latMatr{I}_n}.
            \end{equation*}
            Now, taking the square of both sides (recall that $\latMatr Q$ and $\latMatr L$ commute):
            \begin{align*}
                &4 \latMatr Q^2 + {\latMatr{L}}^2 + \alpha^2 (p-2)^2 \latMatr{I}_n - 4 \latMatr Q{\latMatr{L}}- 4 \alpha (p-2) \latMatr Q + 2 \alpha (p-2) \latMatr L\\
                &= {\latMatr{L}}^2 - 2\alpha (p-2) {\latMatr{L}}+ \alpha^2 p^2 \latMatr{I}_n.
            \end{align*}
            By grouping the common terms and dividing everything by $4$ we get the desired result.
            \item To show the second assertion, it is sufficient to multiply the first equation by ${\grMatr{\Sigma}} = {\latMatr{L}}^-$:
            \begin{align*}
                &{\grMatr{\Sigma}} \latMatr Q^2 - \alpha^2(p-1) {\grMatr{\Sigma}} - \latMatr Q ({\grMatr{\Sigma}} {\latMatr{L}}) - \alpha(p-2) {\grMatr{\Sigma}} \latMatr Q + \alpha (p-2) ({\grMatr{\Sigma}} {\latMatr{L}}) = \matr 0\\
                &{\grMatr{\Sigma}} \latMatr Q^2 - \alpha^2 (p-1) {\grMatr{\Sigma}} - \latMatr Q + \frac{1}{n} \lambda_n(\latMatr Q) \onenn - \alpha (p-2) {\grMatr{\Sigma}} \latMatr Q + \alpha (p-2) \latMatr{I}_n \\
                &\quad - \frac{\alpha (p-2)}{n} \onenn=\matr 0,
            \end{align*}
            where we have used ${\grMatr{\Sigma}} {\latMatr{L}}= \latMatr{I}_n -\frac{1}{n} \onenn$ and $\latMatr Q \onen = \lambda_n(\latMatr Q) \onen= \alpha(p-1) \onen$ (Lemma \ref{lem:eigenvaluesQ}). 
        \end{enumerate}
    \end{proof}

    \begin{lemma}
        \label{lem:sumOfMatricesPD}
        Let $\latMatr X,\latMatr Y\in\R^{n\times n}$ be symmetric matrices. For $i\in\curly{1,\dots,n}$, let $\lambda_i(\latMatr X)$ and $\lambda_i(\latMatr Y)$ denote the $i^\text{th}$ eigenvalues of $\latMatr X$ and $\latMatr Y$, respectively. Let $c\in\R$ such that $\lambda_i(\latMatr X)\geq c$ and $\lambda_i(\latMatr Y)\geq -c$ for all $i\in\curly{1,\dots,n}$. Then $\latMatr X+\latMatr Y$ is positive semidefinite.
    \end{lemma}
    \begin{proof}[Proof of Lemma \ref{lem:sumOfMatricesPD}]
        Let $\latMatr X=\latMatr W_1\grMatr{\Delta}_1 \latMatr W_1^\top$ and $\latMatr Y=\latMatr W_2\grMatr{\Delta}_2 \latMatr W_2^\top$ be the eigendecompositions of $\latMatr X$ and $\latMatr Y$ and let $\dVec v\in\R^n$. In the following lines, we have defined $\dVec a:=\latMatr W_1^\top \dVec v$ and $\dVec b:= \latMatr W_2^\top \dVec v$.
        \begin{align*}
            \dVec v^\top (\latMatr X+\latMatr Y) \dVec v&=\dVec v^\top \latMatr W_1 \Delta_1 \latMatr W_1^\top \dVec v +\dVec v \latMatr W_2 \grMatr{\Delta}_2 \latMatr W_2^\top\dVec v= \dVec a^\top \grMatr{\Delta}_1\dVec a +\dVec b^\top \grMatr{\Delta}_2\dVec b \\
            & = \sum_{i=1}^n\lambda_i(\latMatr X) a_i^2 + \sum_{i=1}^n \lambda_i(\latMatr Y) b_i^2\geq c \sum_{i=1}^n a_i^2 - c \sum_{i=1}^n  b_i^2\\
            &=c(\dVec a^\top \dVec a - \dVec b^\top \dVec b) = c(\dVec v^\top \dVec v-\dVec v^\top \dVec v)=0. 
        \end{align*}
        Videlicet, $\latMatr X+\latMatr Y$ is positive semidefinite.
    \end{proof}
    
    %\begin{proposition}
    	%\label{prop:ThetaQuasiLaplacian}
    	%Let $\latMatr L$ an $n\times n$ Laplacian matrix, ${\grMatr{\Sigma}}:={\latMatr{L}}^-$ its Moore-Penrose inverse and let $\alpha>0$. Then the matrix $\grMatr{\Theta}$ as defined in (\ref{eq:defTheta}) is quasi-Laplacian.
    %\end{proposition}
    
    \begin{lemma}
    	\label{lem:KinverseTheta}
    	Let $\latMatr{L}$ an $n\times n$ Laplacian matrix, let ${\grMatr{\Sigma}}:={\latMatr{L}}^-$ its Moore-Penrose inverse and let $\alpha > 0$. Let $\grMatr{\Theta}$, $\latMatr Q$, $\latMatr X$, $k_1$, $k_2$, $\Tilde{{\grMatr{\Sigma}}}$ and $\Tilde{\latMatr X}$ as described in Subsection \ref{ssec:constructionZ_V}. Finally let:
    	\begin{equation}
            \label{eq:defKTilde}
    		\Tilde {\latMatr K} := \Mv{p}{\Tilde {\grMatr{\Sigma}}}{\Tilde {\latMatr X}}.
    	\end{equation}
    	Then $\Tilde{\latMatr K}=\grMatr{\Theta}^-$.
    \end{lemma}
    %Before proving Proposition \ref{lem:KinverseTheta}, it is in order to formally state its main consequence. The proof is immediate and comes from the above result in concert with the fact that $\grMatr{\Theta}$ is a quasi-Laplacian matrix (see Proposition \ref{prop:ThetaQuasiLaplacian}) and the closure of the set of quasi-Laplacian matrices under Moore-Penrose inversion (Proposition \ref{prop:quasiLaplacianProperties}).
    %\begin{corollary}
        %\label{cor:KTildeIsValidCovariance}
        %$\Tilde{\latMatr K}$ is a quasi-Laplacian matrix as well, and thus a valid covariance matrix.
    %\end{corollary}
    
    \begin{proof}[Proof of Lemma \ref{lem:KinverseTheta}]
        In order to show that $\Tilde{\latMatr K} = \grMatr{\Theta}^-$, by definition of the Moore-Penrose inverse, we need to prove
        \begin{enumerate}
            \item $\grMatr{\Theta} \Tilde{\latMatr K}$ and $\Tilde{\latMatr K} \grMatr{\Theta}$ are symmetric,
            \item $\grMatr{\Theta} \Tilde{\latMatr K} \grMatr{\Theta} = \grMatr{\Theta}$,
            \item $\Tilde{\latMatr K} \grMatr{\Theta} \Tilde{\latMatr K} = \Tilde{\latMatr K}$.
        \end{enumerate}
        To show symmetry, recall from Lemma \ref{lem:usefulIdentities} that all the matrices $\latMatr Q,\alpha \latMatr{I}_n, {\grMatr{\Sigma}}, \latMatr X, \onenn$ commute. As a consequence, $\forall {\latMatr{A}},{\latMatr{B}},{\latMatr{C}},{\latMatr{D}} \in \curly{\latMatr Q,\alpha \latMatr{I}_n, {\grMatr{\Sigma}}, \latMatr X, \onenn}$, it holds (Lemma \ref{lem:MpMp_expr}):
        \begin{align*}
            &\Mv{p}{{\latMatr{A}}}{{\latMatr{B}}}\Mv{p}{{\latMatr{C}}}{{\latMatr{D}}} = \Mv{p}{{\latMatr{A}}{\latMatr{C}}+(p-1){\latMatr{B}}{\latMatr{D}}}{{\latMatr{A}}{\latMatr{D}}+{\latMatr{B}}{\latMatr{C}}+(p-2){\latMatr{B}}{\latMatr{D}}}\\
            &=\Mv{p}{{\latMatr{C}}{\latMatr{A}}+(p-1){\latMatr{D}}{\latMatr{B}}}{{\latMatr{D}}{\latMatr{A}}+{\latMatr{C}}{\latMatr{B}}+(p-2){\latMatr{D}}{\latMatr{B}}}=\Mv{p}{{\latMatr{C}}}{{\latMatr{D}}}\Mv{p}{{\latMatr{A}}}{{\latMatr{B}}}.
        \end{align*}
        In addition, $\Mh{p}{{\latMatr{A}}}{{\latMatr{B}}}^\top=\Mh{p}{{\latMatr{A}}^\top}{{\latMatr{B}}^\top}=\Mh{p}{{\latMatr{A}}}{{\latMatr{B}}}$, as ${\latMatr{A}},{\latMatr{B}}$ are symmetric. Therefore,
        \begin{align*}
            &\round{\grMatr{\Theta} \Tilde{\latMatr K}}^\top = \Tilde{\latMatr K}^\top \grMatr{\Theta}^\top = \Mv{p}{{\grMatr{\Sigma}} + k_1 \onenn}{\latMatr X + k_2 \onenn}^\top \Mv{p}{{\latMatr Q}}{-\alpha \latMatr{I}_n}^\top \\
            &= \Mv{p}{{\grMatr{\Sigma}}^\top + k_1 \onenn^\top}{\latMatr X^\top + k_2 \onenn^\top} \Mv{p}{{\latMatr Q}^\top}{-\alpha \latMatr{I}_n} \\
            &=\Mv{p}{{\grMatr{\Sigma}}}{\latMatr X} \Mv{p}{{\latMatr Q}}{-\alpha \latMatr{I}_n} + \Mv{p}{k_1 \onenn}{k_2 \onenn} \Mv{p}{{\latMatr Q}}{-\alpha \latMatr{I}_n} \\
            &=\Mv{p}{{\latMatr Q}}{-\alpha \latMatr{I}_n} \Mv{p}{{\grMatr{\Sigma}}}{\latMatr X} + \Mv{p}{{\latMatr Q}}{-\alpha \latMatr{I}_n} \Mv{p}{k_1 \onenn}{k_2 \onenn}\\
            &=\Mv{p}{{\latMatr Q}}{-\alpha \latMatr{I}_n} \Mv{p}{{\grMatr{\Sigma}} + k_1 \onenn}{\latMatr X + k_2 \onenn} = \grMatr{\Theta} \Tilde{\latMatr K}.
        \end{align*}
        The same arguments show that $\round{\Tilde{\latMatr K} \grMatr{\Theta}}^\top = \Tilde{\latMatr K} \grMatr{\Theta}$. To prove the other two relations, we start with following equalities. In the next lines, we used Lemmas \ref{lem:MpMp_expr}, \ref{lem:eigenvaluesQ} and \ref{lem:usefulIdentities}.
        \begin{align*}
            \allowdisplaybreaks
            &\grMatr{\Theta} \latMatr{K}=\Mv{p}{{\latMatr Q}}{-\alpha \latMatr{I}_n} \Mv{p}{{\grMatr{\Sigma}}}{\latMatr X}=\Mv{p}{{\latMatr Q}{\grMatr{\Sigma}} - \alpha(p-1)\latMatr X}{{\latMatr Q}\latMatr X-\alpha {\grMatr{\Sigma}}-\alpha(p-2)\latMatr X}\\
            &\quad=\Mv{p}{{\latMatr Q}{\grMatr{\Sigma}} -\alpha(p-1)\frac{1}{-\alpha(p-1)}(\latMatr{I}_n-{\grMatr{\Sigma}} {\latMatr Q})}{\frac{1}{-\alpha(p-1)}({\latMatr Q}-{\grMatr{\Sigma}} {\latMatr Q}^2)-\alpha{\grMatr{\Sigma}} -\alpha(p-2)\latMatr X}=\Mv{p}{\latMatr{I}_n}{-\frac{1}{n(p-1)}\onenn}\\
            &\grMatr{\Theta} \latMatr{K} \grMatr{\Theta} = \Mv{p}{\latMatr{I}_n}{-\frac{1}{n(p-1)}\onenn} \Mv{p}{{\latMatr Q}}{-\alpha \latMatr{I}_n}\\
            &\quad = \Mv{p}{{\latMatr Q}+\alpha (p-1)\frac{1}{n(p-1)}\onenn}{-\alpha \latMatr{I}_n - \frac{1}{n(p-1)}{\latMatr Q}\onenn + \alpha(p-2)\frac{1}{n(p-1)}\onenn}\\
            &\quad =\Mv{p}{{\latMatr Q}+\frac{\alpha}{n}\onenn}{-\alpha \latMatr{I}_n -\frac{\alpha}{n(p-1)}\onenn}\\
            &\grMatr{\Theta}\, \Mv{p}{k_1 \onenn}{k_2 \onenn} = \Mv{p}{k_1 {\latMatr Q} \onenn-\alpha (p-1) k_2 \onenn}{k_2 {\latMatr Q} \onenn - \alpha k_1 \onenn - \alpha k_2 (p-2)\onenn}\\
            &\quad=\Mv{p}{-\alpha(p-1)(k_2-k_1) \onenn}{-\alpha (k_1-k_2)\onenn}\\
            &\grMatr{\Theta}\, \Mv{p}{k_1 \onenn}{k_2 \onenn} \grMatr{\Theta} = \Mv{p}{-\alpha(p-1)(k_2-k_1) \onenn}{-\alpha (k_1-k_2)\onenn} \Mv{p}{{\latMatr Q}}{-\alpha \latMatr{I}_n} \\
            &\quad = \Mv{p}{-\alpha(p-1)(k_2-k_1){\latMatr Q}\onenn + \alpha^2 (p-1)(k_1-k_2)\onenn}{\alpha^2(p-1)(k_2-k_1)\onenn - \alpha(k_1-k_2){\latMatr Q}\onenn +\alpha^2(p-2)(k_1-k_2)\onenn}\\
            &\quad = \Mv{p}{\alpha^2 p (p-1)(k_1-k_2)\onenn}{-\alpha^2 p (k_1-k_2) \onenn}\\
            \intertext{As a consequence:}
            &\grMatr{\Theta}\Tilde{\latMatr K} \grMatr{\Theta} = \grMatr{\Theta} \round{\latMatr K+\Mv{p}{k_1 \onenn}{k_2\onenn}} \grMatr{\Theta}\\
            &\quad = \Mv{p}{{\latMatr Q}+\frac{\alpha}{n}\onenn}{-\alpha \latMatr{I}_n -\frac{\alpha}{n(p-1)}\onenn} + \Mv{p}{\alpha^2 p (p-1)(k_1-k_2)\onenn}{-\alpha^2 p (k_1-k_2) \onenn}\\
            &\quad = \Mv{p}{{\latMatr Q}+\frac{\alpha}{n}\onenn + \alpha^2 p (p-1)(k_1-k_2)\onenn}{ -\alpha \latMatr{I}_n -\frac{\alpha}{n(p-1)}\onenn-\alpha^2 p (k_1-k_2) \onenn}=\Mv{p}{{\latMatr Q}}{-\alpha \latMatr{I}_n} = \grMatr{\Theta}.
            \intertext{Similarly we prove that $\Tilde{\latMatr K} \grMatr{\Theta} \Tilde{\latMatr K} = \Tilde{\latMatr K}$:}
            &\latMatr{K}\grMatr{\Theta} \latMatr{K} = \Mv{p}{{\grMatr{\Sigma}}}{\latMatr X} \Mv{p}{\latMatr{I}_n}{-\frac{1}{n(p-1)}\onenn} = \Mv{p}{{\grMatr{\Sigma}} -\frac{1}{n} \latMatr X \onenn}{-\frac{1}{n(p-1)}{\grMatr{\Sigma}} \onenn + \latMatr X - \frac{p-2}{n(p-2)}\latMatr X \onenn}\\
            &=\Mv{p}{{\grMatr{\Sigma}} + \frac{1}{\alpha n (p-1)}\onenn}{\latMatr X+\frac{p-2}{\alpha n (p-1)^2}\onenn }\\
            &\Mv{p}{k_1 \onenn}{k_2 \onenn} \grMatr{\Theta} \latMatr{K} = \latMatr{K} \grMatr{\Theta} \Mv{p}{k_1 \onenn}{k_2 \onenn}\\
            &\quad = \Mv{p}{k_1 \onenn}{k_2 \onenn}\Mv{p}{\latMatr{I}_n}{-\frac{1}{n(p-1)}\onenn}\\
            &\quad =\Mv{p}{k_1 \onenn -\frac{1}{n}k_2\onenn \onenn}{-\frac{k_1}{n(p-1)}\onenn\onenn+k_2\onenn -\frac{k_2(p-2)}{n(p-1)}\onenn \onenn}\\
            &\quad = \Mv{p}{(k_1-k_2)\onenn}{-\frac{1}{p-1}(k_1-k_2)\onenn }\\
            &\Mv{p}{k_1 \onenn}{k_2 \onenn} \grMatr{\Theta} \Mv{p}{k_1 \onenn}{k_2 \onenn} =\Mv{p}{k_1 \onenn}{k_2 \onenn} \Mv{p}{-\alpha(p-1)(k_2-k_1) \onenn}{-\alpha (k_1-k_2)\onenn}\\
            &\quad = \Mv{p}{-\alpha(p-1)k_1(k_2-k_1) n\onenn-(p-1)\alpha k_2(k_1-k_2)n\onenn}{\big(-\alpha k_1 (k_1-k_2) n  - \alpha(p-1)k_2(k_2-k_1) n -\alpha(p-2)k_2(k_1-k_2) n\big) \onenn}\\
            &\quad = \Mv{p}{\alpha n(p-1)(k_1-k_2)^2 \onenn}{-\alpha n (k_1-k_2)^2\onenn}
            \intertext{Therefore,}
            &\Tilde{\latMatr K} \grMatr{\Theta} \Tilde{\latMatr K} = \latMatr K\grMatr{\Theta} \latMatr K + 2 \, \Mv{p}{k_1 \onenn}{k_2 \onenn} \grMatr{\Theta} \latMatr K + \Mv{p}{k_1 \onenn}{k_2 \onenn} \grMatr{\Theta} \Mv{p}{k_1 \onenn}{k_2 \onenn} = \Tilde{\latMatr K}.
        \end{align*}
    \end{proof}
    
    \begin{proof}[Proof of Proposition \ref{prop:ThetaQuasiLaplacian}]
        We first show that $\grMatr{\Theta}$ is quasi-Laplacian. To this end, we need to show that $\grMatr{\Theta}$ is symmetric, positive semidefinite, $\grMatr{\Theta} \dVec{1}_{np}=\dVec 0$, and $\grMatr{\Theta}$ has only one null eigenvalue.
        \begin{enumerate}
            \item Since $\latMatr{Q}$ and $\alpha \latMatr{I}_n$ are symmetric and $\grMatr{\Theta}:=\Mh{p}{\latMatr{Q}}{-\alpha \latMatr{I}_n}$, $\grMatr{\Theta}$ is symmetric.
            \item First recall that $\grMatr{\Theta} = \matr{\latMatr{I}}_p \otimes (\latMatr{Q}+\alpha \latMatr{I}_n) + \onepp\otimes (-\alpha \latMatr{I}_n)$: we aim to apply Lemma \ref{lem:sumOfMatricesPD} to the matrix $\grMatr{\Theta}$ with $\latMatr X:=\matr{\latMatr{I}}_p \otimes (\latMatr{Q}+\alpha \latMatr{I}_n)$ and $\latMatr Y:=\onepp\otimes (-\alpha \latMatr{I}_n)$. In the following, we write $\dVec{\lambda}({\latMatr{A}})$ to denote the vector of eigenvalues of ${\latMatr{A}}$. In addition, the inequalities apply to all elements: we write $\dVec{\lambda}({\latMatr{A}})\geq c$ as a shortcut of $\forall i\in\curly{1,\dots,n},\,\lambda_i({\latMatr{A}})\geq c$, for an $n\times n$ matrix ${\latMatr{A}}$.
            \begin{align*}
                \dVec{\lambda}(\latMatr{X}) &= \dVec{\lambda}(\matr{\latMatr{I}}_p \otimes (\latMatr{Q}+\alpha \latMatr{I}_n)) = \dVec{\lambda}(\matr{\latMatr{I}}_p) \otimes \dVec{\lambda}(\latMatr{Q}+\alpha \latMatr{I}_n) =\onep \otimes (\dVec{\lambda}(\latMatr{Q}) + \alpha) \\
                &\geq \lambda_n(\latMatr{Q}) +\alpha = \alpha (p-1) +\alpha = \alpha p\\
                \dVec{\lambda}(\latMatr{Y}) &= \dVec{\lambda}(\onepp \otimes (-\alpha \latMatr{I}_n)) = \begin{bmatrix}
                    p\\
                    0\\
                    \vdots\\
                    0
                \end{bmatrix}\otimes \begin{bmatrix}
                    -\alpha \\
                    \vdots \\
                    -\alpha
                \end{bmatrix}\geq -\alpha p.
            \end{align*}
            Therefore $\grMatr{\Theta}$ is positive semidefinite.
            \item To show $\grMatr{\Theta} \dVec{1}_{np} = \dVec 0$, given the special structure of $\grMatr{\Theta}$, it is sufficient to show that 
            \begin{equation*}
                \begin{bmatrix}
                    \latMatr{Q} & -\alpha \latMatr{I}_n & \dots & -\alpha \latMatr{I}_n
                \end{bmatrix} \mathbf{1}_{np} = \dVec 0,
            \end{equation*}
            \ie $\latMatr{Q} \onen = \alpha (p-1) \onen$, {which has already been proved in Lemma \ref{lem:eigenvaluesQ}}. %\xave{(\latMatr{I} think that this stems from Lemma \ref{lem:eigenvalues\latMatr{Q}}; the following could serve as a proof of this lemma.)} \tobi{Agreed! Should be ok now. \latMatr{I} suggest that we removed the following align*, once the proof of Lemma \ref{lem:eigenvalues\latMatr{Q}} has been checked.} \latMatr{I}n the next lines $e_n$ denotes the vector $[0,\dots,0,1]^\top$ and $\star$ denotes any quantity that has no effect on the result. Now:
            %\begin{align*}
            %    \latMatr{Q} \onen &= \frac{1}{2} W \round{\Lambda^- +\alpha (p-2) \latMatr{I} + \sqrt{\Lambda^{-2} - 2\alpha(p-2)\Lambda^{-}+\alpha^2 p^2 \latMatr{I}}} W^\top \onen \\
            %    &=\frac{1}{2} W \round{\Lambda^- +\alpha (p-2) \latMatr{I} + \sqrt{\Lambda^{-2} - 2\alpha(p-2)\Lambda^{-}+\alpha^2 p^2 \latMatr{I}}}(\frac{1}{\sqrt n} e_n)\\
            %    &=\frac{\sqrt{n}}{2} W \round{\Lambda^- e_n +\alpha (p-2) e_n + \sqrt{\Lambda^{-2} - 2\alpha(p-2)\Lambda^{-}+\alpha^2 p^2 \latMatr{I}} \,e_n}\\
            %    &=\frac{\sqrt{n}}{2} W \round{0 +\alpha (p-2) e_n + \begin{bmatrix}
            %        \star & 0\\
            %        0 & \alpha p
            %    \end{bmatrix}\,e_n}\\
            %    &=\frac{\sqrt{n}}{2} W e_n \,(\alpha (p-2) + \alpha p) = \alpha(p-1) \onen.
            %\end{align*}
            %\xave{\latMatr{I} missed the proof of $W^\top \onen = \frac{1}{\sqrt n} e_n$ used in the second equality. Then, in the third equality, the square root $\sqrt n$ changes from denominator to numerator. And in the last equality, \latMatr{I} missed the proof of $W e_n=\frac{2}{\sqrt n} \onen$.}
            \item We now need to show that $\grMatr{\Theta}$ has no other null eigenvalues, \ie the only solution to 
            \begin{equation}
                \label{eq:Theta_v_equal_0}
                \begin{bmatrix}
    			\latMatr{Q} & -\alpha \latMatr{I}_n & \dots & -\alpha \latMatr{I}_n\\
    			-\alpha \latMatr{I}_n & \latMatr{Q} & \dots & -\alpha \latMatr{I}_n\\
    			\vdots & \vdots & \ddots & \vdots\\
    			-\alpha \latMatr{I}_n & -\alpha \latMatr{I}_n & \dots & \latMatr{Q}
    		\end{bmatrix} \begin{bmatrix}
    		    \dVec v_1\\
                \dVec v_2\\
                \vdots\\
                \dVec v_n
    		\end{bmatrix} = \begin{bmatrix}
    		    \dVec 0 \\
                \dVec 0 \\
                \vdots\\
                \dVec 0
    		\end{bmatrix}
            \end{equation}
            is (up to a multiplicative constant) $\dVec v_1=\dots=\dVec v_n=\onen$. By reading Equation (\ref{eq:Theta_v_equal_0}) by blocks, we get the following:  $\forall i \in \curly{1,\dots,n}$
            \begin{align*}
                \latMatr{Q} \dVec v_i - \alpha \sum_{k\ne i} \dVec v_k = \dVec 0.
            \end{align*}
            Taking now the difference for $i=i$ and $i=j$, we get
            \begin{equation*}
                (\latMatr{Q}+\alpha \latMatr{I}_n)(\dVec v_i-\dVec v_j)=\dVec 0
            \end{equation*}
            and, since $\latMatr{Q}+\alpha \latMatr{I}_n$ is strictly positive definite directly by Lemma \ref{lem:eigenvaluesQ}, we obtain that $\dVec v_i=\dVec v_j=:\dVec v$ for all $i,j$. The problem thus becomes
            \begin{equation*}
                (\latMatr{Q}-\alpha(p-1)\latMatr{I}_n) \dVec v = \dVec 0.
            \end{equation*}
            To conclude it is sufficient to notice that the eigenvalues of $\latMatr{Q}-\alpha(p-1) \latMatr{I}_n$ are $\lambda_i(\latMatr{Q}) - \alpha (p-1)$ and, by Lemma \ref{lem:eigenvaluesQ}, all the eigenvalues are not less than $\alpha(p-1)$ and exactly one achieve such a value. As a consequence, the kernel of $\latMatr{Q}-\alpha (p-1)\latMatr{I}_n$ has dimension one and, since we already know that $(\latMatr{Q}-\alpha (p-1)\latMatr{I}_n)\onen=\dVec 0$, the only possibility is $\dVec v=\onen$.
        \end{enumerate}
    
        That $\grMatr{\Theta}^-$ is also a quasi-Laplacian matrix comes from Lemma \ref{lem:KinverseTheta} and the closure of the set of quasi-Laplacian matrices under Moore-Penrose inversion (Proposition \ref{prop:quasiLaplacianProperties}).
        
    \end{proof}
    
    \begin{proof}[Proof of Proposition \ref{prop:covZ_V_isACov}]
          It is sufficient to apply Proposition \ref{prop:equivalentFormPositDef} with $X:=V$. We show that the $np\times np$ matrix $[\latMatr K(v_i,v_j)]_{i,j=1}^{np}$ is positive semidefinite. First, notice that $[\latMatr K(v_i,v_j)]_{i,j=1}^{np}=\Mh{p}{\Tilde{\grMatr{\Sigma}}}{\Tilde{\latMatr X}}=\Tilde{\latMatr K}$. Since $\Tilde{\latMatr K}\in \R^{np \times np}$ is positive semidefinite (see Lemma \ref{lem:KinverseTheta} and Proposition \ref{prop:ThetaQuasiLaplacian}), each matrix ${\latMatr{C}}$ as defined in Proposition \ref{prop:equivalentFormPositDef}, being a principal submatrix of $\matr {\Tilde{\latMatr K}}$, is positive semidefinite as well. This proves the Proposition.
    \end{proof}

    \begin{proof}[Proof of Proposition \ref{prop:actuallyQuasiMetricSpace}]
        \begin{enumerate}
            \item The proof of this point is straightforward, as the pseudo variogram only takes non-negative values.
            \item {The condition $u_1 = u_2$ and $i=j$ is sufficient to have $D_{ij}(u_1,u_2)=0$ due to the fact that the diagonal entries of $\grMatr{\Gamma}_{{\rVec{Z}}}(u_1,u_1)$ are zero. %\textcolor{red}{Xavier: partial proof; need to know if the statement will be corrected}
            Let us show that it is also necessary. From the definition of the pseudo-variogram as per Equation (\ref{eq:pseudoVariogramDef}), $D_{ij}(u_1,u_2) = 0$ if and only if $Z_i(u_1) = Z_j(u_2)$. Since $\rVec{Z} = \rVec{Z}_V+\rVec{Z}_E$ with $\rVec{Z}_V$ and $\rVec{Z}_E$ being mutually independent, a necessary condition 
            %is that $Z_{E,i}(u_1) = Z_{E,j}(u_2)$, which implies $\Cov{Z_{E,i}(u_1)}{Z_{E,j}(u_2)} = \Var \{Z_{E,i}(u_1)\} = \Var \{Z_{E,j}(u_2)\}$. Given the correlation structure (\ref{eq:covZ_E}) of $\rVec{Z}_E$, these equalities hold only if 
            %\begin{itemize}
                %\item both $u_1$ and $u_2$ are vertices, irrespective of $i$ and $j$, or
                %\item $u_i = u_j$ (\ie $e_1=e_2$ and $\delta_1 = \delta_2$), $i \neq j$ and $\beta=1$, or
                %\item $u_i = u_j$ and $i=j$.
            %\end{itemize}
            %A second necessary condition 
            is that $Z_{V,i}(u_1) = Z_{V,j}(u_2)$. Given the correlation structure (\ref{eq:covZ_V_onV}) at the vertices, the variance-covariance matrix $\Tilde{\latMatr{K}}$ of $\rVec{Z}_V\big|_V$ is quasi-Laplacian (Lemma \ref{lem:KinverseTheta} and Proposition \ref{prop:ThetaQuasiLaplacian}) and thus it has only one null eigenvalue associated with the vector $\onenp$. That is, the only non-trivial linear combination of $\rVec{Z}_V\big|_V$ having a zero variance is, up to a constant factor, $\onenp^\top \rVec{Z}_V\big|_V$. Therefore, $Z_{V,i}(u_1) = Z_{V,j}(u_2)$ for $u_1, u_2 \in V$ only if $i=j$ and $u_1=u_2$. Given the linear interpolation on the edges (Equation (\ref{eq:interpZv})), the same holds true for any two points $u_1$ and $u_2$ belonging to $\mathcal G$.  %\emi{Ok this part of the proof}
            \item The claim stems from Equations (\ref{sympseudo}) and (\ref{eq:explicitWritingMultivMetric}). %\textcolor{red}{Xavier: partial proof; need to know if the statement will be corrected}
            \item Let $u_0,u_1,\ldots,u_n \in \mathcal{G}$ and $\dVec{c}_1,\dots, \dVec{c}_n\in\R^p$. Let $\dVec{c}_0=-\latMatr{A}\sum_{i=1}^n \dVec{c}_i$, where $\latMatr{A}$ is a $p \times p$ matrix such that $\latMatr{A}^\top \onep = \onep$. Then, from Equation (\ref{eq:negativesemidef}), one has:
            \begin{equation*}
                \begin{split}    
                    0 &\geq \sum_{i,j=0}^n \dVec{c}_i^\top \grMatr{\Gamma}_{{\rVec{Z}}}(u_i,u_j) \dVec{c}_j \\
                    &= \sum_{i,j=1}^n \dVec{c}_i^\top \grMatr{\Gamma}_{{\rVec{Z}}}(u_i,u_j) \dVec{c}_j + \sum_{i=1}^n \dVec{c}_i^\top \grMatr{\Gamma}_{{\rVec{Z}}}(u_i,u_0) \dVec{c}_0 + \sum_{j=1}^n \dVec{c}_0^\top \grMatr{\Gamma}_{{\rVec{Z}}}(u_0,u_j) \dVec{c}_j + \dVec{c}_0^\top \grMatr{\Gamma}_{{\bm{Z}}}(u_0,u_0) \dVec{c}_0\\
                    &= \sum_{i,j=1}^n \dVec{c}_i^\top \left[\grMatr{\Gamma}_{{\bm{Z}}}(u_i,u_j)-\grMatr{\Gamma}_{{\bm{Z}}}(u_i,u_0)\latMatr{A}-\latMatr{A}^\top \grMatr{\Gamma}_{{\bm{Z}}}(u_0,u_j)+\latMatr{A}^\top\grMatr{\Gamma}_{{\bm{Z}}}(u_0,u_0)\latMatr{A}\right] \dVec{c}_j,
                \end{split}
            \end{equation*}
            \ie, the mapping $(u,v) \mapsto -[\grMatr{\Gamma}_{{\bm{Z}}}(u,v)-\grMatr{\Gamma}_{{\bm{Z}}}(u,u_0)\latMatr{A}-\latMatr{A}^\top\grMatr{\Gamma}_{{\bm{Z}}}(u_0,v)+\latMatr{A}^\top\grMatr{\Gamma}_{{\bm{Z}}}(u_0,u_0)\latMatr{A}]$ is positive semidefinite (Equation \ref{eq:positiveSemidef}). %\emi{This result is closely related to the second statement in \citet[Proposition 4]{zastavnyi_analog_2023}.} \emi{THIS STATEMENT HAS NOTHING TO DO WITH THE PROOF AND SHOULD BE OMITTED.} \tobi{It is true that it is not involved in the proof, yet in my opinion it is better to cite it...} \emi{Not inside the proof.}
            \item Taking $\dVec{c}_1 = \ldots = \dVec{c}_n := \dVec{c}$ and $\latMatr{A} {=\latMatr{I}_p}$, %\emi{identically} \tobi{I'd remove this identically, as $\latMatr A$ is fixed and does not vary.} \emi{variation has nothing to do with identity.} \tobi{so what's the meaning of identically?} equal to the identity matrix, 
            one gets
            \begin{equation*}
                \begin{split}    
                    0 &\geq \dVec{c}^\top \sum_{i,j=1}^n  \left[\grMatr{\Gamma}_{{\bm{Z}}}(u_i,u_j)-\grMatr{\Gamma}_{{\bm{Z}}}(u_i,u_0)-\grMatr{\Gamma}_{{\bm{Z}}}(u_0,u_j)+\grMatr{\Gamma}_{{\bm{Z}}}(u_0,u_0)\right] \dVec{c}.
                \end{split}
            \end{equation*}
            Accounting for the fact that $\grMatr{\Gamma}_{{\bm{Z}}}=\grMatr{\Gamma}_{{\bm{Z}}}^\top$ and $\grMatr{\Gamma}_{{\bm{Z}}}(u_0,u_0)=\grMatr{\Gamma}_{{\bm{Z}}}(u_1,u_1)=\grMatr{\Gamma}_{{\bm{Z}}}(u_2,u_2)$, one obtains, for $n=1$ and $u_0=u_2$:
            \begin{equation*}
                \begin{split}    
                    0 &\geq \dVec{c}^\top \left[\grMatr{\Gamma}_{{\bm{Z}}}(u_0,u_0)-\grMatr{\Gamma}_{{\bm{Z}}}(u_1,u_2)\right] \dVec{c},
                \end{split}
            \end{equation*}  
            and for $n=2$:
            \begin{equation*}
                \begin{split}    
                    0 &\geq \dVec{c}^\top \left[\grMatr{\Gamma}_{{\bm{Z}}}(u_1,u_2)+3\grMatr{\Gamma}_{{\bm{Z}}}(u_0,u_0)-2\grMatr{\Gamma}_{{\bm{Z}}}(u_1,u_0)-2\grMatr{\Gamma}_{{\bm{Z}}}(u_2,u_0)\right] \dVec{c},
                \end{split}
            \end{equation*}
            \ie, $\grMatr{\Gamma}_{{\bm{Z}}}(u_1,u_2)-\grMatr{\Gamma}_{{\bm{Z}}}(u_0,u_0)$ and $2\grMatr{\Gamma}_{{\bm{Z}}}(u_1,u_0)+2\grMatr{\Gamma}_{{\bm{Z}}}(u_2,u_0)-\grMatr{\Gamma}_{{\bm{Z}}}(u_1,u_2)-3\grMatr{\Gamma}_{{\bm{Z}}}(u_0,u_0)$ are positive semidefinite matrices. %emi{PROOF checked.}
            }
        \end{enumerate}
    \end{proof}

    \begin{proof}[Proof of Proposition \ref{prop:ourDistProperties}]
        %\xave{(Xavier: OK with the proof)} 
        \begin{enumerate}    
            \item We split the computation of the distances (the one proposed here, ${\latMatr{D}}$, and the one defined by \cite{anderes_isotropic_2020}, $d_{R}$) in two: the contribution (\ie the variogram) of the $\rVec{Z}_V$ process and the one of $\rVec{Z}_E$, denoted via the superscripts $V$ and $E$ respectively. We show that ${\latMatr{D}}_{ii}^V(u_1,u_2)=d_{R}^V(u_1,u_2)$ and ${\latMatr{D}}_{ii}^E(u_1,u_2)=d_{R}^E(u_1,u_2)$. The latter comes straightforwardly by comparing the diagonal entries of Equation (\ref{eq:covZ_E}), \ie $\indOne\round{e_1=e_2}\ell(e_1)\round{\delta_1 \wedge \delta_2 - \delta_1 \delta_2}$, with the kernel $R_e$ in \citet[Equation 15]{anderes_isotropic_2020}. \par
            Thus, let us show that ${\latMatr{D}}_{ii}^V(u_1,u_2)=d_{R}^V(u_1,u_2)$. Let $u_1,u_2\in \mathcal G$ {and $k_{V,ii}$ be the $i$-th diagonal entry of {$\latMatr K_{\rVec{Z}_V}$}. Then, on account of Proposition \ref{prop:explicit_K_ZV_and_K_ZE}}:
            \begin{align*}
                {\latMatr{D}}_{ii}^V(u_1,u_2)&=k_{V,ii}(u_1,u_1)+k_{V,ii}(u_2,u_2)-2k_{V,ii}(u_1,u_2)\\
                &=\dVec{\delta}_{1,V}^\top \Tilde{\grMatr{\Sigma}} \dVec{\delta}_{1,V}+ \dVec{\delta}_{2,V}^\top \Tilde{\grMatr{\Sigma}} \dVec{\delta}_{2,V} -2 \dVec{\delta}_{1,V}^\top \Tilde{\grMatr{\Sigma}} \dVec{\delta}_{2,V}\\
                &=\round{\dVec{\delta}_{1,V}-\dVec{\delta}_{2,V}}^\top \round{{\grMatr{\Sigma}}+k_1 \onenn} \round{\dVec{\delta}_{1,V}-\dVec{\delta}_{2,V}}\\
                &=\round{\dVec{\delta}_{1,V}-\dVec{\delta}_{2,V}}^\top {\grMatr{\Sigma}} \round{\dVec{\delta}_{1,V}-\dVec{\delta}_{2,V}} + k_1\round{\dVec{\delta}_{1,V}-\dVec{\delta}_{2,V}}^\top \onenn \round{\dVec{\delta}_{1,V}-\dVec{\delta}_{2,V}}\\
                &=\round{\dVec{\delta}_{1,V}-\dVec{\delta}_{2,V}}^\top {\grMatr{\Sigma}} \round{\dVec{\delta}_{1,V}-\dVec{\delta}_{2,V}},
            \end{align*}
            where in the last step we have used $\dVec{\delta}_{1,V}^\top \onen = \dVec{\delta}_{2,V}^\top \onen = 1$ owing to Equation (\ref{eq:delta_i}).\par %$\bm{\delta}_i^\top \onen = \onen^\top \bm{\delta}_i = (1-\delta_i) + (\delta_i) = 1$.\par
            It is left to show that this quantity coincides with the resistance counterpart of \cite{anderes_isotropic_2020}. They define their covariance matrix, here denoted by ${\grMatr{\Sigma}}_{\latMatr{A}}$, as the (proper) inverse of a modified Laplacian matrix $\latMatr{L}$, where a one has been added on an arbitrary diagonal entry. In formulae: ${\grMatr{\Sigma}}_{\latMatr{A}}:=({\latMatr{L}}+\dVec e_j \dVec e_j^\top)^{-1}$, where $\dVec e_j$ is the $j^\text{th}$ vector of the canonical basis of $\R^n$. It follows directly from \citet[Equation 13]{anderes_isotropic_2020} that $d_{R}^V=\round{\dVec{\delta}_{1,V}-\dVec{\delta}_{2,V}}^\top {\grMatr{\Sigma}}_{\latMatr{A}} \round{\dVec{\delta}_{1,V}-\dVec{\delta}_{2,V}}$.\par
            Let us show that, for a given quasi-Laplacian matrix $\latMatr{L}$ and for any $\dVec x \in \R^n$ such that $\onen^\top \dVec x \ne 0$, it holds:
            \begin{equation}
                \label{eq:thesisRank1Update}
                \round{\dVec{\delta}_{1,V}-\dVec{\delta}_{2,V}}^\top \round{{\latMatr{L}}+\dVec x\,\dVec x^\top}^- \round{\dVec{\delta}_{1,V}-\dVec{\delta}_{2,V}} = \round{\dVec{\delta}_{1,V}-\dVec{\delta}_{2,V}}^\top {\latMatr{L}}^- \round{\dVec{\delta}_{1,V}-\dVec{\delta}_{2,V}}.
            \end{equation}
             In order to express the Moore-Penrose inverse of the rank-1 update ${\latMatr{L}}+\dVec x\,\dVec x^\top$, we exploit Theorem 1 of \cite{meyer_generalized_1973} with ${\latMatr{A}}:=\latMatr L$ and $\dVec c:=\dVec d:=\dVec x$. Firstly, let us show that the hypotheses of the theorem are satisfied: we have to prove that $\dVec x\not\in R({\latMatr{L}})$, where $R$ denotes the range or column space. $R({\latMatr{L}})$ is perpendicular to the kernel of $\latMatr L$, that is (recall that $\latMatr L$ is a quasi-Laplacian matrix) $\ker({\latMatr{L}})=\curly{h\,\onen\,:\,h\in\R}$. Now, $\dVec x\in R({\latMatr{L}})\iff \dVec x \perp \ker({\latMatr{L}}) \iff \onen^\top \dVec x=0$, yet $\onen^\top \dVec x \ne 0$ by hypothesis. Therefore $\dVec x\not\in R({\latMatr{L}})$, \ie we can apply the above-mentioned theorem. In the following computations, we have adopted the notation used by \citet[Section 2]{meyer_generalized_1973}, \ie:
             \begin{align*}
                 \dVec k&:={\latMatr{L}}^-\dVec x & \dVec h&:=\dVec x^\top {\latMatr{L}}^-\\
                 \dVec u&:=\round{\latMatr{I}_n-\latMatr L{\latMatr{L}}^-}\dVec x=\frac{1}{n} \onenn \dVec x & \dVec v&:= \dVec x^\top \round{\latMatr{I}_n-{\latMatr{L}}^- \latMatr L} = \frac{1}{n} \dVec x^\top \onenn\\
                 \beta &:= 1+\dVec x^\top \latMatr L^- \dVec x.
             \end{align*}
             In addition, we have $\norm{\dVec u}^2=\norm{\dVec v}^2= \dVec u^\top \dVec u = \frac{1}{n^2}\onen^\top (\onen^\top \dVec x)(\onen^\top \dVec x) \onen = \frac{1}{n}(\onen^\top \dVec x)^2\ne 0$. Let us finally apply Equation (3.1) of \cite{meyer_generalized_1973}:
             \begin{align*}
                 \round{{\latMatr{L}}+\dVec x \dVec x^\top}^-&={\latMatr{L}}^- - \dVec k \dVec u^- - \dVec v^- \dVec h +\beta \dVec v^- \dVec u^-\\
                 &={\latMatr{L}}^--\frac{\dVec k\dVec u^\top}{\norm{\dVec u}^2} - \frac{\dVec v^\top \dVec h}{\norm{\dVec v}^2} + \beta \frac{\dVec v^\top \dVec u^\top}{\norm{\dVec v}^2 \norm{\dVec u}^2}\\
                 &={\latMatr{L}}^- - \frac{{\latMatr{L}}^- \dVec x \onen^\top}{\onen^\top \dVec x} - \frac{\onen \dVec x^\top {\latMatr{L}}^-}{\onen^\top \dVec x} + \round{1+ \dVec x^\top {\latMatr{L}}^- \dVec x} \frac{\onenn}{\round{\onen^\top \dVec x}^2}.
             \end{align*}
             To show (\ref{eq:thesisRank1Update}), it is now sufficient to prove that
             \begin{equation*}
                 \round{\dVec{\delta}_{1,V}-\dVec{\delta}_{2,V}}^\top \round{\round{1+ \dVec x^\top {\latMatr{L}}^- \dVec x} \frac{\onenn}{\round{\onen^\top \dVec x}^2}- \frac{{\latMatr{L}}^- \dVec x \onen^\top}{\onen^\top \dVec x} - \frac{\onen \dVec x^\top {\latMatr{L}}^-}{\onen^\top \dVec x}}\round{\dVec{\delta}_{1,V}-\dVec{\delta}_{2,V}}=0.
             \end{equation*}
             It is possible to see this by noticing that $\onen^\top \round{\dVec{\delta}_{1,V}-\dVec{\delta}_{2,V}}=\onen^\top \dVec{\delta}_{1,V} - \onen^\top \dVec{\delta}_{2,V}= 1-1=0$ and similarly $\round{\dVec{\delta}_{1,V}-\dVec{\delta}_{2,V}}^\top \onen=0$, therefore each term of the above equation equals 0.\par
             As a consequence,
             \begin{align*}
                 d_{R}^V(u_1,u_2)&=\round{\dVec{\delta}_{1,V}-\dVec{\delta}_{2,V}}^\top {\grMatr{\Sigma}}_{\latMatr{A}} \round{\dVec{\delta}_{1,V}-\dVec{\delta}_{2,V}} \\
                 &= \round{\dVec{\delta}_{1,V}-\dVec{\delta}_{2,V}}^\top \round{{\latMatr{L}}+\dVec e_j \dVec e_j^\top}^- \round{\dVec{\delta}_{1,V}-\dVec{\delta}_{2,V}} \\
                 &= \round{\dVec{\delta}_{1,V}-\dVec{\delta}_{2,V}}^\top {\latMatr{L}}^- \round{\dVec{\delta}_{1,V}-\dVec{\delta}_{2,V}} = {D}_{ii}^V(u_1,u_2).
             \end{align*}
             
            \item The proof of this point is straightforward from Equation (\ref{eq:explicitWritingMultivMetric}): notice that ${D}_{ij}(u_1,u_2)$, for $i\ne j$, is the off-diagonal element of ${\latMatr{D}}(u_1,u_2)$, whilst ${D}_{ii}(u_1,u_2)$ is the main diagonal one.
        \end{enumerate}
    \end{proof}

    \begin{proof}[Proof of Proposition \ref{prop:asymptoticResults}]
    %    \tobi{This proof should be complete now, please check it.} \xave{OK for me.}
        We prove each row of the table, in turn. In the following computations, $\sim_a$ stands for ``asymptotically equivalent to''.
        \begin{enumerate}
            \item Fix $\alpha>0$ and let $p\to +\infty$. Considering Equations (\ref{eq:defQEigen}) and (\ref{eq:defXEigen}), we have: 
            \begin{align}
                \label{eq:asymptQ1}
                &\grMatr{\Lambda}^- + \alpha(p-2)\latMatr{I}_n + \sqrt{(\grMatr{\Lambda}^-)^2-2\alpha(p-2)\grMatr{\Lambda}^-+\alpha^2 p^2 \latMatr{I}_n} \\
                &\quad \sim_a \alpha (p-2) \latMatr{I}_n + \alpha p \latMatr{I}_n = 2\alpha (p-1) \latMatr{I}_n\nonumber\\
                \label{eq:asymptX1}
                &\alpha(p-2)\grMatr{\Lambda} + \sqrt{\latMatr{I}_n - \latMatr{J}_n - 2\alpha (p-2) \grMatr{\Lambda} + \alpha^2 p^2 \grMatr{\Lambda}^2} -\latMatr{I}_n - \latMatr{J}_n \\
                &\quad \sim_a \alpha (p-2) \grMatr{\Lambda} + \alpha p \grMatr{\Lambda} = 2 \alpha (p-1) \grMatr{\Lambda}. \nonumber
            \end{align}
            As a consequence, $\latMatr Q \sim_a \alpha (p-1) \latMatr{I}_n$ and $\latMatr X \to \grMatr{\Sigma}$. Furthermore, $k_1,k_2\to 0$ and thus $\Tilde{\grMatr{\Sigma}}, \tilde{\latMatr X} \to \grMatr{\Sigma}$ as $p\to+\infty.$ From this last relation, we obtain that the process $\rVec Z_V$ satisfies $\Corr{Z_{V,i}}{Z_{V,j}}\to 1$ for all $i,j\in\curly{1,\dots,p}$. Throughout the whole proof, the expression for the distance changes only via the first line of Equation (\ref{eq:explicitWritingMultivMetric}), reported below for an neat exposition.
            \begin{align}
                \label{eq:appZ_V_contributionOnD}
                &{\dVec{\delta}_{ {1,V}}^\top \Tilde {\grMatr{\Sigma}} \dVec{\delta}_{ {1,V}} + \dVec{\delta}_{ {2,V}}^\top \Tilde {\grMatr{\Sigma}} \dVec{\delta}_{ {2,V}}-2\dVec{\delta}_{ {1,V}}^\top \Tilde {\latMatr{X}} \dVec{\delta}_{ {2,V}}}\\
                &=\dVec{\delta}_{ {1,V}}^\top {\grMatr{\Sigma}} \dVec{\delta}_{ {1,V}} + \dVec{\delta}_{ {2,V}}^\top {\grMatr{\Sigma}} \dVec{\delta}_{ {2,V}}-2\dVec{\delta}_{ {1,V}}^\top {\latMatr{X}} \dVec{\delta}_{ {2,V}} +2k_1-2k_2.\nonumber
            \end{align}
            In this case, 
            \begin{equation*}
                (\ref{eq:appZ_V_contributionOnD}) \to \round{\dVec{\delta}_{ {1,V}}-\dVec{\delta}_{ {2,V}}}^\top \Tilde {\grMatr{\Sigma}} \round{\dVec{\delta}_{ {1,V}}-\dVec{\delta}_{ {2,V}}}
            \end{equation*}
            If we add the hypothesis $\beta\to 1$, $\Corr{Z_{E,i}}{Z_{E,j}}\to 1$ as well. This implies that all the $p\times p$ entries of $\latMatr D(u_1,u_2)$ converge to the same term, that is $d_R(u_1,u_2)$.
            \item Fix $p\geq 2$ and let $\alpha\to +\infty$. The same asymptotic properties (\ref{eq:asymptQ1}) and (\ref{eq:asymptX1}) hold (notice that $p$ and $\alpha$ always appear with the same exponent in each monomial). Therefore, the same reasoning applies.
            \item Fix $p\geq 2$ and let $\alpha\to 0^+.$ From Equation (\ref{eq:defQ}), it is patent that $\latMatr Q \to \frac{1}{2}(\latMatr L + \sqrt{\latMatr L^2})=\latMatr L$; whilst from Equation (\ref{eq:defXEigen}), we have:
            \begin{align*}
                \latMatr X = \frac{1}{2(p-1)} \latMatr W \round{(p-2)\grMatr{\Lambda} + \frac{1}{\alpha}\sqrt{\latMatr I_n - \latMatr J_n - 2\alpha (p-2) \grMatr{\Lambda} + \alpha^2 p^2 \grMatr{\Lambda}} - \frac{1}{\alpha}\latMatr I_n - \frac{1}{\alpha}\latMatr J_n} \latMatr W^\top.
            \end{align*}
            Consider the diagonal matrix in the big brackets: it is possible to expand the square root term by means of the Taylor expansion $\sqrt{a-b\cdot\alpha+c\cdot \alpha^2}=\sqrt{a}-\frac{b}{2\sqrt{a}} \alpha + o(\alpha)$ as $\alpha \to 0^+$: the $i^\text{th}$ element on the main diagonal of $\sqrt{\latMatr I_n - \latMatr J_n - 2\alpha (p-2) \grMatr{\Lambda} + \alpha^2 p^2 \grMatr{\Lambda}}$ is therefore:
            \begin{equation*}
                \begin{dcases*}
                    \sqrt{1} - \frac{2(p-2)\lambda_i}{2\sqrt{1}}\, \alpha + o(\alpha)=1-(p-2)\lambda_i\,\alpha +o(\alpha) \quad &\text{if $i\ne n$}\\
                    0 & \text{if $i=n$.}
                \end{dcases*}    
            \end{equation*}
            As a consequence, we have:
            \begin{align*}
                &(p-2)\grMatr{\Lambda} + \frac{1}{\alpha}\sqrt{\latMatr I_n - \latMatr J_n - 2\alpha (p-2) \grMatr{\Lambda} + \alpha^2 p^2 \grMatr{\Lambda}} - \frac{1}{\alpha}\latMatr I_n - \frac{1}{\alpha}\latMatr J_n \\
                &\quad \sim_a (p-2)\grMatr{\Lambda} + \frac{1}{\alpha} (\latMatr I_n - \latMatr J_n) - (p-2)\grMatr{\Lambda} - \frac{1}{\alpha}\latMatr I_n - \frac{1}{\alpha}\latMatr J_n\\
                &\quad = -\frac{2}{\alpha} \latMatr J_n,
            \end{align*}
            that is:
            \begin{equation*}
                \latMatr X \sim \frac{1}{2(p-1)}\latMatr W \round{-\frac{2}{\alpha}\latMatr J_n} \latMatr W^\top=-\frac{1}{\alpha(p-1)} \frac{1}{n} \onenn,
            \end{equation*}
            where in the last step we have used from the proof of Lemma \ref{lem:usefulIdentities} that $\latMatr W \latMatr J_n \latMatr W^\top = \frac{1}{n} \onenn$ (Equation (\ref{eq:WJWtop_expr})). In this case
            \begin{align*}
                (\ref{eq:appZ_V_contributionOnD}) &\sim_a \dVec{\delta}_{ {1,V}}^\top {\grMatr{\Sigma}} \dVec{\delta}_{ {1,V}} + \dVec{\delta}_{ {2,V}}^\top {\grMatr{\Sigma}} \dVec{\delta}_{ {2,V}}-2\dVec{\delta}_{ {1,V}}^\top \round{-\frac{1}{\alpha n (p-1)}\onenn} \dVec{\delta}_{ {2,V}} -2\frac{1}{\alpha n p (p-1)}\\
                &=\dVec{\delta}_{ {1,V}}^\top {\grMatr{\Sigma}} \dVec{\delta}_{ {1,V}} + \dVec{\delta}_{ {2,V}}^\top {\grMatr{\Sigma}} \dVec{\delta}_{ {2,V}}+\frac{2}{\alpha n p}.
            \end{align*}
            %which explodes as $\alpha \to 0^+$. 
            Clearly, the contribution of the process $\rVec Z_E$ in Equation (\ref{eq:explicitWritingMultivMetric}) and the above term $\dVec{\delta}_{ {1,V}}^\top {\grMatr{\Sigma}} \dVec{\delta}_{ {1,V}} + \dVec{\delta}_{ {2,V}}^\top {\grMatr{\Sigma}} \dVec{\delta}_{ {2,V}}$ do not depend on $\alpha$. As a consequence, $D_{ij}(u_1,u_2)\sim_a \frac{2}{\alpha n p}$ for $i\ne j$ {as $\alpha \to 0^+$}.
            \item {Let $\alpha \to 0^+$ and $p \to +\infty$ such that $\alpha p \tau \to 1$ with $0 < \tau \leq \lambda_{n-1}$ being fixed. Then, from Equations (\ref{eq:defQEigen}) and (\ref{eq:eigL2}), one has $\latMatr{Q} \to \frac{\latMatr{I}_n}{\tau}$. Likewise, from Equation (\ref{eq:defXEigen}), it comes $\latMatr{X} \to \grMatr{\Sigma}-\tau\latMatr{I}_n$.} {The computation of the distances $D_{ij}(u_1,u_2)$ is straightforward: following the same reasoning of the previous derivation, we have 
            \begin{align*}
                (\ref{eq:appZ_V_contributionOnD})
                &=\dVec{\delta}_{ {1,V}}^\top {\grMatr{\Sigma}} \dVec{\delta}_{ {1,V}} + \dVec{\delta}_{ {2,V}}^\top {\grMatr{\Sigma}} \dVec{\delta}_{ {2,V}}-2\dVec{\delta}_{ {1,V}}^\top \round{\grMatr{\Sigma} - \tau \latMatr I_n} \dVec{\delta}_{ {2,V}} -2\frac{1}{\alpha n p (p-1)}\\
                &\to \round{\dVec \delta_{1,V}-\dVec \delta_{2,V}}^\top \grMatr{\Sigma} \round{\dVec \delta_{1,V}-\dVec \delta_{2,V}} + 2\tau \dVec \delta_{1,V}^\top \dVec \delta_{2,V}.
            \end{align*}
            The first addendum, together with the contribution of the process $\rVec Z_E$ (if $\beta \to 1$), gives the distance $d_R(u_1,u_2)$.
            } 
            \item {Let $\alpha \to 0^+$ and $p\to + \infty $, with $\alpha p \tau \to 1$ with $\tau \geq \lambda_1$ fixed. Equation (\ref{eq:defQEigen}) gives
            \begin{align*}
                \latMatr Q \to \latMatr W \diag{[\lambda_1^{-1},\ldots,\lambda_{n-1}^{-1},\frac{1}{\tau}]^\top} \latMatr W^\top = \latMatr W \round{\grMatr{\Lambda}^- + \frac{1}{\tau} \latMatr J_n} \latMatr W^\top = \latMatr L + \frac{1}{n\tau}\onenn.
                %\latMatr Q \to \latMatr W \begin{bmatrix}
                %    \lambda_1^{-1} &&&\\
                %    & \ddots & & \\
                %    & & \lambda_{n-1}^{-1}\\
                %    & & & \tau
                %\end{bmatrix} \latMatr W^\top = \latMatr W \round{\grMatr{\Lambda}^- + \tau \latMatr J_n} \latMatr W^\top = \latMatr L + \frac{\tau}{n}\onenn.
            \end{align*}
            Therefore, from Equation (\ref{eq:defX}):
            \begin{align*}
                \latMatr X = \frac{p\tau}{p-1} \round{\grMatr{\Sigma} \round{\latMatr L + \frac{1}{n\tau}\onenn.}-\latMatr I_n} \to -\frac{\tau}{n} \onenn, \qquad \text{as $p\to +\infty$.} 
            \end{align*}
            Hence,
            \begin{align*}
                (\ref{eq:appZ_V_contributionOnD})&\to \dVec{\delta}_{ {1,V}}^\top {\grMatr{\Sigma}} \dVec{\delta}_{ {1,V}} + \dVec{\delta}_{ {2,V}}^\top {\grMatr{\Sigma}} \dVec{\delta}_{ {2,V}}-2\dVec{\delta}_{ {1,V}}^\top \round{-\frac{\tau}{n}\onenn} \dVec{\delta}_{ {2,V}} -2\frac{1}{\alpha n p (p-1)}\\
                &\to \dVec{\delta}_{ {1,V}}^\top {\grMatr{\Sigma}} \dVec{\delta}_{ {1,V}} + \dVec{\delta}_{ {2,V}}^\top {\grMatr{\Sigma}} \dVec{\delta}_{ {2,V}}+\frac{2\tau}{n}.
            \end{align*}
            The result for the distance $D_{ij}(u_1,u_2)$ is obtained from the last equation by adding and subtracting $2\dVec{\delta}_{1,V}^\top \latMatr X \dVec{\delta}_{2,V}$ and by using the fact that $\beta\to 1$.
            } 
        \end{enumerate}
    \end{proof}
    
    \begin{proof}[Proof of Theorem \ref{theo:covFunctionBuilding}]
        %\xave{(Xavier: OK with the proof)} 
        The result stems from Theorem \ref{theo:MultivariateSchoenberg} and the fact that the pseudo-variogram ${\latMatr{D}}$ satisfies ${\latMatr{D}}^\top(u_1,u_2)={\latMatr{D}}(u_2,u_1)$ and is $\onep$-conditionally negative semidefinite.%To prove the result it is sufficient to apply Theorem \ref{theo:MultivariateSchoenberg} to our case. Indeed, the pseudo-variogram $d=\Gamma$ satisfies $\Gamma^\top(u_1,u_2)=\Gamma(u_2,u_1)$ and it is $\onep$-conditionally negative semidefinite.
    \end{proof}

    \begin{proof}[Proof of Theorem \ref{theo:la_raja}]   
    {The result stems from Theorem \ref{theo:covFunctionBuilding} with $\psi(t)=\exp(-t)$ and {$g(t)=t^\theta$} and the fact that the set of positive semidefinite matrix-valued functions is closed under sums, point-wise limits and element-wise products with symmetric positive semidefinite matrices.}
    \end{proof}

    \begin{proof}[Proof of Theorem \ref{theo:dios}]
        {From \citet[Equations (28)-(30) and Theorem 3]{EPW}, $\boldsymbol{{\cal M}}_p$ has a representation of the form (\ref{eq:mapping}) with ${\rm d}F_{ij}(\xi)=R_{ij}(\xi) {\rm d}\xi$, where $[R_{ij}(\xi)]_{i,j=1}^p$ is positive semidefinite under the stated conditions (A) or (B). The theorem thus results from the application of Theorem \ref{theo:la_raja}.}
    \end{proof}
    
    {\begin{proof}[Proof of Theorem \ref{theo:dios2}]
    Using formula 3.381.4 of \cite{grad}, one can write 
    \begin{equation*}
        \boldsymbol{{\cal C}}_p(\latMatr{D}) = \Bigg[ \frac{\sigma_{ij} \beta_{ij}^\nu}{\Gamma(\nu)} \int_0^{+\infty} \xi^{\nu-1} \exp(-\beta_{ij} \xi) \exp(-\xi D_{ij}^\theta) {\rm d}\xi \Bigg]_{i,j=1}^p.
    \end{equation*}
    The conditional negative semidefiniteness of $[\beta_{ij}]_{i,j=1}^p$ implies that $\big[\frac{\xi^{\nu-1}}{\Gamma(\nu)} \exp(-\beta_{ij} \xi)\big]_{i,j=1}^p$ is positive semidefinite for any $\xi \geq 0$ \citep[Lemma 2.5]{Reams}. The theorem results from the application of the Schur product theorem and Theorem \ref{theo:la_raja}.
    \end{proof}}

    \begin{proof}[Proof of Theorem \ref{theo:dios3}]
        One can write \citep{porcu2018shkarofsky}
        \begin{equation*}
            \boldsymbol{{\cal SG}}_p(\latMatr{D}) = \Bigg[ \frac{\sigma_{ij} 2^{\nu-1} (\alpha_{ij} \beta_{ij})^{\nu/2}}{{\cal K}_\nu\left(\sqrt{\frac{\beta_{ij}}{\alpha_{ij}}}\right)} \int_0^{+\infty} \xi^{\nu-1} \exp\left(- \frac{1}{4\alpha_{ij} \xi}-\beta_{ij} \xi \right)  \exp(-\xi D_{ij}^\theta) {\rm d}\xi \Bigg]_{i,j=1}^p.
        \end{equation*}
        The fact that $[\alpha^{-1}_{ij}]_{i,j=1}^p$ and $[\beta_{ij}]_{i,j=1}^p$ are conditionally negative semidefinite implies that $$\big[\exp\left(- \frac{1}{4\alpha_{ij} \xi}-\beta_{ij} \xi \right)\big]_{i,j=1}^p$$ is positive semidefinite for any $\xi \geq 0$. The theorem results from the application of the Schur product theorem and Theorem \ref{theo:la_raja}.
    \end{proof}

    \begin{proof}[Proof of Proposition \ref{prop:p-Schoen}] %\xave{OK FOR THE FIRST PART; THE SECOND PART IS UNCLEAR TO ME...} \emi{can you be more specific?} \xave{I fixed the end of the proof: OK now for me}
    We start with the case $p=1$ and observe that, by Bochner's theorem \citep{bochner1955harmonic}, a continuous mapping $C: \R^{{m^\prime}} \to \R$ is positive semidefinite if and only if 
    $$ C(\dVec{x}) = \int_{\R^{{m^\prime}}} \exp(\mathsf{i} \dVec{x}^{\top} \dVec{\omega}) G({\rm d}\dVec{\omega}), \qquad \dVec{x} \in \R^{{m^\prime}}, $$
    for $G$ being a non-decreasing and bounded measure in $\R^m$. Here, $\mathsf{i}$ denotes the unit complex number. Using the above representation, \cite{cambanis1983symmetric} show that a mapping $K: [0,+\infty) \to \R$ belongs to the class $\Phi_m^1$ if and only if Equation (\ref{eq:alpha-1}) holds for $p=
    1$ and for $F$ being non-decreasing and bounded. {In particular, choosing a Dirac measure for $F$, one finds that $t \mapsto \omega_{m^\prime}(rt)$ belongs to $\Phi_m^1$ for any $r \geq 0$, that is, $\omega_{m^\prime}(r d_M(\cdot,\cdot))$ is positive semidefinite on $\R^{m^\prime}$.}\par
    We now address the general case with $p \geq 1$. First, consider a matrix-valued mapping $\latMatr{K}$ belonging to the class $\Phi_m^p$. We invoke the extension of Bochner's theorem to the $p$-variate case \citep{cramer1940theory}, for which a continuous mapping $\latMatr{C}: \R^{{m^\prime}} \to \R^{p \times p}$ is positive semidefinite if and only if 
    $$ \latMatr{C}(\dVec{x}) = \int_{\R^{{m^\prime}}} \exp(\mathsf{i} \dVec{x}^{\top} \dVec{\omega})  \latMatr{G}({\rm d}\dVec{\omega}), \qquad \dVec{x} \in \R^{{m^\prime}}, $$
    for $\latMatr{G}$ a bounded {complex Hermitian} matrix-valued measure. Hence, the same calculations as much as in \cite{cambanis1983symmetric} show that the representation (\ref{eq:alpha-1}) holds for a bounded {symmetric} real matrix-valued measure $\latMatr{F}$ {(that $\latMatr{F}$ is bounded stems from the fact that $\latMatr{K}(0)$ is finite)}. To finish this part of the proof we need to show that ${\rm d}\latMatr{F}(r)$ is positive semidefinite for any $r \geq 0$, which happens if and only if the mapping
    $$ H_{{\dVec{c}}}(r) := \dVec{c}^\top {\latMatr{F}}({\rm d}r) \dVec{c}, \quad r \in [0,+\infty),$$
    is non-negative for any arbitrary $\dVec{c} \in \R^p$. 
    %Now, let $\latMatr{K}$ belong to the class $\Phi_m^p$. 
    To this end, consider a $p$-variate random field $\rVec{Z}$ in $\R^{{m^\prime}}$ with isotropic covariance $\latMatr{K}$, and define $Z_{\dVec{c}} := \dVec{c}^\top \dVec{Z}$, which is a univariate random field with isotropic covariance $K_{\dVec{c}} = \dVec{c}^\top \latMatr{K} \dVec{c}$. Hence, $K_{\dVec{c}}$ obeys to the representation provided by \cite{cambanis1983symmetric}, with a measure $F_{\dVec{c}}$ such that 
    $$ {F_{\dVec{c}} = \dVec{c}^\top \latMatr{F} \dVec{c}}. $$
    This proves the sufficiency part of the theorem by noting that $F_{\dVec{c}}$ is non-decreasing. \par
    {Reciprocally, assume that the representation (\ref{eq:alpha-1}) holds with a bounded symmetric matrix-valued measure $\latMatr{F}$ such that ${\rm d}\latMatr{F}(r)$ is positive semidefinite for any $r \geq 0$. Being positive semidefinite on $\R^{m^\prime}$ for any $r \geq 0$, $\omega_{m^\prime}(r d_M(\cdot,\cdot))$ is the covariance function of a univariate Gaussian random field $Z_r$ on $\R^{m^\prime}$. Accordingly, the matrix-valued mapping $\omega_{m^\prime}(r d_M(\cdot,\cdot)) \onepp$ belongs to $\Phi_m^p$ as the covariance function of $Z_r \onep$, and so do the mappings $\omega_{m^\prime}(r d_M(\cdot,\cdot)) {\rm d}\latMatr{F}(r)$ and $\latMatr{K}(d_M(\cdot,\cdot))$ owing to the fact that the set of positive semidefinite functions is closed under sums, element-wise products, and pointwise limits.}
    \end{proof}
    
    \begin{proof}[Proof of Proposition \ref{prop:PleaseEmilioChooseWiseLabels}]
    We provide a proof by direct construction. {Let ${m^\prime}$ be a positive integer.} Arguments in Theorem 2.1 of \cite{gneiting1998onalpha} provide the identity 
    \begin{equation}
     \label{omega2Omega}
        \omega_{{m^\prime}} (t) = \kappa_{{m^\prime}} I^{({{m^\prime}}-1)} \Omega_{2{{m^\prime}} -1} (t), \qquad t \ge 0,
    \end{equation}
    for $\kappa_{{m^\prime}}$ a strictly positive constant. Hence, the proof comes from the integral representation (\ref{eq:alpha-1}) %of the function $\latMatr{H}$ according to the $p$-valued version of Schoenberg's theorem \citep{AlonsoMalaver2015251} in concert with Proposition \ref{prop:p-Schoen} and 
    in concert with a straight use of Fubini's theorem. 
    \end{proof}
    
    \begin{proof}[Proof of Proposition \ref{prop:infarto_miocardio}]
        The results comes from a straight application of Proposition 2 in \cite{emery2023schoenberg} in concert with Lemma \ref{lem:tree}, Proposition \ref{prop:PleaseEmilioChooseWiseLabels} and Fubini's theorem. 
    \end{proof}    
\newpage
\section{A Worked Example} \label{A_B}

    Here we provide a simple example of the multivariate distance properties.
    \begin{example}
        \label{ex:counterExampleDistance}
        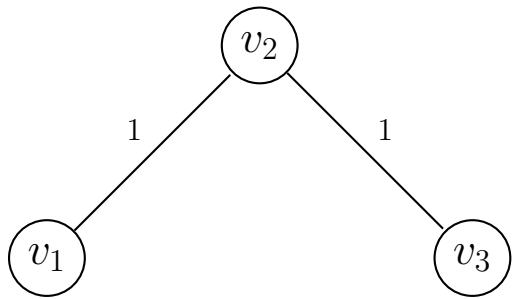
\begin{figure}
            \centering
            \begin{tikzpicture}[shorten >=1pt,auto,node distance=4cm, thick, main node/.style={circle, draw, font=\Large\bfseries}]
                
                \node[main node] (2) {$v_2$};
                \node[main node] (1) [below left of=2] {$v_1$};
                \node[main node] (3) [below right of=2] {$v_3$};
                
                \path
                (1) edge node [above left]{1} (2)
                (2) edge node [above right] {1} (3);      
            \end{tikzpicture}
            \caption{An example of graph for which the distance $D_{11}(v_1,v_3) \not \leq D_{12}(v_1,v_3)$, for $p\geq 3$ and $\alpha$ wisely chosen.}
            \label{fig:counterExampleGraph}
        \end{figure}
        Let $p\geq 3$ and $\alpha>0$. Consider the graph depicted in Figure \ref{fig:counterExampleGraph} with $V:=\curly{v_1,v_2,v_3}$ and whose Laplacian matrix $\latMatr{L}$ and $\grMatr{\Sigma}={\latMatr{L}}^-$ are 
        \begin{equation*}
            \latMatr{L}=\begin{bmatrix}
                1 & -1 & 0\\
                -1 & 2 & -1\\
                0 & -1 & 1
            \end{bmatrix},
            \qquad
            {\grMatr{\Sigma}}=\frac{1}{9}\begin{bmatrix}
                5 & -1 & -4\\
                -1 & 2 & -1\\
                -4 & -1 & 5
            \end{bmatrix}.
        \end{equation*}
        In addition, the eigendecomposition ${\grMatr{\Sigma}}={\latMatr{W}} {\grMatr{\Lambda}} {\latMatr{W}}^\top$ is
        \begin{equation*}
            {\latMatr{W}}=\begin{bmatrix}
                -\frac{1}{\sqrt{2}} & \frac{1}{\sqrt{6}} & \frac{1}{\sqrt{3}} \\
                0 & -\frac{2}{\sqrt{6}} & \frac{1}{\sqrt{3}}\\
                \frac{1}{\sqrt{2}} & \frac{1}{\sqrt{6}} & \frac{1}{\sqrt{3}}
            \end{bmatrix},
            \qquad
            {\grMatr{\Lambda}} = \begin{bmatrix}
                1 & & \\
                & \frac{1}{3} & \\
                & & 0
            \end{bmatrix}.
        \end{equation*}
        Let us consider the distance between $v_1$ and $v_3$ for the same variable and for different variables, \ie $D_{11}(v_1,v_3)$ and $D_{12}(v_1,v_3)$. Since both $v_1$ and $v_3$ are vertices, the contribution of $\boldsymbol{Z}_E$ is null, therefore, by Equation (\ref{eq:explicitWritingMultivMetric}), we have:
        \begin{align*}
            D_{11}(v_1,v_3)&=\round{\boldsymbol{\delta}_{ {1,V}}-\boldsymbol{\delta}_{ {3,V}}}^\top \Tilde {\grMatr{\Sigma}} \round{\boldsymbol{\delta}_{ {1,V}}-\boldsymbol{\delta}_{ {3,V}}} \\&= \square{1,0,-1}\Tilde {\grMatr{\Sigma}} \begin{bmatrix}
                1\\
                0\\
                -1
            \end{bmatrix} \\
            &=\Tilde{{{\Sigma}}}_{11}+\Tilde{{{\Sigma}}}_{33}-2\Tilde{{{\Sigma}}}_{13}\\
            &={{\Sigma}}_{11}+{{\Sigma}}_{33}-2{{\Sigma}}_{13}=2,\\
            D_{12}(v_1,v_3)&=\boldsymbol{\delta}_{ {1,V}}^\top \Tilde {\grMatr{\Sigma}} \boldsymbol{\delta}_{ {1,V}} + \boldsymbol{\delta}_{ {3,V}}^\top \Tilde {\grMatr{\Sigma}} \boldsymbol{\delta}_{ {3,V}}-2\boldsymbol{\delta}_{ {1,V}}^\top \Tilde {\latMatr{X}} \boldsymbol{\delta}_{ {3,V}}\\
            &=\Tilde{{\Sigma}}_{11}+\Tilde{{\Sigma}}_{33}-2\Tilde {{X}}_{13}\\
            &={{\Sigma}}_{11}+{{\Sigma}}_{33}+2k_1-2 {{X}}_{13} -2k_2\\
            &=\frac{10}{9}-2\cdot \frac{1}{3\alpha p (p-1)}-2{{X}}_{13}.
        \end{align*}
        Now it is possible to compute ${{X}}_{13}$ by means of Equation (\ref{eq:defX}): %\tobi{please adjust the notation}:
        \begin{align*}
            \allowdisplaybreaks
            &{{X}}_{13}=\frac{1}{2\alpha (p-1)} {{W}}_{1\star} \round{\alpha(p-2){\grMatr{\Lambda}} + \sqrt{ {\latMatr{I}_n - \latMatr{J}_n} - 2 \alpha (p-2){\grMatr{\Lambda}} + \alpha^2 p^2 {\grMatr{\Lambda}}^2} -{\latMatr{I}_n}-\latMatr{J}_n} {{W}}^\top_{\star3}\\
            &=\frac{1}{2\alpha (p-1)} \sum_{i=1}^3 {{W}}_{1i}{{W}}_{3i} \round{\alpha(p-2){\grMatr{\Lambda}} + \sqrt{ {\latMatr{I}_n - \latMatr{J}_n} - 2 \alpha (p-2){\grMatr{\Lambda}} + \alpha^2 p^2 {\grMatr{\Lambda}}^2} -{\latMatr{I}_n}-\latMatr{J}_n}_{ii}\\
            &=\frac{1}{2\alpha (p-1)}\Bigg(
            -\frac{1}{2}\round{\alpha (p-2)\cdot 1+\sqrt{1-2\alpha(p-2)\cdot 1 +\alpha^2 p^2 \cdot 1^2}-1}\\
            & \qquad + \frac{1}{6} \round{\alpha (p-2)\cdot \frac{1}{3}+\sqrt{1-2\alpha(p-2)\cdot \frac{1}{3} +\alpha^2 p^2 \cdot \frac{1}{3^2}}-1} \\
            & \qquad + \frac{1}{3} \round{\alpha (p-2) \cdot 0 + \sqrt{0-2\alpha (p-2)\cdot 0 + \alpha^2 p^2 \cdot 0^2}-2}
            \Bigg)\\
            &=\frac{1}{2\alpha (p-1)} \Bigg(-\frac{1}{2} \sqrt{1-2\alpha(p-2)+\alpha^2 p^2}\\
            &\qquad +\frac{1}{6}\sqrt{1-\frac{2}{3}\alpha(p-2)+\frac{\alpha^2 p^2}{9}} -\frac{4}{9} \alpha (p-2) -\frac{1}{3}
            \Bigg).
        \end{align*}
        {
        With $p=3$, the maximum discrepancy $d(\alpha):=D_{11}(v_1,v_3)-D_{12}(v_1,v_3)$ is reached at $\alpha^\star=1.32092$. In such a case, ${X}_{13}=-0.48598$ and thus:
        \begin{equation*}
            D_{12}(v_1,v_3)=\frac{10}{9} - 2\cdot \frac{1}{18\cdot 1.32092} - 2\cdot (-0.48598)=1.99896<2=D_{11}(v_1,v_3),
        \end{equation*}
        that is, $d(\alpha^\star)=0.00103.$ Notice that for $\alpha:=\hat \alpha=0.97222$ and for $\alpha\to+\infty$, we achieve %the equality %$D_{12}(v_1,v_3)=D_{11}(v_1,v_3)$ (\ie 
        $d(\hat\alpha)=0$, and for $0<\alpha<\hat\alpha$ we have $d(\alpha)<0$.  %$D_{12}(v_1,v_3)>D_{11}(v_1,v_3)$. \\
        %For the above choice, the discrepancy $d(\alpha^\star)$ is very small. However, it is possible to achieve more significant values of 
        One can obtain larger values for $d(\alpha)$ by taking $\alpha:=\frac{1}{p}$ and letting $p\to +\infty$. In this setting, we have ${X}_{13}\to -\frac{1}{3}$ and $d(\alpha)\to \frac{2}{9}$}.% Thus,
        %\begin{align*}
            %D_{12}(v_1,v_3)&\to \frac{10}{9} + \frac{2}%{3}=\frac{16}{9}<2=D_{11}(v_1,v_3),
        %\end{align*}
        %that is $d(\alpha)\to \frac{2}{9}$.}
    \end{example}

\end{document}